\documentclass[a4paper,10pt]{scrartcl}
\pdfoutput=1
\usepackage{ifpdf}
\usepackage[english]{babel}
\usepackage[utf8]{inputenc}
\usepackage[T1]{fontenc}
\usepackage{lmodern}
\usepackage{geometry}
\usepackage{amsfonts,amsmath,amssymb,amsopn}
\usepackage{amsthm}
\usepackage{mathtools}
\usepackage{nicefrac}
\usepackage{leftindex}
\usepackage{exscale}
\usepackage{empheq}
\usepackage{accents}

\usepackage[numbers,square]{natbib}
\usepackage{url} 
\usepackage[colorlinks=true,pdfpagelabels,unicode]{hyperref}
\hypersetup{linkcolor=blue, urlcolor=blue, citecolor=blue}
\usepackage[dvipsnames,svgnames,table]{xcolor}
\usepackage{orcidlink}

\allowdisplaybreaks 

\theoremstyle{plain}

\newtheoremstyle{theo}
	{3pt} 
	{3pt} 
	{\itshape} 
	{} 
		{\bfseries} 
	{\\} 
	{ } 
	{\thmname{#1}\thmnumber{ #2.}\thmnote{ - #3}} 
\theoremstyle{theo}

\newtheorem{definition}{Definition}[section]
\newtheorem{lemma}[definition]{Lemma}
\newtheorem{theorem}[definition]{Theorem}
\newtheorem{corollary}[definition]{Corollary}

\newenvironment{bew}{\begin{proof}[\bfseries Proof:]}{\end{proof}}

\newtheoremstyle{remark}
	{3pt} 
	{3pt} 
	{} 
	{} 
		{\bfseries} 
	{} 
	{ } 
	{\thmname{#1}\thmnumber{ #2.}\thmnote{ - #3}} 
\theoremstyle{remark}
\newtheorem{remark}[definition]{Remark}

\DeclareMathOperator{\bomega}{\overline{\Omega}}
\DeclareMathOperator{\romega}{\partial\Omega}

\DeclareMathOperator{\supp}{supp}
\DeclareMathOperator{\intd}{d\!}

\DeclareMathOperator{\dist}{dist}

\newcommand{\rd}{r}
\newcommand{\sig}{\varpi}

\newcommand{\GNI}{Gagliardo--Nirenberg inequality}

\newcommand{\into}[1]{\int_0^{#1}\!}

\newcommand{\intoT}{\into{T}}

\newcommand{\intomega}{\int_{\Omega}\!} 
\newcommand{\intoTomega}{\intoT\!\intomega}

\newcommand{\Lo}[1][1]{L^{#1}(\Omega)} 
\newcommand{\W}[1][1,2]{W^{#1}(\Omega)}

\newcommand{\LSp}[2]{L^{#1\;\!}\!\left(#2\right)} 
\newcommand{\LSpn}[2]{L^{#1\;\!}\!(#2)}
\newcommand{\LSpb}[2]{L^{#1\;\!}\!\big(#2\big)}
\newcommand{\LSploc}[2]{L_{loc}^{#1}\!\left(#2\right)} 

\newcommand{\LSplocb}[2]{L_{loc}^{#1}\big(#2\big)}

\newcommand{\CSp}[2]{C^{#1}\!\left(#2\right)}

\newcommand{\CSploc}[2]{C_{loc}^{#1}\!\left(#2\right)}

\newcommand{\R}{\mathbb{R}}
\newcommand{\N}{\mathbb{N}}

\newcommand{\nfrac}[2]{{\nicefrac{#1}{#2}}}

\newcommand{\dimN}{N}

\DeclareMathOperator*{\essosc}{ess\,osc}
\DeclareMathOperator*{\essinf}{ess\,inf}
\newcommand{\intQ}[1][\rd]{\iint\limits_{Q_{#1}}}
\newcommand{\OmT}{\Omega_T}
\newcommand{\bOmT}{\overline{\Omega_T}}
\newcommand{\intQO}[1][\rd]{\iint\limits_{Q_{#1}\cap\OmT}}
\DeclareMathOperator*{\esssup}{ess\, sup}
\DeclareMathOperator*{\esslim}{ess\, lim}

\newcommand{\dels}{\delta_\star}
\newcommand{\deld}{\delta_\diamond}
\newcommand{\cx}{\check{x}}

\newcommand{\js}{j_{\star}}
\newcommand{\mK}{\mathcal{K}}

\makeatletter
\def\@fnsymbol#1{\ensuremath{\ifcase#1\or *\or \ddagger\or \#\or
   \mathsection\or \mathparagraph\or \|\or **\or \dagger\dagger
   \or \ddagger\ddagger \else\@ctrerr\fi}}
\makeatletter
\def\@fnsymbol#1{\ensuremath{\ifcase#1\or *\or \ddagger\or \#\or
   \mathsection\or \mathparagraph\or \|\or **\or \dagger\dagger
   \or \ddagger\ddagger \else\@ctrerr\fi}}
   
 \def\@lbibitem[#1]#2#3{%
  \if\relax\@extra@b@citeb\relax\else
    \@ifundefined{br@#2\@extra@b@citeb}{}{%
     \@namedef{br@#2}{\@nameuse{br@#2\@extra@b@citeb}}%
    }%
  \fi
  \@ifundefined{b@#2\@extra@b@citeb}{%
   \def\NAT@num{}%
  }{%
   \NAT@parse{#2}%
  }%
  \def\NAT@tmp{#1}%
  \expandafter\let\expandafter\bibitemOpen\csname NAT@b@open@#2\endcsname
  \expandafter\let\expandafter\bibitemShut\csname NAT@b@shut@#2\endcsname
  \@ifnum{\NAT@merge>\@ne}{%
   \NAT@bibitem@first@sw{%
    \@firstoftwo
   }{%
    \@ifundefined{NAT@b*@#2}{%
     \@firstoftwo
    }{%
     \expandafter\def\expandafter\NAT@num\expandafter{\the\c@NAT@ctr}%
     \@secondoftwo
    }%
   }%
  }{%
   \@firstoftwo
  }%
  {%
   \global\advance\c@NAT@ctr\@ne
   \@ifx{\NAT@tmp\@empty}{\@firstoftwo}{%
    \@secondoftwo
   }%
   {%
    \expandafter\def\expandafter\NAT@num\expandafter{\the\c@NAT@ctr}%
    \global\NAT@stdbsttrue
   }{}%
   \bibitem@fin
   \item[\href{#3}{\hfil\NAT@anchor{#2}{\NAT@num}}]
   \global\let\NAT@bibitem@first@sw\@secondoftwo
   \NAT@bibitem@init
  }%
  {%
   \NAT@anchor{#2}{}%
   \NAT@bibitem@cont
   \bibitem@fin
  }%
  \@ifx{\NAT@tmp\@empty}{%
    \NAT@wrout{\the\c@NAT@ctr}{}{}{}{#2}%
  }{%
    \expandafter\NAT@ifcmd\NAT@tmp(@)(@)\@nil{#2}%
  }%
}
\makeatother

\author{
Tobias Black\footnote{tblack@math.upb.de}\ \orcidlink{0000-0001-9963-0800}\\
{\small Institute of Mathematics, Paderborn University,}\\[-5pt]
{\small 33098 Paderborn, Germany}
}
\title{Refining H\"older regularity theory in degenerate drift-diffusion equations}
\date{}

\begin{document}
\maketitle
\begin{abstract}
\noindent
{\textbf{Abstract.} We establish the Hölder continuity of bounded nonnegative weak solutions to
\begin{align*}
\big(\Phi^{-1}(w)\big)_t=\Delta w+\nabla\cdot\big(a(x,t)\Phi^{-1}(w)\big)+b\big(x,t,\Phi^{-1}(w)\big),
\end{align*}
with convex $\Phi\in C^0([0,\infty))\cap C^2((0,\infty))$ satisfying $\Phi(0)=0$, $\Phi'>0$ on $(0,\infty)$ and $$s\Phi''(s)\leq C\Phi'(s)\quad\text{for all }s\in[0,s_0]$$
for some $C>0$ and $s_0\in(0,1]$. The functions $a$ and $b$ are only assumed to satisfy integrability conditions of the form 
\begin{align*}
a&\in L^{2q_1}\big((0,T);L^{2q_2}(\Omega;\mathbb{R}^N)\big),\\
b&\in M\big(\Omega_T\times\mathbb{R}\big)\ \text{such that }\big|b(x,t,\xi)\big|\leq \hat{b}(x,t)\ \text{a.e. for some }\hat{b}\in L^{q_1}\big((0,T);L^{q_2}(\Omega)\big)
\end{align*}
with $q_1,q_2>1$ such that $$\frac{2}{q_1}+\frac{N}{q_2}=2-N\kappa\quad\text{for some }\kappa\in(0,\tfrac{2}{N}).$$
Letting $w=\Phi(u)$ and assuming the inverse $\Phi^{-1}:[0,\infty)\to[0,\infty)$ to be locally Hölder continuous, this entails Hölder regularity for bounded weak solutions of $$u_t=\Delta\Phi(u)+\nabla\cdot\big(a(x,t)u\big)+b(x,t,u)$$ and, accordingly, covers a wide array of taxis type structures. In particular, many chemotaxis frameworks with nonlinear diffusion, which cannot be covered by the standard literature, fall into this category.

After rigorously treating local Hölder regularity, we also extend the regularity result to the associated initial-boundary value problem for boundary conditions of flux-type.
}\medskip



\end{abstract}

\newpage
\section{Introduction}\label{sec1:intro}
Plenty of microorganisms, although seeming to move in an unorganized and random fashion individually, tend to exhibit distinct taxis-schemes on population scales by reacting to external environmental stimuli as a group. One archetypical process, where the movement is provoked by concentration gradients of a signal chemical dissolved in the surrounding habitat, is often referred to as chemotaxis. Experimental observations, e.g.\ on specimens of \emph{Escherichia coli} (\cite{Adler66}) and \emph{Bacillus subtilis} (\cite{Dombrowski04,tuval2005bacterial}), indicate that the ability to adapt motion in response to perceived signal may lead to precedents with peculiar dynamics and a profound interest was sparked in the community of mathematical biology, when Keller and Segel proposed the system
\begin{align*}
\begin{array}{r@{\ }l}
u_t&=\Delta u-\nabla\cdot\big(u\chi(v)\nabla v\big),\\
v_t&=\Delta v+f(u,v),
\end{array}
\end{align*}
to describe the experimentally observed aggregation process of the slime mold \emph{Dictyostelium discoideum} (\cite{KS70}). Herein, $u$ and $v$ denote the density of the cells and the concentration of a chemoattractant, respectively. Ever since its introduction, colorful new mathematical facets beyond the enigmatic aggregation mechanism have been discovered in the Keller--Segel system and its related variants. (We refer to the surveys \cite{Ho03,BBWT15} for an overview on both historical and analytical context.)\smallskip

PDE-systems within the framework of chemotaxis equations are abundant and modifications range from assuming a movement mechanism governed by nonlinear diffusion motivated by enhanced cell stress in densely packed areas (\cite{HP-volumefilling-CAMQ02,TaoWin-quasilinear_JDE12}), over including logistic source terms to incorporate cell reproduction and death (\cite{Os02-chemologatract,TW07}) to coupling with additional equations like the (Navier--)Stokes-equations to allow feedback from the surrounding environment to the cell dynamics (\cite{Dombrowski04,lorz-M3AS10}). The common building block of all relatives, however, is interaction between cells and signal by means of a cross-diffusive taxis term like $-\nabla\cdot\big(u\chi(v)\nabla v\big)$. In many cases, depending on the context and the precise choice for $\chi$, even the global existence of a solution is a priori uncertain. Moreover, even the expected class of solutions to chemotaxis systems is far from clear-cut. Solutions to slightly different systems feature quite diverse regularity properties. In some instances classical solutions may be obtained (\cite{TW07,Winkler_MMAS11}), while in others already the moderate requirements for the standard concept of weak solutions can not be verified with current methods and only (less regular) generalized solutions are available (\cite{win15_chemorot,Win16CS1,LanWin2017}). (See also the survey \cite{LanWin_JDMV_19} for further discussion on low-regularity solutions in chemotaxis settings.)\smallskip

In those settings, unfortunately, some qualitative information on the solution behavior may remain undetected. Take for instance the recent study in \cite{xuChemotaxisModelDegenerate2020}, where Hölder regularity was the key to exemplify interesting properties like the initial shrinking of the support and also its eventual expansion for a system featuring a diffusion of porous medium type, which contrasts the standard nonlinear heat equation (\cite{Vazquez-PME-07}) without taxis. Interest in regularity questions, however, is far from a singular occurrence in chemotaxis contexts alone. Throughout the whole field of PDE the tools to (easily) decide whether the solution to a certain PDE can in fact be Hölder continuous (or better) have been highly sought after and the literature is accordingly rich on previous results, see e.g. \cite{LSU,lieberman,dibenedettoHolderEstimatesNonlinear1985,dibenedettoDegenerateParabolicEquations1993,
PorzVesp93,gianazzaNewProofHolder2010,liaoUnifiedApproachHolder2020,hissinkmullerInteriorHolderContinuity2022} and references therein. Since these results, however, mostly work with parabolic equations of a quite general form and try to cover a broad variety of structures, it would not be too surprising if these results over-restrict in some situations. Especially when considering how intricate the influence of the taxis term in chemotaxis frameworks on the regularity of solutions are identified to be.\smallskip

\textbf{Coping with drift terms which are large relative to diffusion.} The underlying motivation in our investigation is the possibility of restrained parameter choices appearing in available results when degenerate diffusion is in interplay with chemotaxis-type advection terms. Let us, for instance, take a more detailed look at the celebrated and widely utilized Hölder regularity results of Porzio and Vespri (\cite{PorzVesp93}) for parabolic equations with principal part in divergence form, i.e.
$$u_t-\nabla\cdot a(x,t,u,\nabla u)-b(x,t,u,\nabla u)=0.$$
For their result, $a,b$ have to satisfy certain structure conditions. To make our motivational point, we will only have to highlight one of these structural conditions, which can be stated as $$a(x,t,u,\nabla u)\cdot\nabla u\geq C\phi(|u|)|\nabla u|^p-\psi_0(x,t)$$ with $p\geq2$, some continuous function $\phi$ satisfying a power-type growth restriction and a function $\psi_0$ lying in a suitable $L^{q}L^r$-space. Now, in the context of a chemotaxis system with diffusion of porous medium type, which typically would amounts to setting $a(x,t,u,\nabla u)=u^{m-1}\nabla u-u\nabla v$ with $m\geq 1$, even under the strong assumption that perchance we have arbitrarily good control over both $u$ and $\nabla v$, we could only expect $a(x,t,u,\nabla u)$ to satisfy an inequality of the form $$a(x,t,u,\nabla u)\cdot \nabla u\geq \frac{1}{2}u^{m-1}|\nabla u|^2-\frac{1}{2}u^{3-m}|\nabla v|^2,$$ consequently restricting the range for permissible $m$ to $[1,3]$ unless additional restrictions on positivity properties of $u$ are imposed, which, however, would inconveniently narrow the set of admissible data. (See the results on nontrivial support evolution in \cite{fischerAdvectionDrivenSupportShrinking2013} and \cite{stevensTaxisdrivenPersistentLocalization2022}.) From intuition this restriction seems inappropriate for the coupled chemotaxis system, as precedents in chemotaxis-fluid frameworks indicate that a higher diffusion exponent generally permits the solution to exhibit enhanced regularity properties, so that, essentially, stronger diffusion grants more leeway in withstanding the destabilizing effects of the cross-diffusive term (\cite{Win-ct_fluid_3d-CPDE15,WangLi-ZAMP17,TB2017_nonlindiff,TB2018_NA}). This crude observation suggests, that the formidable results of \cite{PorzVesp93} for the general divergence form still have some blind spots when it comes to the convoluted interplay between degenerate diffusion and convection terms, and chemotaxis systems provide one archetypical example of the fact.\smallskip 

Unsuspecting, one might argue that straightforward particularization of the equation towards a more closely fitted chemotaxis framework should address this problem, and that the powerful methods for proving Hölder regularity (see e.g. \cite[Chapter 1]{BenedettoGianazzaVespri-Harnack} for some historical background) established throughout the last decades would then also immediately allow more freedom in the choice of $m$. Unfortunately, this does not appear to be the case. The upper bound on the admissible parameter appears not only to be a complication arising from generalization in the structure, but also by the approach itself. This is emphasized by the recent result in \cite{MRVV22}, where by means of a slightly adjusted classical approach to the intrinsic scaling by DiBenedetto (\cite{DeGiorgi57,dibenedettoHolderEstimatesNonlinear1985,dibenedettoDegenerateParabolicEquations1993,BenedettoGianazzaVespri-Harnack}), the local Hölder regularity of bounded weak solutions to a chemotaxis system with signal production was established for a divergence term of the form $a(x,t,u,\nabla u)=\nabla u^m-\chi u^{q-1}\nabla v$ with $m\geq1,\chi>0$ and $q\geq\max\{\tfrac{m+1}{2},2\}$. For $m>3$ this leads to $q-1>1$, so again we observe that a linear impact of $u$ in the sensitivity is not covered by the result whenever $m>3$.\smallskip

\textbf{Flexibility arising from structural separation of diffusion and convection.}
While the problems discussed above cannot simply be handled by taking a less general structure, we can still hope to at least earn another degree of flexibility by slightly reshaping the equation considered in \cite{PorzVesp93}. Being mostly interested in chemotaxis applications, we will shift towards a particularized form still encompassing most of the common chemotaxis equations. To be precise, we will consider a general degenerate equation of the form
\begin{align}\label{eq:gen-eq}
\left.
\begin{array}{r@{\ }l@{\quad}l@{\quad}l@{\,}c}
u_t&=\Delta\Phi(u)+\nabla\cdot\big(a(x,t)u\big)+b(x,t,u),\ &x\in\Omega,& t\in(0,T),
\end{array}\right.
\end{align}
in a smoothly bounded domain $\Omega\subset\R^\dimN$. Herein, $\Phi\in\CSp{0}{[0,\infty)}\cap\CSp{2}{(0,\infty)}$ is a convex function satisfying $\Phi(0)= 0$, $\Phi'>0$ on $(0,\infty)$ and being such that there is $C_\Phi>0$ with the property that for some small $s_0\in(0,1]$ the inequality
\begin{align}\label{eq:cond-phi}
s\Phi''(s)\leq C_\Phi\Phi'(s)\quad\text{for all }s\in[0,s_0]
\end{align}
holds. The functions $a,b$ are assumed to be such that
\begin{align}\label{eq:cond-a-b}
\begin{aligned}
a&\in L^{2q_1}\big((0,T);L^{2q_2}(\Omega;\R^\dimN)\big),\\
b&\in M\big(\OmT\times \R\big)\text{ such that }|b(x,t,\xi)|\leq\hat{b}(x,t)\ \text{a.e. for some }\hat{b}\in\LSpb{q_1}{(0,T);\LSp{q_2}{\Omega}\!}
\end{aligned}
\end{align}
for some $q_1,q_2>1$ such that
\begin{align}\label{eq:cond-kappa}
\frac{2}{q_1}+\frac{\dimN}{q_2}=2-\dimN\kappa,\quad\text{for some }\kappa\in(0,\tfrac{2}{\dimN}),
\end{align}
where $M(X)$ denotes the vector space of equivalence classes of measurable functions from $X$ to $\R$. Essentially, we separated the diffusion and convection terms and thereby also got rid of the dependence on $\nabla u$ for the convection term. We will combine this structural adaptation with a slightly different angle on the conventional approach to establish results that on one hand facilitate an improvement on the range for $m$ in the porous medium example above, and on the other also cover closely related diffusion terms not precisely of porous medium type. The fundamental ideas of the approach are rooted in DiBenedetto's observations on continuity of solutions to certain Stefan problems in \cite{DB82,DiBenedetto83} and DiBenedetto's and Friedman's methods in \cite{dibenedettoHolderEstimatesNonlinear1985}.\smallskip

Observe that certainly diffusion terms of porous medium type, i.e. $\Phi(s)=s^m$, $m\geq 1$ and, more generally, all convex polynomials with $\Phi(0)=0$ and $\Phi'(s)>0$ for $s>0$ satisfy the above conditions imposed on $\Phi$. In particular, all polynomials of the form $\Phi(s)=\sum_{k=1}^m a_ks^k$ with $a_k\geq0$ for $k\in\{1,\dots,m\}$ are included. The nonnegativity of $a_k$ in said polynomials, however, is not a necessary condition. A prototypical example going beyond the realm of nonnegative coefficients is given by $\Phi(s)=s^2(s^2-s+3)$.\smallskip

\textbf{Details on the approach.} The main difference to the conventional approach will consist of investigating a suitably transformed equation instead of \eqref{eq:gen-eq}. As it turns out, it will be beneficial to rewrite \eqref{eq:gen-eq} in terms of the inverse of $\Phi$. We note that by the conditions on $\Phi$, the inverse function $\Phi^{-1}:[0,\infty)\to[0,\infty)$ is well-defined, continuous, concave and strictly monotonically increasing. Moreover, as $\Phi'$ is positive in $(0,\infty)$, $\Phi^{-1}$ is differentiable in $(0,\infty)$. In fact, we are going to consider
\begin{align*}
\left(\Phi^{-1}\big(\Phi(u)\big)\right)_t=\Delta \Phi(u)+\nabla\cdot\left(a(x,t)\Phi^{-1}\big(\Phi(u)\big)\right)+b\left(x,t,\Phi^{-1}\big(\Phi(u)\big)\right),\quad x\in\Omega,\quad t\in(0,T),
\end{align*}
throughout most of the paper. Writing $w(x,t)=\Phi\big(u(x,t)\big)$, so that $u(x,t)=\Phi^{-1}\big(w(x,t)\big)$, we hence obtain the transformed version of \eqref{eq:gen-eq} given by
\begin{align}\label{eq:trans-gen-eq}
\big(\Phi^{-1}(w)\big)_t=\Delta w+\nabla\cdot\big(a(x,t)\Phi^{-1}(w)\big)+b\big(x,t,\Phi^{-1}(w)\big),\quad x\in\Omega,\quad t\in(0,T).
\end{align}
Written like this, \eqref{eq:trans-gen-eq} resembles the equations studied in \cite{DB82,DiBenedetto83}, albeit with a more regular function inside the time-derivative. This transformation has also been employed in \cite{HwangZhang19} for the special case of porous medium type diffusion. Before giving some details on the approach, let us first specify what we will consider a weak solution of \eqref{eq:trans-gen-eq} in $\OmT:=\Omega\times(0,T]$.

\begin{definition}\label{def:weaksol}
A nonnegative function $$w\in \CSploc{0}{(0,T);\LSploc{2}{\Omega}}\cap\LSplocb{2}{(0,T);W_{loc}^{1,2}(\Omega)}$$ will be called a weak solution of \eqref{eq:trans-gen-eq} in $\OmT$, if
\begin{align*}
\Phi^{-1}(w)\in\CSploc{0}{(0,T);\LSploc{2}{\Omega}}\quad\text{and}\quad a\Phi^{-1}(w)\in\LSploc{2}{\OmT;\R^\dimN}
\end{align*}
and if for every subdomain $\Omega'$ compactly contained in $\Omega$ and every interval $[t_1,t_2]\subset(0,T)$ the equality
\begin{align}\label{eq:weak-sol-eq}
\int_{\Omega'}\Phi^{-1}(w)\Psi\Big\vert_{t_1}^{t_2}+\int_{t_1}^{t_2}\!\int_{\Omega'}\nabla w\cdot\!\nabla\Psi&+\int_{t_1}^{t_2}\!\int_{\Omega'}a(x,t)\Phi^{-1}(w)\cdot\!\nabla\Psi\\&=\int_{t_1}^{t_2}\!\int_{\Omega'}\Phi^{-1}(w)\Psi_t+\int_{t_1}^{t_2}\!\int_{\Omega'}b\big(x,t,\Phi^{-1}(w)\big)\Psi\nonumber
\end{align}
holds for all $\Psi\in\CSp{0}{[t_1,t_2];\LSp{\infty}{\Omega'}}\cap \LSpb{2}{(t_1,t_2);W_{0}^{1,2}(\Omega')}$ with $\Psi_t\in\LSp{2}{(t_1,t_2);\LSp{2}{\Omega'}}$.
\end{definition}

To derive the energy estimates in Section~\ref{sec3}, we will rely on a slightly reformulated version of the solution concept above by introducing the notion of Steklov averages (\cite{LSU,chagasPropertiesSteklovAverages2017}). For $v\in\LSp{1}{\OmT}$ and $0<h<T$ the Steklov average of $v$ is defined as
$$[v(\cdot,t)]_h:=\begin{cases}\frac{1}{h}\int_t^{t+h}v(\cdot,\tau)\intd\tau,\quad&\text{if }t\in(0,T-h],\\0,&\text{if }t>T-h.\end{cases}$$

Let us note that different regularization processes are possible here. See e.g. \cite{MR0771666,MR4500278} for precedents of mollification using the exponential function.

\begin{remark}\label{rem:steklov_1}
The weak solution concept in Definition~\ref{def:weaksol} can be rephrased in terms of Steklov averages by requiring for every subdomain $\Omega'$ compactly contained in $\Omega$ and every $t\in(0,T-h)$ that
\begin{align}\label{eq:weak-sol-steklov}
\int_{\Omega'}\!\Big(\big[\Phi^{-1}\big(w(\cdot,t)\big)\big]_h\Big)_t\Psi+\int_{\Omega'}\!\big[\nabla w(\cdot,t)\big]_h\cdot\nabla\Psi&+\int_{\Omega'}\!\big[a(\cdot,t)\Phi^{-1}\big(w(\cdot,t)\big)\big]_h\cdot\nabla\Psi\nonumber\\&=\int_{\Omega'}\!\big[b\big(\cdot,t,\Phi^{-1}(w(\cdot,t)\big)\big]_h\Psi
\end{align}
for all $t\in(0,T-h)$ and all $\Psi\in W_0^{1,2}(\Omega')\cap\LSp{\infty}{\Omega'}$. (See also \cite[Section II.1]{dibenedettoDegenerateParabolicEquations1993}.)
\end{remark}

Starting with a bounded weak solution of \eqref{eq:trans-gen-eq}, we can adapt the local energy estimates used in the conventional approach for the transformed equation (Section~\ref{sec3}). Herein, the structural condition \eqref{eq:cond-phi} on $\Phi''$ and $\Phi'$ will be key to control some of the ill-signed integral terms. The all-underlying main objective will remain the same as in the standard approach, that is to establish a decay estimate for the essential oscillation when considered over nested cylinders (see Corollary~\ref{cor:iteration}). The main ingredients, the so-called alternatives, for the transformed equation have to be thoroughly checked for the intricate influence of $\Phi^{-1}$, which most of the time appears in the form of $(\Phi^{-1}(w))'$. This introduces some difficulty into the process, which have to be treated in a quite subtle fashion (we refer to Section~\ref{sec4:alt1} and \ref{sec4:alt2}). Building on these alternatives, the rigorous iteration process encompassing the construction of the nested cylinders will then be undertaken in Section~\ref{sec4:iter}. Having the essential oscillation on nested cylinders under control, Hölder regularity in compact subsets of $\OmT$ can be concluded in standard manner (Section~\ref{sec5}). The final Section~\ref{sec6:boundary} will then be devoted to extending the regularity results up to the boundary.\smallskip

\textbf{Main results.} Our first main result on the Hölder regularity of bounded weak solutions to \eqref{eq:trans-gen-eq} inside compact subsets of $\OmT$ reads as follows.

\begin{theorem}\label{theo:1}
Assume that $\Phi\in\CSp{0}{[0,\infty)}\cap\CSp{2}{(0,\infty)}$ is convex with $\Phi(0)=0$ and $\Phi'>0$ on $(0,\infty)$ and satisfies \eqref{eq:cond-phi}. Suppose \eqref{eq:cond-a-b} is valid with $q_1,q_2>1$ fulfilling \eqref{eq:cond-kappa} for some $\kappa\in(0,\frac{2}{\dimN})$. Let $w$ be a nonnegative weak solution of \eqref{eq:trans-gen-eq} in the sense of Definition \ref{def:weaksol} with $M>0$ such that
\begin{align}\label{eq:linfty-bound-w}
\big\|w(\cdot,t)\big\|_{\Lo[\infty]}\leq M\quad\text{ for all }t\in(0,T).
\end{align}
Then, up to a re-definition of $w$ on a null subset of $\OmT$, $(x,t)\mapsto w(x,t)$ is locally Hölder continuous in $\OmT$. In particular, there is $\alpha=\alpha(\Phi,M,q_1,q_2,a,\hat{b})$ $\in(0,1)$ such that, up to a re-definition of $w$ on a null subset of $\OmT$, for any compact $\Sigma_\tau\subset \OmT$ one can find $C=C(\Phi,M,q_1,q_2,a,\hat{b},\Sigma_\tau)>0$ such that the inequality
\begin{align*}
\big|w(x_0,t_0)-w(x_1,t_1)\big|\leq C\left(|x_0-x_1|+|t_0-t_1|^\frac{1}{2}\right)^{\alpha}
\end{align*}
holds for every pair $(x_0,t_0),(x_1,t_1)\in \Sigma_\tau$.
\end{theorem}

\begin{remark}
The dependencies on $\Phi$ in Theorem~\ref{theo:1} are quite general and could in principle be quantified more explicitly by taking into account other key properties characterizing $\Phi$. To obtain a more explicit specification one could, for instance, replace the dependence of $\alpha$ and $C$ on $\Phi$ by a dependence on the constants $s_0$ and $C_\Phi$ from \eqref{eq:cond-phi}.
\end{remark}

If $\Phi^{-1}:[0,M]\to\R$ is also assumed to be Hölder continuous we readily obtain the following evident consequence for $u$.

\begin{corollary}\label{cor:theo1}
If, additionally, $\Phi^{-1}:[0,M]\to\R$, $s\mapsto \Phi^{-1}(s)$ is Hölder continuous with Hölder exponent $\beta\in(0,1)$, then $(x,t)\mapsto u(x,t)$ is also Hölder continuous with exponent $\alpha\beta\in(0,1)$.
\end{corollary}

\begin{remark}
The Hölder continuity of $\Phi^{-1}$ can conveniently be verified in the following ways:\\
i) If $s_0>0$ and $q>0$ are such that $\int_0^{s_0}\frac{1}{|\Phi'(s)|^q}\intd s<\infty$, then $\Phi^{-1}$ is Hölder continuous with $\beta=\frac{q}{q+1}$.\\
ii) If $s_0>0$, $p\geq 1$ and $C_{\Phi'}>0$ are such that $\Phi'(s)\geq  C_{\Phi'} s^{p-1}$ on $[0,s_0]$, then  $\Phi^{-1}$ is Hölder continuous with $\beta=\frac{1}{p}$.\\
iii) If $\lim_{s\searrow0}\Phi'(s)>0$, then $\Phi^{-1}$ is Hölder continuous with $\beta=1$.	
\end{remark}

Theorem~\ref{theo:1} and Corollary~\ref{cor:theo1} encompass the Hölder regularity results of both \cite{HwangZhang19} and \cite{MRVV22} and are not restricted to porous medium-type diffusion, thereby allowing to treat a wider array of nonlinear diffusion processes in the equation.

In a second step we will consider a corresponding initial-boundary value problem with no-flux boundary condition. In fact, supposing that $\romega$ is of class $C^{1}$ and assuming that the nonnegative initial data satisfy $w_0\in\CSp{\beta_0}{\bomega}$, we are going to consider
\begin{align}\label{eq:bvp-w}
\left\{
\begin{array}{r@{\,}l@{\quad}l}
\big(\Phi^{-1}(w)\big)_t=\Delta w+\nabla\cdot\big(&\!a(x,t)\Phi^{-1}(w)\big)+b\big(x,t,\Phi^{-1}(w)\big),\ &(x,t)\in\OmT,\\
\big(\nabla w+a(x,t)\Phi^{-1}(w)\big)\cdot\nu&=g\big(x,t,\Phi^{-1}(w)\big),\ &(x,t)\in \romega\times[0,T),\\
w(x,0)&=w_0(x),&x\in\Omega,
\end{array}\right.
\end{align}
where $\nu$ denotes the outward unit normal vector and the boundary data $g\in \CSp{0}{\romega\times[0,T]\times\R}$ is assumed to admit an extension onto $\Omega$ for a.e. $t\in(0,T)$ denoted by $\hat{g}$, such that given some $M>0$ the inequalities
\begin{align}\label{eq:cond-g}
\big|\hat{g}\big(x,t,\xi\big)\big|&\leq g_0(x,t)\Phi(\xi)\quad\text{for a.e. } (x,t)\in\OmT\text{ and }\xi\in[0,M],\nonumber\\
\big|\partial_\xi \hat{g}\big(x,t,\xi\big)\big|&\leq g_0(x,t)\Phi'\big(\xi\big)\quad\text{for a.e. } (x,t)\in\OmT\text{ and }\xi\in[0,M],\\
\big|\partial_{x_i} \hat{g}\big(x,t,\xi\big)\big|&\leq g_0^2(x,t),\ i\in\{1,\dots,\dimN\}\ \text{for a.e. } (x,t)\in\OmT\text{ and }\xi\in[0,M],\nonumber
\end{align}
hold with some nonnegative $g_0\in\LSp{2q_1}{(0,T);\Lo[2q_2]}$, where $q_1,q_2>1$ satisfy \eqref{eq:cond-kappa} for some $\kappa\in(0,\frac{2}{\dimN})$. A weak solution to \eqref{eq:bvp-w} in $\bomega\times[0,T)$ will be specified in the following way.
\begin{definition}\label{def:weaksol-bv}
A nonnegative function $$w\in \CSp{0}{[0,T];\LSp{2}{\Omega}}\cap\LSpb{2}{(0,T);W^{1,2}(\Omega)}$$ will be called a weak solution of \eqref{eq:trans-gen-eq}, if
\begin{align*}
\Phi^{-1}(w)\in\CSp{0}{[0,T];\LSp{2}{\Omega}}\quad\text{and}\quad a\Phi^{-1}(w)\in\LSp{2}{\OmT;\R^\dimN}
\end{align*}
and if for every bounded domain $\mathcal{B}\subset \R^\dimN$ the equality
\begin{align}\label{eq:weak-sol-bv-eq}
&\intomega\Phi^{-1}\big(w(\cdot,T)\big)\Psi(\cdot,T)+\int_{0}^{T}\!\int_{\mathcal{B}\cap\Omega}\!\!\nabla w\cdot\!\nabla\Psi+\int_{0}^{T}\!\int_{\mathcal{B}\cap\Omega}\!\!a(x,t)\Phi^{-1}(w)\cdot\!\nabla\Psi\\
=&\intomega\Phi^{-1}(w_0)\Psi(\cdot,0)+\int_{0}^{T}\!\int_{\mathcal{B}\cap\Omega}\!\!\Phi^{-1}(w)\Psi_t+\int_{0}^{T}\!\int_{\mathcal{B}\cap\Omega}\!\!b\big(x,t,\Phi^{-1}(w)\big)\Psi+\int_{0}^{T}\!\int_{\mathcal{B}\cap\romega}\!\!g\big(x,t,\Phi^{-1}(w)\big)\Psi\nonumber
\end{align}
holds for all $\Psi\in\CSp{0}{[0,T];\LSp{\infty}{\mathcal{B}}}\cap \LSplocb{2}{(0,T);W_{0}^{1,2}(\mathcal{B})}$ with $\Psi_t\in\LSp{2}{(0,T);\LSp{2}{\mathcal{B}}}$.
\end{definition}

\begin{remark}\label{rem:steklov_2}
The solution concept of Definition~\ref{def:weaksol} can also be rephrased in terms of Steklov averages by replacing \eqref{eq:weak-sol-bv-eq} with
\begin{align}\label{eq:weak-sol-bv-steklov}
\int_{\mathcal{B}\cap\Omega}\!\Big(\big[\Phi^{-1}\big(w(\cdot,t)\big)\big]_h\Big)_t\Psi&+\int_{\mathcal{B}\cap\Omega}\!\big[\nabla w(\cdot,t)\big]_h\cdot\nabla\Psi+\int_{\mathcal{B}\cap\Omega}\! \big[a(\cdot,t)\Phi^{-1}\big(w(\cdot,t)\big)\big]_h\cdot\nabla\Psi\nonumber\\
&=\int_{\mathcal{B}\cap\Omega}\!\Big[b\big(\cdot,t,\Phi^{-1}\big(w(\cdot,t)\big)\big)\Big]_h\Psi+\int_{\mathcal{B}\cap\romega}\!\Big[g\big(\cdot,t,\Phi^{-1}\big(w(\cdot,t)\big)\big)\Big]_h\Psi
\end{align}
for all $t\in(0,T-h)$ and all $\Psi\in W_0^{1,2}(\mathcal{B})$. (Compare \cite[Section II.2-ii)]{dibenedettoDegenerateParabolicEquations1993}.)
\end{remark}

For the initial-boundary-value problems in \eqref{eq:bvp-w} we then obtain the following main result.

\begin{theorem}\label{theo:2}
Let $M>0$. Assume that $\Phi\in\CSp{0}{[0,\infty)}\cap\CSp{2}{(0,\infty)}$ is convex with $\Phi(0)\geq0$ and $\Phi'>0$ on $(0,\infty)$ and satisfies \eqref{eq:cond-phi}. Suppose that $\romega$ is of class $C^{1}$, that \eqref{eq:cond-a-b} and \eqref{eq:cond-g} are valid with $q_1,q_2>1$ satisfying \eqref{eq:cond-kappa} for some $\kappa\in(0,\frac{2}{\dimN})$ and that $w_0\in\CSp{\beta_0}{\bomega}$ for some $\beta_0\in(0,1)$ is nonnegative. Let $w$ be a nonnegative weak solution of \eqref{eq:bvp-w} in the sense of Definition~\ref{def:weaksol-bv} such that
\begin{align}\label{eq:linfty-bound-w-2}
\big\|w(\cdot,t)\big\|_{\Lo[\infty]}\leq M\quad\text{for all }t\in(0,T]
\end{align}
holds. Then, up to a re-definition of $w$ on a null subset of $\bomega\times[0,T]$, $(x,t)\mapsto w(x,t)$ is Hölder continuous in $\bomega\times[0,T]$. In particular, there are $\alpha=\alpha(\Phi,M,q_1,q_2,a,\hat{b},g_0,\beta_0)\in(0,1)$ and $C=C(\Phi,M,q_1,q_2,a,\hat{b},g_0,w_0)>0$ such that, up to a re-definition of $w$ on a null subset of $\bomega\times[0,T]$,
\begin{align*}
\big|w(x_0,t_0)-w(x_1,t_1)\big|\leq C\left(|x_0-x_1|+|t_0-t_1|^\frac{1}{2}\right)^{\alpha}
\end{align*}
holds for every pair $(x_0,t_0),(x_1,t_1)\in \bomega\times[0,T]$.
\end{theorem}

As a counterpart to Corollary~\ref{cor:theo1}, we find that if $\Phi^{-1}:[0,M]\to\R$ is Hölder continuous and $u_0$ is given by $u_0:=\Phi^{-1}(w_0)$, then the result above also entails a regularity result for the associated initial-boundary value problem for $u$ given by
\begin{align*}
\left\{
\begin{array}{r@{\;}l@{\quad}l}
u_t=\Delta\Phi(u)&+\,\nabla\cdot\big(a(x,t)u\big)+b(x,t,u),\ &(x,t)\in\OmT,\\
\big(\nabla\Phi(u)&+\,a(x,t)u\big)\cdot\nu=g(x,t,u),\ &(x,t)\in \romega\times[0,T],\\
u(x,0)&=u_0(x),&x\in\Omega.
\end{array}\right.
\end{align*}

\begin{corollary}\label{cor:theo2}
If, additionally, $\Phi^{-1}:[0,M]\to\R$, $s\mapsto \Phi^{-1}(s)$ is Hölder continuous with Hölder exponent $\beta\in(0,1)$, then $(x,t)\mapsto u(x,t)$ is also Hölder continuous in $\bomega\times[0,T]$ with exponent $\alpha\beta\in(0,1)$.
\end{corollary}

With Theorem~\ref{theo:2} and Corollary~\ref{cor:theo2} we rigorously extend the results of \cite{HwangZhang19,MRVV22} beyond compact subsets of $\OmT$. Specifically, the results given above allow the treatment of no-flux boundary conditions commonly encountered in the frameworks of chemotaxis equations.

\begin{remark}\label{rem:boundary-data}
If in addition to the previous assumptions we also have $\Phi\in C^2([0,\infty))$, minor adjustments in the proof allow also for treatment of boundary data where the extension $\hat{g}$ fulfills
\begin{align*}
\big|\hat{g}(x,t,\xi)\big|&\leq g_0(x,t)\Phi(\xi)+g_1(x,t)\quad\text{for a.e. }(x,t)\in\OmT\text{ and } \xi\in[0,M],\\
\big|\partial_{\xi}\hat{g}(x,t,\xi)\big|&\leq g_0(x,t)\Phi'(\xi)\quad\text{for a.e. }(x,t)\in\OmT\text{ and } \xi\in[0,M],\\
\big|\partial_{x_i}\hat{g}(x,t,\xi)\big|&\leq g_1^2(x,t),\quad i\in\{1,\dots,\dimN\}\ \text{for a.e. }(x,t)\in\OmT\text{ and } \xi\in[0,M]
\end{align*}
with nonnegative $g_0,g_1\in\LSp{2q_1}{(0,T);\Lo[2q_2]}$ and $q_1,q_2>1$ satisfying \eqref{eq:cond-kappa} for some $\kappa\in(0,\frac{2}{\dimN})$, instead of \eqref{eq:cond-g}. In this case one has to replace the exponent $-\frac{q_1-1}{q_1}$ in the alternatives corresponding to Lemma~\ref{lem:first-alt-init} and \ref{lem:second-alt-init} with $-\frac{2q_1-1}{q_1}$, which leads to smaller admissible values for $\beta$ in Lemma~\ref{lem:iteration-1-init} and hence also slightly smaller Hölder exponents.
\end{remark}

\setcounter{equation}{0} 
\section{Preliminary estimates and well-established results}\label{sec2:prelim}
Let us begin with gathering some additional properties of the function $\Phi$ and other well-known results we will use throughout the paper. The estimates presented in the following lemma are direct consequences of the condition \eqref{eq:cond-phi}. Together with the convexity assumption on $\Phi$ -- and the entailed concavity of $\Phi^{-1}$ -- these will lay the groundwork for the treatment of $\Phi$ and $\Phi^{-1}$ in the local energy estimates of Section~\ref{sec3} and the alternatives discussed in Section~\ref{sec4}.
\begin{lemma}\label{lem:phi}
Let $\Phi\in\CSp{0}{[0,\infty)}\cap\CSp{2}{(0,\infty)}$ be a convex function satisfying $\Phi(0)= 0$ and $\Phi'>0$ on $(0,\infty)$ and assume that \eqref{eq:cond-phi} holds for some $C>0$ and $s_0\in(0,1]$. Then, there is $C_1=C_1(\Phi)>1$ such that
\begin{align}\label{eq:phi-prop-1}
s\Phi'(s)\leq C_1\Phi(s)\quad\text{for all }s\in[0,s_0],
\end{align}
and for all $\alpha>1$ there is $C_2=C_2(\Phi,\alpha)>1$ such that
\begin{align}\label{eq:phi-prop-2}
\Phi'(\alpha s)\leq C_2\Phi'(s)\quad\text{for all }s\in\big[0,\frac{s_0}{\alpha}\big].
\end{align}
Moreover, for all $\alpha>1$ and $s_1>0$ one can find $C=C(\Phi,\alpha,s_1)>1$ such that
\begin{align}\label{eq:phi-prop-3}
\Phi'(\alpha s)\leq C\Phi'(s),\quad s\Phi''(s)\leq C\Phi'(s)\quad\text{and}\quad s\Phi'(s)\leq C\Phi(s)\quad\text{for all }s\in[0,s_1].
\end{align}
\end{lemma}

\begin{bew}
Integrating \eqref{eq:cond-phi} we find that since $\Phi(0)=0$,
\begin{align*}
C_{\Phi}\Phi(s)= C_\Phi\big(\Phi(s)-\Phi(0)\big)=C_\Phi\int_0^s \Phi'(\sigma)\intd\sigma\geq\int_0^s \sigma \Phi''(\sigma)\intd\sigma\quad\text{for }s\in[0,s_0].
\end{align*}
In view of an integration by parts on the right, again using $\Phi(0)=0$, entails
\begin{align*}
C_\Phi\Phi(s)\geq s\Phi'(s)-\int_0^s \Phi'(\sigma)\intd\sigma= s\Phi'(s)-\Phi(s)
\end{align*}
for $s\in[0,s_0]$. Therefore, we conclude that
\begin{align*}
s\Phi'(s)\leq (C_\Phi+1)\Phi(s)=:C_1\Phi(s)\quad\text{for all }s\in[0,s_0]
\end{align*}
and hence obtain \eqref{eq:phi-prop-1}. As for \eqref{eq:phi-prop-2}, we note that the claim is obviously true for $s=0$, while for $s>0$ \eqref{eq:cond-phi} is equivalent to
\begin{align*}
\big(\ln\big(\Phi'(s)\big)\big)'\leq \big(\ln\big(s^{C_\Phi}\big)\big)'\quad\text{for all }s\in(0,s_0].
\end{align*}
Integrating this from $s$ to $\alpha s$ for $s\in\big(0,\frac{s_0}{\alpha}\big]$ we find that
\begin{align*}
\ln\big(\Phi'(\alpha s)\big)-\ln\big(\Phi'(s)\big)\leq \ln\big((\alpha s)^{C_\Phi}\big)-\ln\big(s^{C_\Phi}\big)\quad\text{for all }s\in\big(0,\frac{s_0}{\alpha}\big],
\end{align*}
which entails $$\frac{\Phi'(\alpha s)}{\Phi'(s)}\leq\alpha^{C_\Phi}=:C_2\quad\text{for all }s\in\big(0,\frac{s_0}{\alpha}\big].$$
Finally, assuming $s_1>s_0$, the extension of the inequalities \eqref{eq:phi-prop-1} and \eqref{eq:phi-prop-2} for $s\in[\frac{s_0}{\alpha},s_1]$ and $s\in[s_0,s_1]$, respectively, is then an evident consequence of the fact that by the convexity of $\Phi$ the number $C_{3}=C_{3}(\Phi,s_1):=\max_{s\in[s_0,s_1]}\Phi''(s)\geq 0$ is well-defined, and the fact that $\Phi'$ and $\Phi$ are both monotonically increasing. Taking $C=C(\Phi,\alpha,s_1):=\max\big\{C_1,C_2,\frac{\Phi'(\alpha s_1)}{\Phi'(\tfrac{s_0}{\alpha})},\frac{C_3s_1}{\Phi'(s_0)},\frac{\Phi'(s_1)s_1}{\Phi(s_0)}\big\}$ we arrive at the asserted inequalities in \eqref{eq:phi-prop-3}.
\end{bew}

The following general Poincaré-type inequality was introduced by De Giorgi (\cite{DeGiorgi57}). The version given here can be found in e.g. \cite[Lemma~I.2.2 and Remark~I.2.2]{dibenedettoDegenerateParabolicEquations1993} and \cite[Lemma~2.2.2 and Remark~2.1]{BenedettoGianazzaVespri-Harnack}.
\begin{lemma}\label{lem:DeGiorgi}
Let $\dimN\geq1$ be an integer. There exists a constant $C=C(\dimN)>0$ such that whenever $x_0\in\R^\dimN$, $\rd>0$, $L,K\in\R$ satisfy $L>K$ and $\varphi\in W^{1,1}(\mK_\rd(x_0))$ for the $\dimN$-dimensional cube centered at $x_0$ with edge length $2\rd$ denoted by $\mK_\rd(x_0)$, then 
\begin{align*}
(L-K)\big|\big\{x\in \mK_\rd(x_0): \varphi(x)>L\big\}\big|\leq \frac{C\rd^{\dimN+1}}{\big|\big\{x\in \mK_\rd(x_0):\varphi(x)<K\big\}\big|}\int_{\big\{x\in \mK_\rd(x_0):K<\varphi(x)<L\big\}}\big|\nabla\varphi(x)\big|\intd x.
\end{align*}
\end{lemma}

When establishing our alternatives for the level sets we will encounter parabolic spaces of the form $\LSp{\infty}{(0,T);\Lo[p]}\cap\LSpb{p}{(0,T);\W[1,p]}$ for $p\geq 1$. These spaces embed into $\LSp{q}{\OmT}$ for $q=p(1+\tfrac{p}{\dimN})$ and in particular we have the following. (See also \cite[Corollary 2.4.1]{BenedettoGianazzaVespri-Harnack}.)
\begin{corollary}\label{cor:sob}
Let $\dimN\geq1$ be an integer, $\Omega\subset\R^\dimN$ be a bounded domain with Lipschitz boundary, $T\in(0,1]$ and $p\geq 1$. There is $C=C(\dimN,\Omega,p)>0$ such that for any $\varphi\in\LSp{\infty}{(0,T);\Lo[p]}\cap\LSpb{p}{(0,T);\W[1,p]}$ the inequality
\begin{align*}
\intoTomega|\varphi|^p\leq C\Big|\big\{(x,t)\in\OmT:|\varphi(x,t)|>0\big\}\Big|^\frac{p}{p+\dimN}\bigg(\esssup_{t\in(0,T)}\|\varphi(\cdot,t)\|_{\Lo[p]}+\Big(\intoTomega|\nabla \varphi|^p\Big)^\frac{1}{p}\bigg)^p
\end{align*}
holds.
\end{corollary}

\begin{bew}
Denoting by $\bar{\varphi}:=\frac{1}{|\Omega|}\intomega\varphi$ the spatial average of $\varphi$, we make use of Hölder's inequality to estimate
\begin{align*}
\intoTomega|\varphi|^p=\iint\limits_{\{|\varphi|>0\}}|\varphi|^p&\leq\Big(\iint\limits_{\{|\varphi|>0\}}|\varphi|^\frac{p(p+\dimN)}{\dimN}\Big)^\frac{\dimN}{p+\dimN}\Big(\iint\limits_{\{|\varphi|>0\}}1\Big)^\frac{p}{p+\dimN}\\
&\leq 2^p\big|\big\{|\varphi|>0\big\}\big|^\frac{p}{p+\dimN}\Big(\intoTomega|\varphi-\bar{\varphi}|^\frac{p(p+\dimN)}{\dimN}+\intoTomega|\bar{\varphi}|^\frac{p(p+\dimN)}{\dimN}\Big)^\frac{\dimN}{p+\dimN}.
\end{align*} 
Drawing on the \GNI\ we find $C_1=C_1(\Omega,\dimN,p)>0$ such that for all $\varphi$ of said regularity class,
\begin{align*}
\intoTomega|\varphi|^p\leq 2^p\big|\big\{|\varphi|>0\big\}\big|^\frac{p}{p+\dimN}\Big(C_1\intoT\big\|\varphi-\bar{\varphi}\big\|_{\W[1,p]}^p\big\|\varphi-\bar{\varphi}\|_{\Lo[p]}^\frac{p^2}{\dimN}+\intoT\big\|\bar{\varphi}\big\|_{\Lo[\frac{p(p+\dimN)}{\dimN}]}^\frac{p(p+\dimN)}{\dimN}\Big)^\frac{\dimN}{p+\dimN}.
\end{align*}
Herein, we can rely on a Poincaré--Wirtinger inequality to obtain a constant $C_2=C_2(\dimN,\Omega,p)>0$ such that $\|\varphi-\bar{\varphi}\|_{\W[1,p]}^p\leq C_2\|\nabla\varphi\|_{\Lo[p]}^p$, so that with $C_3=C_3(\dimN,\Omega,p):=2^p(1+C_1C_2)^\frac{\dimN}{p+\dimN}>0$
\begin{align}\label{eq:emb-eq-1}
\intoTomega|\varphi|^p&\leq C_3\big|\big\{|\varphi|>0\big\}\big|^\frac{p}{p+\dimN}\Big(\intoT\big\|\nabla\varphi\big\|_{\Lo[p]}^p\big\|\varphi-\bar{\varphi}\|_{\Lo[p]}^\frac{p^2}{\dimN}+\intoT\big\|\bar{\varphi}\big\|_{\Lo[\frac{p(p+\dimN)}{\dimN}]}^\frac{p(p+\dimN)}{\dimN}\Big)^\frac{\dimN}{p+\dimN}\nonumber\\
&=:C_3\big|\big\{|\varphi|>0\big\}\big|^\frac{p}{p+\dimN}\big(I_1+I_2\big)^\frac{\dimN}{p+\dimN}.
\end{align}
To estimate further, we note that for arbitrary $q\geq1$
\begin{align}\label{eq:est-barphi-eq1}
\|\bar{\varphi}\|_{\Lo[q]}\leq |\Omega|^{\frac{1}{q}-\frac{1}{p}}\|\varphi\|_{\Lo[p]}\quad\text{ and }\quad \big\|\varphi-\bar{\varphi}\big\|_{\Lo[p]}\leq 2\|\varphi\|_{\Lo[p]}.
\end{align}
Accordingly, using the assumption $T\leq 1$ we find that
\begin{align*}
I_2\leq\esssup_{t\in(0,T)}\big\|\varphi(\cdot,t)\|_{\Lo[p]}^\frac{p(p+\dimN)}{\dimN}|\Omega|^{-\frac{p}{\dimN}},
\end{align*}
whereas for $I_1$ we draw on Young's inequality in addition to \eqref{eq:est-barphi-eq1} and conclude that there is $C_4=C_4(\dimN,\Omega,p)>0$ such that
\begin{align*}
I_1&\leq 2^\frac{p^2}{\dimN}\esssup_{t\in(0,T)}\big\|\varphi(\cdot,t)\big\|_{\Lo[p]}^\frac{p^2}{\dimN}\intoTomega\big|\nabla\varphi\big|^p\leq C_4\esssup_{t\in(0,T)}\big\|\varphi(\cdot,t)\big\|_{\Lo[p]}^\frac{p(p+\dimN)}{\dimN}+\Big(\intoTomega\big|\nabla\varphi\big|^p\Big)^{\frac{p+\dimN}{\dimN}}.
\end{align*}
Plugging these estimates back into \eqref{eq:emb-eq-1} and rearranging exponents we obtain $C_5=C_5(\dimN,\Omega,p)>0$ satisfying
\begin{align*}
\intoTomega|\varphi|^p&\leq C_5\big|\big\{|\varphi|>0\big\}\big|^\frac{p}{p+\dimN}\Big(\esssup_{t\in(0,T)}\big\|\varphi(\cdot,t)\big\|_{\Lo[p]}^\frac{p(p+\dimN)}{\dimN}+\Big(\intoTomega\big|\nabla\varphi\big|^p\Big)^{\frac{p+\dimN}{\dimN}}\Big)^{\frac{\dimN}{p+\dimN}\cdot\frac{p}{p}}\\
&\leq C_5\big|\big\{|\varphi|>0\big\}\big|^\frac{p}{p+\dimN}\Big(\esssup_{t\in(0,T)}\big\|\varphi(\cdot,t)\|_{\Lo[p]}+\Big(\intoTomega\big|\nabla\varphi\big|^p\Big)^\frac{1}{p}\Big)^{p}
\end{align*}
for any $\varphi\in\LSp{\infty}{(0,T);\Lo[p]}\cap\LSpb{p}{(0,T);\W[1,p]}$.
\end{bew}

More generally, one can establish the following embedding result (cf. \cite[Proposition~I.3.3]{dibenedettoDegenerateParabolicEquations1993} and \cite[Proposition~2.4.2]{BenedettoGianazzaVespri-Harnack}), which -- following steps similar to the proof above -- is a direct consequence of the Gagliardo--Nirenberg and Poincaré--Wirtinger inequalities.

\begin{lemma}\label{lem:lemmaX}
Let $\Omega\subset\R^\dimN$, $\dimN\geq 1$, be a bounded domain with Lipschitz boundary, $T\in(0,1]$ and $p>1$. There is a constant $C=C(\dimN,\Omega,p,q_1,q_2)>0$ such that whenever $q_1,q_2\geq 1$ satisfy
\begin{align*}\begin{cases}
q_1\in[p^2,\infty),\quad q_2\in(p,\infty]&\quad\text{for }\dimN=1,\\
q_1\in(p,\infty),\quad q_2\in(p,\tfrac{\dimN p}{\dimN-p})&\quad\text{for }1<p<\dimN,\\
q_1\in(\tfrac{p^2}{\dimN},\infty),\quad q_2\in(p,\infty)&\quad\text{for }1<\dimN\leq p\end{cases}
\end{align*}
and $$\frac{p}{q_1}+\frac{\dimN}{q_2}=\frac{\dimN}{p},$$
then
\begin{align*}
\big\|\varphi\big\|_{\LSp{q_1}{(0,T);\Lo[q_2]}}\leq C\bigg(\esssup_{t\in(0,T)}\big\|\varphi(\cdot,t)\big\|_{\Lo[p]}+\Big(\intoTomega|\nabla \varphi|^p\Big)^\frac{1}{p}\bigg)
\end{align*}
for all $\varphi\in\LSp{\infty}{(0,T);\Lo[p]}\cap\LSpb{p}{(0,T);\W[1,p]}$.
\end{lemma}

\begin{bew}
Similarly to the steps in the corollary above, we make use of the Minkowski, Gagliardo--Nirenberg and Poincaré--Wirtinger inequalities, to first find $C_1=C_1(\dimN,\Omega,p,q_1,q_2)>0$ such that
\begin{align*}
\big\|\varphi\big\|_{\LSp{q_1}{(0,T);\Lo[q_2]}}&\leq2^{\frac{q_1-1}{q_1}}\Big(\intoT\|\varphi-\bar{\varphi}\|_{\Lo[q_2]}^{q_1}+\intoT\|\bar{\varphi}\|_{\Lo[q_2]}^{q_1}\Big)^\frac{1}{q_1}\\
&\leq C_1\Big(\intoT\|\nabla\varphi\|_{\Lo[p]}^{q_1\alpha}\|\varphi-\bar{\varphi}\|_{\Lo[p]}^{q_1(1-\alpha)}+\intoT\big\|\bar{\varphi}\big\|_{\Lo[q_2]}^{q_1}\Big)^{\frac{1}{q_1}},
\end{align*}
with $$\alpha=\frac{\dimN}{p}-\frac{\dimN}{q_2}=\frac{p}{q_1}\quad\text{and}\quad 1-\alpha=1-\frac{p}{q_1}.$$
Here, the assumption $T\leq 1$ and the estimates 
\begin{align*}
\|\bar{\varphi}\|_{\Lo[q]}\leq |\Omega|^{\frac{1}{q}-\frac{1}{p}}\|\varphi\|_{\Lo[p]},\ q\geq1, \quad\text{and}\quad \big\|\varphi-\bar{\varphi}\big\|_{\Lo[p]}\leq 2\|\varphi\|_{\Lo[p]},
\end{align*}
entail that there is $C_2=C_2(\dimN,\Omega,p,q_1,q_2)>0$ such that
\begin{align*}
\big\|\varphi\big\|_{\LSp{q_1}{(0,T);\Lo[q_2]}}&\leq C_2\Big(\intoT\|\nabla\varphi\|_{\Lo[p]}^p\|\varphi\|_{\Lo[p]}^{q_1-p}+\intoT\|\varphi\|_{\Lo[p]}^{q_1}\Big)^\frac{1}{q_1}\\
&\leq C_2\Big(\esssup_{t\in(0,T)}\big\|\varphi(\cdot,t)\big\|_{\Lo[p]}^{q_1-p}\intoTomega\big|\nabla\varphi\big|^p+\esssup_{t\in(0,T)}\big\|\varphi(\cdot,t)\big\|_{\Lo[p]}^{q_1}\Big)^\frac{1}{q_1}.
\end{align*}
Making use of Young's inequality with exponents $\frac{q_1}{p},\frac{q_1}{q_1-p}>1$ yields $C_3=C_3(\dimN,\Omega,p,q_1,q_2)>0$ satisfying
\begin{align*}
\big\|\varphi\big\|_{\LSp{q_1}{(0,T);\Lo[q_2]}}&\leq C_3\Big(\esssup_{t\in(0,T)}\big\|\varphi(\cdot,t)\big\|_{\Lo[p]}^{q_1}+\Big(\intoTomega\big|\nabla\varphi\big|^p\Big)^\frac{q_1}{p}\Big)^\frac{1}{q_1},
\end{align*}
from which we conclude the claim upon rearranging exponents.
\end{bew}

As final preparatory result, we state a fast geometric convergence results for two sequences of positive numbers obeying certain kind of recursive inequalities. This result can be found e.g. in \cite[Lemma I.4.2]{dibenedettoDegenerateParabolicEquations1993} or \cite[Lemma~C.3.3]{BenedettoGianazzaVespri-Harnack}.
\begin{lemma}\label{lem:fast-geo-conv}
Let $\{Y_i\}$ and $\{Z_i\}$, $i=0,1,\dots$, be sequences of positive numbers satisfying the inequalities
\begin{align*}
Y_{i+1}&\leq C b^i\big(Y_i^{1+\alpha}+Y_i^\alpha Z_i^{1+\kappa}\big)\quad\text{and}\quad
Z_{i+1}\leq C b^i\big(Y_i+Z_i^{1+\kappa}\big),
\end{align*}
for constants $C,b>1$ and $\kappa,\alpha>0$. If  
\begin{align*}
Y_0+Z_0^{1+\kappa}\leq (2C)^{-\frac{1+\kappa}{\sigma}}b^{-\frac{1+\kappa}{\sigma^2}}\quad\text{with }\sigma=\min\{\kappa,\alpha\},
\end{align*}
then $Y_i\to0$ and $Z_i\to 0$ as $i\to\infty$.
\end{lemma}

\setcounter{equation}{0} 
\section{Energy estimates localized in time-space cylinders}\label{sec3}
In preparation of the utmost important Section~\ref{sec4}, we will prepare some energy-type estimates by certain testing procedures. These suitably localized and truncated estimates will represent our main source of information on the solution behavior. The general form of these estimates is quite standard and has been used in many contexts. (See also \cite{dibenedettoDegenerateParabolicEquations1993}.) In the setting of \eqref{eq:trans-gen-eq}, however, we have to account for necessary changes regarding $\Phi^{-1}$. Since most of these modifications explicitly involve the derivative $(\Phi^{-1})'$, the degeneracy of $\Phi$ provides an even more difficult hurdle to overcome in this framework and makes the required alterations to the procedure non-obvious. The approach we are undertaking here draws on some arguments in \cite{HwangZhang19}, where local Hölder regularity was considered for the easier and less general case $\Phi(s)=s^m$, $m\geq1$ and $b=0$, and the main ideas go back to DiBenedetto's continuity results regarding a Stefan problem in \cite{DB82,DiBenedetto83}.\smallskip

For convenience in presentation, let us once more gather some of the standing assumptions we will impose throughout the paper. We will always assume the following:
\begin{align*}\label{eq:standing_assump}\tag{$A$}
\begin{cases}
&\Phi\in\CSp{0}{[0,\infty)}\cap\CSp{2}{(0,\infty)}\ \text{is convex and satisfies}\ \Phi(0)= 0, \Phi'>0\text{ on }(0,\infty)\text{ and }\eqref{eq:cond-phi},\nonumber\\
&\Phi^{-1}\in\CSp{0}{[0,\infty)}\cap\CSp{1}{(0,\infty)}\text{ is well-defined, concave and monotonically increasing},\\
&\eqref{eq:cond-a-b}\text{ is valid with } q_1,q_2>1\text{ satisfying }\eqref{eq:cond-kappa}\text{ with some }\kappa\in(0,\tfrac{2}{\dimN}).\nonumber
\end{cases}
\end{align*}
Moreover, whenever we call $w$ a weak solution of \eqref{eq:trans-gen-eq} or \eqref{eq:bvp-w}, we mean a weak solution in the sense of Definition~\ref{def:weaksol} and Definition~\ref{def:weaksol-bv}, respectively. As common for investigations concerning Hölder regularity of $w$, our aim will be to establish a decay estimate for the oscillation of $w$ with respect to a sequence of shrinking nested cylinders. By cylinders we actually mean boxes, which is a standard naming convention in the literature (see e.g. \cite{dibenedettoDegenerateParabolicEquations1993,BenedettoGianazzaVespri-Harnack}). In fact, we let $x_0\in\Omega$, $0\leq T_0<T_1\leq T$ and $\sqrt{\dimN}\rd<\dist(x_0,\romega)$ and introduce the notations
\begin{align*}
\mK_\rd(x_0)&:=\big\{x\in\R^\dimN\ \big\vert\ \|x-x_0\|_{\infty}<\rd\big\},\nonumber\\
Q_\rd(x_0,T_0,T_1)&:=\mK_\rd(x_0)\times(T_0,T_1],\\
\partial_p Q_\rd(x_0,T_0,T_1)&:=\partial \mK_\rd(x_0)\times[T_0,T_1]\cup \mK_\rd(x_0)\times\{T_0\}\nonumber
\end{align*}
for the open spatial cube centered at $x_0$ with edge length $2\rd$, the cylindrical box with top point $(x_0,T_1)$ and the parabolic boundary of $Q_\rd(x_0,T_0,T_1)$, respectively. Note that the assumption $\sqrt{\dimN}\rd<\dist(x_0,\romega)$ implies that $\mK_\rd(x_0)\subset\Omega$. We say $\rd$ is the radius of such a cylinder and $\rho=T_1-T_0$ its height. Additionally, for $K\in\R$ we define the sub and superlevel sets of $w$ in $\mK_\rd(x_0)$ and $Q_\rd(x_0,T_0,T_1)$, respectively, by
\begin{align}\label{eq:def-A}\begin{aligned}
A_{K,\rd,x_0}^{\pm}(\tau)&:=\big\{ x\in \mK_\rd(x_0)\,:\,\big(w(x,\tau)-K\big)_\pm>0\big\}\\
\text{and }\ A_{K,\rd,x_0,T_0,T_1}^\pm&:=Q_\rd(x_0,T_0,T_1)\cap\big\{(x,t)\in \OmT\,:\,\big(w(x,t)-K\big)_\pm>0\big\}.
\end{aligned}
\end{align}
Herein, $f_+(x,t):=\max\{f(x,t),0\}$ and $f_-(x,t):=\max\{-f(x,t),0\}$ denote the positive and negative part of $f$. If $x_0,T_0$ and $T_1$ are clear, we will drop these from the subscript.

\subsection{Energy estimates with level truncation}
The first localized energy estimate will be obtained from using $\Psi=\pm\big(\Phi\big([\Phi^{-1}(w)]_h\big)-K\big)_\pm\psi$ as a test function in \eqref{eq:weak-sol-steklov}, herein $\psi$ is a smooth cutoff function compactly supported in some closed parabolic cylinder.  Introducing for $\mu^+>\mu^-\geq0$ and $K>0$ the averaged truncations
\begin{align}\label{eq:Lambda}\Lambda^\pm(s)&=\Lambda^\pm(s,K):=\pm\int_{\Phi^{-1}(K)}^{\Phi^{-1}(s)}\big(\Phi(\tilde{\sigma})-K\big)_\pm\intd\tilde{\sigma},\quad s\in[\mu^-,\mu^+],\\
\label{eq:Lambda_h}
\text{and}\quad \Lambda^\pm_h(s)&=\Lambda_h^\pm(s,K):=\pm\int_{\Phi^{-1}(K)}^{[\Phi^{-1}(s)]_h}\big(\Phi(\tilde{\sigma})-K\big)_\pm\intd\tilde{\sigma},\quad s\in[\mu^-,\mu^+],
\end{align}
it will be possible to rewrite the term which contains the time-derivative first in terms of $\Lambda_h^\pm$ and $\psi$ additional manipulations thereafter establish the following lemma upon passing to the limit $h\to0$. We should note that in \eqref{eq:loc-energy} below, the terms containing $\Lambda^\pm$ are not really favorable terms yet, since they contain an implicit dependence on $\Phi^{-1}$. Obtaining sufficiently good estimates for $\Lambda^\pm$ from above and below is one important task undertaken at the beginning of Section~\ref{sec4}.

\begin{lemma}\label{lem:local-energy}
Let $w$ be a nonnegative bounded weak solution of \eqref{eq:trans-gen-eq} with $M>0$ such that \eqref{eq:linfty-bound-w} holds. There is $C=C(\Phi,M,q_1,q_2,a,\hat{b})>0$ with the following property: Assume that $x_0\in\Omega$, $0<T_0<T_1<T$ and $0<\sqrt{\dimN}\rd<\dist(x_0,\romega)$ are such that $Q_\rd=Q_\rd(x_0,T_0,T_1)\subset \OmT$. Then for each $K\in(0,2M]$, the functions $\Lambda^\pm$ as defined in \eqref{eq:Lambda} and every cut-off function $\psi\in C_0^\infty(\OmT)$ fulfilling $0\leq \psi\leq 1$ and $\psi\equiv0$ on $\Omega\setminus\overline{Q_\rd}$ as well as $\psi_{\vert\partial_p Q_\rd}=0$, satisfy
\begin{align}\label{eq:loc-energy}
&\int_{\mK_\rd}\Lambda^\pm\big(w(\cdot,T_1)\big)\psi^2(\cdot,T_1)+\frac{1}{2}\intQ \big|\nabla\big((w-K)_\pm\psi\big)\big|^2\nonumber\\
\leq\ &\intQ\Lambda^\pm(w)\big(\psi^2\big)_t+\frac{3}{2}\intQ(w-K)_\pm^2|\nabla\psi|^2+C\Big(\int_{T_0}^{T_1}|A_{K,\rd}^\pm(t)|^{\frac{q_1(q_2-1)}{q_2(q_1-1)}}\Big)^\frac{q_1-1}{q_1}.
\end{align}
\end{lemma}

\begin{bew}
We pick any cut-off function $\psi\in C_0^\infty(\OmT)$ satisfying $0\leq \psi\leq 1$ in $\OmT$ and $\psi\equiv0$ in $\OmT\setminus \overline{Q_\rd}$ with $\psi_{\vert\partial_p Q_\rd}=0$ and set $\Psi_h:=\pm\big(\Phi\big([\Phi^{-1}(w)]_h\big)-K\big)_\pm\psi^2$. Notice that by the assumptions on $\Phi$, $\Phi^{-1}$ and $w$ we may use $\Psi_h$ as test function in \eqref{eq:weak-sol-steklov}. Integrating from $T_0$ to $T_1$ and recalling $\Lambda_h^\pm$ from \eqref{eq:Lambda_h} we find that, akin to the reasoning undertaken in \cite[Proposition 3.1]{dibenedettoDegenerateParabolicEquations1993}, writing
$$\pm\iint\limits_{Q_\rd}\Big(\big[\Phi^{-1}(w)\big]_h\Big)_t\Big(\Phi\big(\big[\Phi^{-1}(w)\big]_h\big)-K\Big)_\pm\psi^2=\iint\limits_{Q_\rd}\big(\Lambda_h^\pm(w)\big)_t\psi^2,$$
integrating by parts with respect to time and taking the limit $h\to0$, we find that that with $\Lambda^\pm$ given by \eqref{eq:Lambda} we have
\begin{align}\label{eq:loc-en-+-est1}
0&=\int_{\mK_\rd}\Lambda^\pm\big(w(\cdot,T_1)\big)\psi^2(\cdot,T_1)-\intQ \Lambda^\pm(w)\big(\psi^2\big)_t\pm\intQ\nabla w\cdot\nabla\big((w-K)_\pm\psi^2)\nonumber\\
&\hspace*{3.9cm}\pm\intQ a\Phi^{-1}(w)\cdot\nabla\big((w-K)_\pm\psi^2\big)\mp\intQ b\big(\cdot,\cdot,\Phi^{-1}(w)\big)(w-K)_\pm\psi^2\nonumber\\
&=:I^\pm_1+ I^\pm_2+ I^\pm_3+ I^\pm_4+I^\pm_5.
\end{align}
Making use of $\nabla w\cdot\nabla(w-K)_\pm=\pm\big|\nabla(w-K)_\pm\big|^2$ and $(w-K)_\pm\nabla w=\pm(w-K)_\pm\nabla(w-K)_\pm$, we can rewrite $I^\pm_3$ as
\begin{align}\label{eq:loc-en-+-est3}
I^\pm_3&=\pm\intQ\nabla w\cdot\nabla(w-K)_\pm\psi^2 \pm 2\intQ(w-K)_\pm\psi\big(\nabla w\cdot\nabla\psi\big)\nonumber\\
&=\intQ \Big(\nabla\big((w-K)_\pm\psi\big)-(w-K)_\pm\nabla\psi\Big)\Big(\nabla\big((w-K)_\pm\psi\big)+(w-K)_\pm\nabla\psi\Big)\nonumber\\
&=\intQ\big|\nabla\big((w-K)_\pm\psi\big)\big|^2 - \intQ (w-K)_\pm^2|\nabla\psi|^2.
\end{align}
 For the integral $I^\pm_4$ in \eqref{eq:loc-en-+-est1} we draw on Young's inequality so that writing $A^\pm_{K,\rd}:=A^\pm_{K,\rd,x_0,T_0,T_1}$ as in \eqref{eq:def-A}, we have
\begin{align}\label{eq:loc-en-+-est4}
|I^\pm_4|&=\Big\vert\intQ a\Phi^{-1}(w)\cdot \nabla\big((w-K)_\pm\psi^2\big)\Big\vert\nonumber\\
&\leq\intQ \big\vert a\Phi^{-1}(w)\nabla\big((w-K)_\pm\psi\big)\psi\big\vert+\intQ \big\vert a\Phi^{-1}(w)(w-K)_\pm\psi\nabla\psi\big\vert\nonumber\\
&\leq\frac{1}{2}\intQ\big|\nabla\big((w-K)_\pm\psi\big)\big|^2+\frac{1}{2}\intQ(w-K)_\pm^2|\nabla\psi|^2+\iint\limits_{A^\pm_{K,\rd}}|a|^2\big|\Phi^{-1}(w)\big|^2\psi^2\nonumber\\
&\leq\frac{1}{2}\intQ\big|\nabla\big((w-K)_\pm\psi\big)\big|^2+\frac{1}{2}\intQ(w-K)_\pm^2|\nabla\psi|^2+\big(\Phi^{-1}(M)\big)^2\iint\limits_{A^\pm_{K,\rd}}|a|^2,
\end{align}
due to $\Phi^{-1}$ being monotonically increasing, $0\leq\psi\leq1$ and $0\leq w\leq M$. Moreover, we conclude from $0\leq\psi\leq 1$, $(w-K)_\pm\leq 2M$ and $|b(x,t,\xi)|\leq \hat{b}(x,t)$ for a.e. $(x,t,\xi)\in\OmT\times\R$, that $I^\pm_5$ satisfies
\begin{align}\label{eq:loc-en-+-est5}
|I^\pm_5|=\Big\vert\intQ b\big(\cdot,\cdot,\Phi^{-1}(w)\big)(w-K)_\pm\psi^2\Big\vert\leq 2M\iint\limits_{A^\pm_{K,\rd}}|\hat{b}|.
\end{align}
Therefore, collecting \eqref{eq:loc-en-+-est1}--\eqref{eq:loc-en-+-est5} shows that
\begin{align*}
&\int_{\mK_\rd}\Lambda^\pm\big(w(\cdot,T_1)\big)\psi^2(\cdot,T_1)+\frac{1}{2}\intQ\big|\nabla\big((w-K)_\pm\psi\big)\big|^2\\
\leq\ &\!\intQ\Lambda^\pm(w)\big(\psi^2\big)_t+\frac{3}{2}\!\intQ(w-K)_\pm^2|\nabla\psi|^2+\big(\Phi^{-1}(M)\big)^2\!\iint\limits_{A^\pm_{K,\rd}}|a|^2+2M\iint\limits_{A^\pm_{K,\rd}}|\hat{b}|.\nonumber
\end{align*}
To further estimate the last two integrals on the right hand side, we make use of two applications of Hölder's inequality to see that for $q_1,q_2>1$ we have
\begin{align*}
\iint\limits_{A^\pm_{K,\rd}}|a|^2&\leq \int_{T_0}^{T_1}\!\Big(\int_{A^\pm_{K,\rd}(t)}\!\!|a(\cdot,t)|^{2q_2}\Big)^\frac{1}{q_2}\Big(\int_{A^\pm_{K,\rd}(t)}\!\! 1\Big)^{\frac{q_2-1}{q_2}}\\
&\leq\int_{T_0}^{T_1}\|a(\cdot,t)\|_{\LSp{2q_2}{\mK_\rd}}^2\big|A_{K,\rd}^\pm(t)\big|^{\frac{q_2-1}{q_2}}\\
&\leq \Big(\int_{T_0}^{T_1}\!\|a(\cdot,t)\|_{\LSp{2q_2}{\mK_\rd}}^{2q_1}\Big)^\frac{1}{q_1}\Big(\int_{T_0}^{T_1}\!\big|A_{K,\rd}^\pm(t)\big|^{\frac{q_1(q_2-1)}{q_2(q_1-1)}}\Big)^\frac{q_1-1}{q_1}\\
&\leq\|a\|_{\LSp{2q_1}{(0,T);\LSpn{2q_2}{\Omega;\R^n}}}^2\Big(\int_{T_0}^{T_1}\!\big|A^\pm_{K,\rd}(t)\big|^{\frac{q_1(q_2-1)}{q_2(q_1-1)}}\Big)^\frac{q_1-1}{q_1}
\end{align*}
and
\begin{align*}
\iint\limits_{A^\pm_{K,\rd}}|\hat{b}|\leq\|\hat{b}\|_{\LSp{q_1}{(0,T);\LSp{q_2}{\Omega}}}\Big(\int_{T_0}^{T_1}\!\big|A^\pm_{K,\rd}(t)\big|^{\frac{q_1(q_2-1)}{q_2(q_1-1)}}\Big)^\frac{q_1-1}{q_1},
\end{align*}
which proves \eqref{eq:loc-energy} with constant $C=C(\Phi,M,q_1,q_2,a,\hat{b})>0$ explicitly given by 
\[
C:=\big(\Phi^{-1}(M)\big)^2\|a\|^2_{\LSpn{2q_1}{(0,T);\LSp{2q_2}{\Omega;\R^\dimN}}}+2M\|\hat{b}\|_{\LSp{q_1}{(0,T);\LSp{q_2}{\Omega}}}.
\qedhere\]
\end{bew}

\subsection{An estimate of logarithmic type}
For a second testing procedure we introduce the following logarithmic function. Given $L>0$, $K\in(0,L]$ and $\delta\in(0,1)$ we set
\begin{align}\label{eq:xi-def}
\xi(s):=\xi_{L,K,\delta}(s):=\Big(\ln\Big(\frac{K}{(1+\delta)K-(s-L+K)_+}\Big)\Big)_+,\quad s\in[0,L].
\end{align}
Clearly, $\xi$ satisfies $0\leq \xi\leq \ln\big(\frac{1}{\delta}\big)$ on $[0,L]$ with $\xi\equiv0$ on $[0,L-(1-\delta)K]$ as well as the properties
\begin{align}\label{eq:xi-prop}
0\leq \xi' \leq \frac{1}{\delta K},\qquad \big(\xi'\big)^2=\xi{''},\quad\text{and}\quad \big(\xi^2\big){''}=2(1+\xi)\big(\xi'\big)^2\quad\text{on }(0,L).
\end{align}
Moreover, given a cylinder $Q_\rd=\mK_\rd(x_0)\times(T_0,T_1]\subset\OmT$ we pick any smooth spatial cutoff function $\zeta=\zeta_{x_0,\rd}\in C_0^\infty(\Omega)$ compactly supported in $\mK_\rd=\mK_\rd(x_0)$ satisfying $0\leq \zeta\leq 1$. In the standard setting the test function for the logarithmic estimate would now be of the form $\Psi(w)=\big(\xi^2(w)\big)'\zeta^2$ (\cite{dibenedettoDegenerateParabolicEquations1993}). For our application, however, this would not be sufficient to treat the time-derivative term. 
One obvious approach to cope with this problem would be to take $\Psi\big(\Phi^{-1}(w)\big)$ instead of just $\Psi(w)$. Then, however, the integral originating from the convection portion of the equation would contain a derivative of $\Phi^{-1}$, which cannot be treated in an obvious manner when $w$ is close to zero. Instead, we use a function of the form
\begin{align}\label{eq:tpsi-def}
\Psi(w):=\Psi_{\Phi,L,K,\delta,\rd}\big(w(x,t)\big):=\Phi'\big(\Phi^{-1}\big(w(x,t)\big)\big)\big(\xi^2\big)'\big(w(x,t)\big)\zeta^2(x),\quad (x,t)\in \OmT
\end{align}
as test function in \eqref{eq:weak-sol-eq}. From the properties of $\Phi$, $\xi$ and $\zeta$ we conclude that $\Psi(w)\in W_0^{1,2}(Q_\rd)$ and direct calculations involving the inversion rule and the identity $(\xi^2)''=2(1+\xi)(\xi')^2$ show that on $\big\{(x,t)\in Q_\rd\,\vert\,w(x,t)>L-(1-\delta)K\big\}$ we have
\begin{align}\label{eq:tpsi-grad}
\nabla\Psi(w)=2\frac{\Phi''\big(\Phi^{-1}(w)\big)}{\Phi'\big(\Phi^{-1}(w)\big)}\xi(w)\xi'(w)\zeta^2\nabla w&+2\Phi'\big(\Phi^{-1}(w)\big)\big(1+\xi(w)\big)\big(\xi'(w)\big)^2\zeta^2\nabla w\nonumber\\&+4\Phi'\big(\Phi^{-1}(w)\big)\xi(w)\xi'(w)\zeta\nabla\zeta.
\end{align}
As one can see, the derivative of $\Phi^{-1}$ is still present in one of the terms of $\nabla\Psi(w)$. Nonetheless, this time, because of the assumed convexity of $\Phi$, we also gain an additional well-signed and beneficial term from partial integration in the diffusion integral, which we can use to consume parts of the convection term in such a way that the assumption in \eqref{eq:cond-phi} takes care of the troublesome portions, allowing us to derive the following estimate.

\begin{lemma}\label{lem:log-est}
Let $w$ be a nonnegative bounded weak solution of \eqref{eq:trans-gen-eq} with $M>0$ such that \eqref{eq:linfty-bound-w} holds. Then there is $C=C(\Phi,M,q_1,q_2,a,\hat{b})>0$ with the following property: Assume that $x_0\in\Omega$, $0<T_0<T_1<T$ as well as $0<\sqrt{\dimN}\rd<\dist(x_0,\romega)$ are such that $Q_\rd=Q_\rd(x_0,T_0,T_1)\subset\OmT$. Suppose that $L>0$ is such that $L\geq \esssup_{Q_\rd} w$ and $K\in(0,L]$. Then, given any $\zeta=\zeta_{x_0,\rd}\in C_0^{\infty}(\Omega)$ compactly supported in $\mK_\rd=\mK_\rd(x_0)$ with $0\leq \zeta\leq 1$ and any $\delta\in(0,1)$, the function $\xi=\xi_{L,K,\delta}$ provided by \eqref{eq:xi-def} satisfies
\begin{align}\label{eq:log-est}
\int_{\mK_\rd}\xi^2\big(w(\cdot,T_1)\big)\zeta^2\nonumber&\leq\int_{\mK_\rd}\xi^2\big(w(\cdot,T_0)\big)\zeta^2+6\ln\Big(\frac{1}{\delta}\Big)\iint\limits_{Q_\rd}\Phi'\big(\Phi^{-1}(w)\big)|\nabla\zeta|^2\nonumber\\&\qquad+C\Big(\frac{1+\ln\big(\frac{1}{\delta}\big)}{\delta^2 K^2}+\frac{\ln\big(\frac1\delta\big)}{\delta K}\Big)\Big(\int_{T_0}^{T_1}\!\big|A^+_{L-K,\rd}(t)\big|^{\frac{q_1(q_2-1)}{q_2(q_1-1)}}\Big)^\frac{q_1-1}{q_1}.
\end{align}
\end{lemma}

\begin{bew}
Again, first arguing formally using Steklov averages, i.e. taking $$\Psi_h(w):=\Phi'\big(\big[\Phi^{-1}(w)\big]_h\big)\big(\xi^2\big)'\big(\Phi\big([\Phi^{-1}(w)]_h\big)\big)\zeta^2$$ as test function in \eqref{eq:weak-sol-steklov} and going to the limit, we find that with $\Psi(w)$ as described in \eqref{eq:tpsi-def}, we have
\begin{align}\label{eq:log-est-0}
\int_{\mK_\rd}\xi^2\big(w(\cdot,T_0)\big)\zeta^2=\int_{\mK_\rd}\xi^2\big(w(\cdot,T_1)\big)\zeta^2&+\intQ\nabla w\cdot\nabla \Psi(w)+\intQ a\Phi^{-1}(w)\cdot\nabla\Psi(w)\nonumber\\&-\intQ b\big(\cdot,\cdot,\Phi^{-1}(w)\big)\Psi(w)=:\tilde{I}_1+\tilde{I}_2+\tilde{I}_3+\tilde{I}_4.
\end{align}
To estimate $\tilde{I}_2$ we draw on \eqref{eq:tpsi-grad}, the nonnegativity of $\Phi',\Phi^{-1},\zeta,\xi$ and $\xi'$ as well as Young's inequality to find that for all $\eta_1>0$
\begin{align}\label{eq:log-est-2}
\tilde{I}_2&=2\intQ\!\frac{\Phi''\big(\Phi^{-1}(w)\big)}{\Phi'\big(\Phi^{-1}(w)\big)}\xi(w)\xi'(w)\zeta^2|\nabla w|^2+2\intQ\!\Phi'\big(\Phi^{-1}(w)\big)\big(1+\xi(w)\big)\big(\xi'(w)\big)^2\zeta^2|\nabla w|^2\nonumber\\
&\hspace*{6.3cm}+4\intQ\!\Phi'\big(\Phi^{-1}(w)\big)\xi(w)\xi'(w)\zeta(\nabla\zeta\cdot\nabla w)\nonumber\\
&\geq 2\intQ\!\frac{\Phi''\big(\Phi^{-1}(w)\big)}{\Phi'\big(\Phi^{-1}(w)\big)}\xi(w)\xi'(w)\zeta^2|\nabla w|^2+2\intQ\!\Phi'\big(\Phi^{-1}(w)\big)\big(1+(1-\eta_1)\xi(w)\big)\big(\xi'(w)\big)^2\zeta^2|\nabla w|^2\nonumber\\&\hspace{6.3cm}-\frac{2}{\eta_1}\!\intQ\Phi'\big(\Phi^{-1}(w)\big)\xi(w)|\nabla\zeta|^2.
\end{align}
Similarly, relying on the nonnegativity of $\Phi''$ as entailed by the assumed convexity and three applications of Young's inequality, we estimate $\tilde{I}_3$ to obtain that for all $\eta_2,\eta_3,\eta_4>0$ we have
\begin{align}\label{eq:log-est-3}
\tilde{I}_3&=2\intQ\!\Phi^{-1}(w)\frac{\Phi''\big(\Phi^{-1}(w)\big)}{\Phi'\big(\Phi^{-1}(w)\big)}\xi(w)\xi'(w)\zeta^2(a\cdot\nabla w)\nonumber\\
&\hspace*{2.1cm}+2\intQ\!\Phi^{-1}(w)\Phi'\big(\Phi^{-1}(w)\big)\big(1+\xi(w)\big)\big(\xi'(w)\big)^2\zeta^2(a\cdot\nabla w)\nonumber\\
&\hspace*{4.3cm}+4\intQ\!\Phi^{-1}(w)\Phi'\big(\Phi^{-1}(w)\big)\xi(w)\xi'(w)\zeta(a\cdot\nabla\zeta)\nonumber\\
&\geq -2\eta_2\intQ\!\frac{\Phi''\big(\Phi^{-1}(w)\big)}{\Phi'\big(\Phi^{-1}(w)\big)}\xi(w)\xi'(w)\zeta^2|\nabla w|^2-\frac{1}{2\eta_2}\intQ\!\frac{\Phi''\big(\Phi^{-1}(w)\big)\big(\Phi^{-1}(w)\big)^2}{\Phi'\big(\Phi^{-1}(w)\big)}\xi(w)\xi'(w)\zeta^2|a|^2\nonumber\\
&\hspace*{0.8cm}-2\eta_3\intQ\!\Phi'\big(\Phi^{-1}(w)\big)\big(1+\xi(w)\big)\big(\xi'(w)\big)^2\zeta^2|\nabla w|^2-8\eta_4\intQ\!\Phi'\big(\Phi^{-1}(w)\big)\xi(w)|\nabla\zeta|^2\nonumber\\
&\hspace*{1.6cm}-\intQ\!\big(\Phi^{-1}(w)\big)^2\Phi'\big(\Phi^{-1}(w)\big)\Big(\frac{1}{2\eta_3}+\Big(\frac{1}{2\eta_3}+\frac{1}{2\eta_4}\Big)\xi(w)\Big)\big(\xi'(w)\big)^2\zeta^2|a|^2.
\end{align}
Moreover, we clearly have
\begin{align}\label{eq:log-est-4}
\tilde{I}_4\geq -2\intQ\big|b(\cdot,\cdot,\Phi^{-1}(w))\big|\Phi'\big(\Phi^{-1}(w)\big)\xi(w)\xi'(w)\zeta^2,
\end{align}
so that by combining \eqref{eq:log-est-0}--\eqref{eq:log-est-4} we conclude that for all $\eta_1,\eta_2,\eta_3,\eta_4>0$
\begin{align*}
&\int_{\mK_\rd}\xi^2\big(w(\cdot,T_1)\big)\zeta^2+2(1-\eta_2)\intQ\!\frac{\Phi''\big(\Phi^{-1}(w)\big)}{\Phi'\big(\Phi^{-1}(w)\big)}\xi(w)\xi'(w)\zeta^2|\nabla w|^2\\
&\hspace*{3.1cm}+2\intQ\!\Phi'\big(\Phi^{-1}(w)\big)\big(1-\eta_3+(1-\eta_1-\eta_3)\xi(w)\big)\big(\xi'(w)\big)^2\zeta^2|\nabla w|^2\\
\leq\ &\int_{\mK_\rd}\xi^2\big(w(\cdot,T_0)\big)\zeta^2+\Big(\frac{2}{\eta_1}+8\eta_4\Big)\!\intQ\Phi'\big(\Phi^{-1}(w)\big)\xi(w)|\nabla\zeta|^2
\\&\hspace*{0.6cm}
+\frac{1}{2\eta_2}\intQ\!\frac{\Phi''\big(\Phi^{-1}(w)\big)\big(\Phi^{-1}(w)\big)^2}{\Phi'\big(\Phi^{-1}(w)\big)}\xi(w)\xi'(w)\zeta^2|a|^2\\
&\hspace*{1.2cm}+\intQ\!\big(\Phi^{-1}(w)\big)^2\Phi'\big(\Phi^{-1}(w)\big)\Big(\frac{1}{2\eta_3}+\big(\frac{1}{2\eta_3}+\frac{1}{2\eta_4}\big)\xi(w)\Big)\big(\xi'(w)\big)^2\zeta^2|a|^2\\
&\hspace*{1.8cm}+2\intQ\big|b(\cdot,\cdot,\Phi^{-1}(w))\big|\Phi'\big(\Phi^{-1}(w)\big)\xi(w)\xi'(w)\zeta^2.
\end{align*}
Letting $\eta_1=\frac{1}{2}$, $\eta_2=1$, $\eta_3=\frac12$ and $\eta_4=\frac14$ this entails that
\begin{align}\label{eq:log-est-5}
&\int_{\mK_\rd}\xi^2\big(w(\cdot,T_1)\big)\zeta^2\nonumber\\\leq\ &\int_{\mK_\rd}\xi^2\big(w(\cdot,T_0)\big)\zeta^2+6\intQ\!\Phi'\big(\Phi^{-1}(w)\big)\xi(w)|\nabla\zeta|^2\nonumber
+\frac{1}{2}\intQ\!\frac{\Phi''\big(\Phi^{-1}(w)\big)\big(\Phi^{-1}(w)\big)^2}{\Phi'\big(\Phi^{-1}(w)\big)}\xi(w)\xi'(w)\zeta^2|a|^2\nonumber
\\&\hspace*{0.6cm}
+3\intQ\!\big(\Phi^{-1}(w)\big)^2\Phi'\big(\Phi^{-1}(w)\big)\big(1+\xi(w)\big)\big(\xi'(w)\big)^2\zeta^2|a|^2\nonumber
\\&\hspace*{1.2cm}
+2\intQ\!\big|b(\cdot,\cdot,\Phi^{-1}(w))\big|\Phi'\big(\Phi^{-1}(w)\big)\xi(w)\xi'(w)\zeta^2.
\end{align}
To further treat the last four integrals on the right-hand side and establish a connection to suitable superlevel sets, we observe that by our construction of $\xi$ we have
\begin{align}\label{eq:log-est-7}
\xi(w)=0\text{ on }\big\{w\leq L-K\big\}\subseteq \big\{w\leq L-(1-\delta)K\big\}.
\end{align}
Furthermore, we note that since $w\leq M$ the fact that both $\Phi^{-1}$ and $\Phi'$ are monotonically increasing entails that 
\begin{align}\label{eq:log-est-5.1}
\Phi'\big(\Phi^{-1}(w)\big)\leq \Phi'\big(\Phi^{-1}(M)\big)\quad\text{and}\quad\big(\Phi^{-1}(w)\big)^2\Phi'\big(\Phi^{-1}(w)\big)\leq \big(\Phi^{-1}(M)\big)^2\Phi'\big(\Phi^{-1}(M)\big)
\end{align}
and from Lemma \ref{lem:phi} we find that since $\Phi^{-1}(w)\leq \Phi^{-1}(M)$ on $Q_\rd$ there is $C_1=C_1(\Phi,M)$ such that $\Phi''(\Phi^{-1}(w))\Phi^{-1}(w)\leq C_1\Phi'(\Phi^{-1}(w))$, so that we may estimate
\begin{align}\label{eq:log-est-5.2}
\frac{\Phi''\big(\Phi^{-1}(w)\big)\big(\Phi^{-1}(w)\big)^2}{\Phi'\big(\Phi^{-1}(w)\big)}\leq  C_1\Phi^{-1}(M)=:C_2(\Phi,M)\quad\text{on }Q_\rd.
\end{align}
Hence, aggregating \eqref{eq:log-est-5}--\eqref{eq:log-est-5.2} and drawing on the estimates $\xi\leq\ln(\tfrac{1}{\delta})$, $\xi'\leq\frac{1}{\delta K}$ of \eqref{eq:xi-prop} as well as the properties of $\zeta$, we obtain $C_3=C_3(\Phi,M)>0$ satisfying
\begin{align*}
\int_{\mK_\rd}\xi^2\big(w(\cdot,T_1)\big)\zeta^2\leq&\int_{\mK_\rd}\xi^2\big(w(\cdot,T_0)\big)\zeta^2+6\ln\left(\tfrac{1}{\delta}\right)\iint\limits_{Q_\rd}\Phi'\big(\Phi^{-1}(w)\big)|\nabla\zeta|^2\\
&+C_3\left(\frac{\ln(\frac1\delta)}{\delta K}+\frac{\left(1+\ln(\frac1\delta)\right)}{\delta^2K^2}\right)\iint\limits_{Q_\rd\cap\{w>L-K\}}\!\!|a|^2+C_3\frac{\ln(\frac1\delta)}{\delta K}\iint\limits_{Q_\rd\cap\{w>L-K\}}\!\!|\hat{b}|.
\end{align*}
Herein, we have in view of Hölder's inequality that
\begin{align*}
&\iint\limits_{Q_\rd\cap\{w>L-K\}}\!\!|a|^2\leq \|a\|_{\LSpn{2q_1}{(0,T);\LSp{2q_2}{\Omega;\R^\dimN}}}\Big(\int_{T_0}^{T_1}\big|A_{L-K,\rd}^+(t)\big|^\frac{q_1(q_2-1)}{q_2(q_1-1)}\Big)^\frac{q_1-1}{q_1}
\intertext{and}
&\iint\limits_{Q_\rd\cap\{w>L-K\}}\!\!|\hat{b}|\leq \|\hat{b}\|_{\LSp{q_1}{(0,T);\Lo[q_2]}}\Big(\int_{T_0}^{T_1}\big|A_{L-K,\rd}^+(t)\big|^{\frac{q_1(q_2-1)}{q_2(q_1-1)}}\Big)^\frac{q_1-1}{q_1},
\end{align*}
so that letting $C=C(\Phi,M,q_1,q_2,a,\hat{b}):=C_3\big(\|a\|_{\LSpn{2q_1}{(0,T);\LSp{2q_2}{\Omega;\R^\dimN}}}+\|\hat{b}\|_{\LSp{q_1}{(0,T);\Lo[q_2]}}\big)>0$ yields \eqref{eq:log-est}, completing the proof.
\end{bew}

One minor caveat from our adjusted test function $\Psi$ is, that -- in contrast to the standard procedure (e.g. \cite{dibenedettoDegenerateParabolicEquations1993,BenedettoGianazzaVespri-Harnack}) -- an inequality with a truncation using the negative part is not possible, as then $\xi'$ would have a negative sign, turning the essentially beneficial term from before into a problematic term. For the most part this is inconsequential, but it will influence the approach when discussing regularity up to the boundary in Section~\ref{sec6:boundary}.

\setcounter{equation}{0} 
\section{Tracking the decay of essential oscillation within nested cylinders}\label{sec4}
The most vital elements when establishing the decay of oscillation are the two so-called alternatives, which -- in essence -- when starting from some cylinder  provide quantitative information on the size of the sets inside a subcylinder where $w$ is close to its essential infimum or essential supremum over the starting cylinder. In fact, the first alternative (cf. Lemma~\ref{lem:first-alt-1}) says that if the set where $w$ is significantly larger than its essential infimum is not too big then either the radius of the cylinder is of a certain size, or there is a subcylinder where $w$ is bounded away from the infimum of the outer cylinder. Similar information can be derived for the set where $w$ is close to the essential supremum via the second alternative (see Lemma~\ref{lem:second-alt}). These results are obtained by iterating our local energy estimates from Section~\ref{sec3} for suitable choices for the radius $\rd$ of the cylinders and levels $K$. One very essential ingredient for the inequalities with different levels and radii to interlock is adequate control over $\Lambda^\pm$, which we prepare in the following Lemma.

\begin{lemma}\label{lem:Lambda}
Let $\mu^+>\mu^-\geq0$. For each $K>0$, the functions 
\begin{align*}
\Lambda^\pm(s)=\Lambda^\pm(s,K)=\pm\int_{\Phi^{-1}(K)}^{\Phi^{-1}(s)}\big(\Phi(\tilde{\sigma})-K\big)_\pm\intd\tilde{\sigma},\quad s\in[\mu^-,\mu^+],
\end{align*}
as defined in \eqref{eq:Lambda}, satisfy
\begin{align}\label{eq:Lambda-est-+}
\frac{1}{2\Phi'\big(\Phi^{-1}(\mu^+)\big)}(s-K)_+^2\leq\Lambda^+(s)\leq\frac{1}{2\Phi'\big(\Phi^{-1}(K)\big)}(s-K)_+^2\quad\text{for all }s\in[\mu^-,\mu^+]
\end{align}
and
\begin{align}\label{eq:Lambda-est--}
\frac{1}{2\Phi'\big(\Phi^{-1}(K)\big)}(s-K)_-^2\leq\Lambda^-(s)\leq\frac{\Phi^{-1}(K)}{K}(s-K)_-^2\quad\text{for all }s\in[\mu^-,\mu^+].
\end{align}
If, moreover, $\mu^- >0$, then $\Lambda^-(s)$ additionally satisfies
\begin{align}\label{eq:Lambda-est--_mu-}
\Lambda^-(s)\leq \frac{1}{2\Phi'\big(\Phi^{-1}(\mu^-)\big)}(s-K)_-^2\quad\text{for all }s\in[\mu^-,\mu^+].
\end{align}
\end{lemma}

\begin{bew}
Noticing that for $s\leq K$ we have $\Lambda^+(s)=0$ and $(s-K)_+=0$, the inequalities in \eqref{eq:Lambda-est-+} are obviously valid for $s\leq K$, so we are left with verifying them for $s\in(K,\mu^+]$. Similarly, for $s\geq K$ we have $\Lambda^-(s)=0$ and $(s-K)_-=0$, so clearly \eqref{eq:Lambda-est--} and \eqref{eq:Lambda-est--_mu-} hold for $s\geq K$ and we only have to consider \eqref{eq:Lambda-est--} and \eqref{eq:Lambda-est--_mu-} for $s\in[\mu^-,K)$.

By change of variables $\sigma=\Phi(\tilde{\sigma})$, we find that
\begin{align*}
\Lambda^+(s)=\int_{K}^{s}\!\big(\sigma-K\big)_+\big(\Phi^{-1}\big)'(\sigma)\intd\sigma\quad\text{for }s\in(K,\mu^+]
\end{align*}
and
\begin{align*}
\Lambda^-(s)=\int_{s}^{K}\!\big(\sigma-K\big)_-\big(\Phi^{-1}\big)'(\sigma)\intd\sigma\quad\text{for }s\in[\mu^-,K).
\end{align*}
Due to $\Phi$ being convex, $\Phi'$ is monotonically increasing and hence we conclude from the monotonically increasing property of $\Phi^{-1}$ and the inversion rule that $$\frac{1}{\Phi'\big(\Phi^{-1}(\mu^+)\big)}\leq \big(\Phi^{-1}\big)'(\sigma)\leq\frac{1}{\Phi'\big(\Phi^{-1}(K)\big)}\quad\text{ for all }\sigma\in(K,\mu^+].$$ Accordingly, we find that
\begin{align*}
\Lambda^+(s)\geq\frac{1}{\Phi'\big(\Phi^{-1}(\mu^+)\big)}\int_{K}^{s}(\sigma-K)_+\intd\sigma=\frac{(s-K)_+^2}{2\Phi'\big(\Phi^{-1}(\mu^+)\big)}\quad\text{for all }s\in(K,\mu^+]
\end{align*}
as well as
\begin{align*}
\Lambda^+(s)\leq \frac{1}{\Phi'\big(\Phi^{-1}(K)\big)}\int_K^s(\sigma-K)_+\intd\sigma=\frac{(s-K)_+^2}{2\Phi^{-1}\big(\Phi(K)\big)}\quad\text{for all }s\in(K,\mu^+].
\end{align*}
In a similar fashion, we conclude from $\big(\Phi^{-1}\big)'(\sigma)\geq\frac{1}{\Phi'(\Phi^{-1}(K))}$ for all $\sigma\in[\mu^-,K)$ that
\begin{align*}
\Lambda^-(s)\geq \frac{(s-K)_-^2}{2\Phi'\big(\Phi^{-1}(K)\big)}\quad\text{for all }s\in[\mu^-,K).
\end{align*}
Let us turn to the upper bound on $\Lambda^-$ in \eqref{eq:Lambda-est--}. Here, we note that $(\sigma-K)_-\leq (s-K)_-$ for all $\sigma\in[s,K)$ and thus
\begin{align}\label{eq:proof-lambda-1}
\Lambda^-(s)&\leq (s-K)_-\int_s^{K}\big(\Phi^{-1}\big)'(\sigma)\intd\sigma=(s-K)_-\Big(\Phi^{-1}(K)-\Phi^{-1}(s)\Big).
\end{align}
In view of $\Phi(0)=0$, we can make use of the convexity of $\Phi$ and the fact that $\Phi^{-1}$ is monotonically increasing to estimate
\begin{align*}
\frac{K}{\Phi^{-1}(K)}=\frac{\Phi\big(\Phi^{-1}(K)\big)-\Phi\big(\Phi^{-1}(0)\big)}{\Phi^{-1}(K)-\Phi^{-1}(0)}\leq \frac{\Phi\big(\Phi^{-1}(K)\big)-\Phi\big(\Phi^{-1}(s)\big)}{\Phi^{-1}(K)-\Phi^{-1}(s)}=\frac{K-s}{\Phi^{-1}(K)-\Phi^{-1}(s)}
\end{align*}
for $K> s\geq \mu^-\geq 0$, so that actually
\begin{align*}
\Phi^{-1}(K)-\Phi^{-1}(s)\leq \frac{\Phi^{-1}(K)}{K}(K-s)\quad\text{for }\mu^-\leq s< K.
\end{align*}
Plugging this into \eqref{eq:proof-lambda-1}, we obtain for $K> s\geq \mu^-$ that
\[
\Lambda^-(s)\leq \frac{\Phi^{-1}(K)}{K}(s-K)_-(K-s)= \frac{\Phi^{-1}(K)}{K}(s-K)_-^2.
\]
The estimate \eqref{eq:Lambda-est--_mu-} in the special case when $\mu^->0$ follows from arguments along the line of the estimate from above for $\Lambda^+$, utilizing that in this case we have $0<\Phi'\big(\Phi^{-1}(\mu^-)\big)\leq \Phi'(\Phi^{-1}(s))$ for all $s\geq \mu^->0$ again by the monotonicity properties of $\Phi$ and $\Phi^{-1}$.
\end{bew}

\subsection{The first alternative}\label{sec4:alt1}
In what follows, we always denote by $Q_\rho$ the cylindrical box of radius $\rho$ spatially centered at $x_0\in\Omega$, with top at $t_0\in(0,T)$ and height $\rho^2$, to be precise $Q_\rho:=Q_\rho(x_0,t_0):=\mK_{\rho}(x_0)\times(t_0-\rho^2,t_0]$. Moreover, we will also rely on a time-rescaling of the cylinder, leading to its so-called intrinsic geometry. Rescaled cylinders will be denoted by $$Q_\rho^\theta:=Q_\rho^\theta(x_0,t_0):=\mK_{\rho}(x_0)\times(t_0-\theta\rho^2,t_0],$$ where $\theta$ is a conveniently chosen number depending on the level sets under consideration. Additionally, to shorten notation, we set $$\theta_c:=\big(\Phi^{-1}\big)'(c)=\big(\Phi'\big(\Phi^{-1}(c)\big)\big)^{-1}\quad\text{ for }c>0.$$ Note that by the assumed monotonicity properties of $\Phi'$ and $\Phi^{-1}$ the map $c\mapsto\theta_c$ is monotonically decreasing. The symbols $\mu^-$ and $\mu^+$ will always denote the essential infimum and essential supremum, respectively, of $w$ in some cylinder. Note that without loss of generality we may always assume that $\mu^-<\mu^+$, as otherwise $w$ would be constant throughout the cylinders and hence certainly obey the decay of oscillation captured in Corollary~\ref{cor:iteration}. For convenience of notation also recall that according to \eqref{eq:standing_assump} and \eqref{eq:cond-kappa} we have $2-\frac{2}{q_1}-\frac{\dimN}{q_2}=\dimN\kappa$. The first alternative now reads as follows:

\begin{lemma}\label{lem:first-alt-1}
Let $M>0$. There is $\sig_0=\sig_0(\Phi,M,q_1,q_2,a,\hat{b})\in(0,\tfrac12)$ such that for all nonnegative bounded weak solutions $w$ of \eqref{eq:trans-gen-eq} satisfying \eqref{eq:linfty-bound-w} the following holds: Given $x_0\in\Omega$, $0<t_0<T$, $0<\vartheta$ and $0<\rd<\min\big\{\frac{1}{\sqrt{\dimN}}\dist(x_0,\partial\Omega),\sqrt{\tfrac{t_0}{4\vartheta}},1\big\}$ set $\mu^-:=\essinf_{Q_{2\rd}^{\vartheta}}w<\esssup_{Q_{2\rd}^{\vartheta}}w=:\mu^+$. Let $k\in(0,M]$ and introduce $\theta:=\min\big\{\theta_k,\theta_{\mu^-}\big\}$. Then, if $\vartheta\geq\theta$ and
\begin{align}\label{eq:first-alt-1-ass1}
\big|\big\{(x,t)\in Q_{\rd}^{\theta}:(w- \mu^--k)_->0\big\}\big|\leq \sig_0 \big|Q_{\rd}^{\theta}\big|
\end{align}
then either
\begin{align}\label{eq:first-alt-1-r}
k^{2}\theta^{-\frac{q_1-1}{q_1}}\leq \rd^{\dimN\kappa}
\end{align}
or
\begin{align}\label{eq:first-alt-concl}
\Big|\Big\{(x,t)\in Q_{\frac{\rd}{2}}^{\theta}:\big(w-\mu^--\tfrac{k}{2}\big)_->0\Big\}\Big|=0.
\end{align}
\end{lemma}

\begin{bew}
Setting $\hat{\mu}_-:=\essinf_{Q_{2\rd}^{\theta}}w$ and $\esssup_{Q_{2\rd}^{\theta}}w=:\hat{\mu}_+$, we first note that the assumption $\vartheta\geq\theta$ entails that
\begin{align*}
Q_{2r}^\theta\subseteq Q_{2r}^{\vartheta}\quad\text{and hence}\quad \mu^-\leq\hat{\mu}_-\leq\hat{\mu}_+\leq\mu^+.
\end{align*}
Our aim is to employ the geometric convergence Lemma~\ref{lem:fast-geo-conv} for a suitable sequence of sublevel sets of nested cylinders. To obtain the recursive inequality required by said lemma, we will use the local energy estimate presented in Lemma~\ref{lem:local-energy} for $\Lambda^-$. Let us start by preparing the sequence of levels as well as the nested cylinders. Set $\rd_0=\rd$ and $k_0=\mu^-+k>0$ and for $i\in\N$ define
$$\rd_i=\frac{\rd}{2}+\frac{\rd}{2^{i+1}}\quad\text{and}\quad k_i=\mu^-+\frac{k}{2}+\frac{k}{2^{i+1}}.$$
Evidently, for all $i\in\N_0$ we have $\rd_i\in(\frac{\rd}{2},\rd]$ and $k_i\in(\mu^-+\frac{k}{2},\mu^-+k]$ and 
\begin{align}\label{eq:first-alt-1-rinfty}
\rd_\infty:=\lim_{i\to\infty}\rd_i=\frac{\rd}{2}\quad\text{and}\quad k_\infty:=\lim_{i\to\infty} k_i=\mu^-+\frac{k}{2}.
\end{align}
We introduce $Q_i:=Q^{\theta}_{\rd_i}=\mK_{\rd_i}\times(t_0-\theta\rd_i^2,t_0]$ and fix a sequence of smooth cutoff functions $\psi_i\in C^{\infty}_0\big(\OmT\big)$ such that 
\begin{align}\label{eq:first-alt-1-psi-est}
&\psi_i\equiv 1\ \text{ in }\ Q_{i+1},\quad \psi_i\equiv 0\ \text{ in }\ \OmT\setminus \overline{Q_i},\ \text{ and }\ \psi_{\vert\partial_p Q_i}=0\nonumber\\&\text{with }\ |\nabla\psi_i|\leq \frac{2^{i+3}}{\rd},\ \text{ as well as }\ \big|(\psi_{i})_t\big|\leq\frac{2^{2(i+2)}}{\theta\rd^2}\ \text{ on }Q_i.
\end{align}
Since our estimations for $$J_i:=\sup_{t_0-\theta\rd_i^2\leq t\leq t_0}\int_{\mK_{\rd_i}}\big(w(\cdot,t)-k_i\big)_-^2\psi^2_i(\cdot,t)+\frac{1}{\theta_{k_i}}\intQ[i]\big|\nabla\big((w-k_i)_-\psi_i\big)\big|^2$$ differ slightly in each of the cases for $\theta$, let us briefly split our argument up in order to verify that the specific choice for $\theta$ does not impact the constants appearing below.

First we consider the instance where $\theta=\theta_k$. We recall that according to Lemma~\ref{lem:Lambda} for each $i\in\N_0$ we have $$\Lambda^-_i(w):=\Lambda^-(w,k_i)\leq\frac{\Phi^{-1}(k_i)}{k_i}(w-k_i)_-^2\quad\text{ and }\quad\Lambda_i^-(w)\geq \frac{1}{2\Phi'(\Phi^{-1}(k_i))}(w-k_i)_-^2$$
a.e. in $Q_i$.
Hence, we find from employing Lemma~\ref{lem:local-energy} for $\Lambda_i^-(w)$ inside the cylinders $Q_i$ that there is $\Gamma_1=\Gamma_1(\Phi,M,q_1,q_2,a,\hat{b})>0$ such that for all $i\in\N_0$ we have
\begin{align*}
&\sup_{t_0-\theta\rd_i^2\leq t\leq t_0}\frac{1}{\Phi'(\Phi^{-1}(k_i))}\int_{\mK_{\rd_i}}\big(w(\cdot,t)-k_i\big)_-^2\psi^2_i(\cdot,t)+\intQ[i]\big|\nabla\big((w-k_i)_-\psi_i\big)\big|^2\\
\leq\ &\Gamma_1\frac{\Phi^{-1}(k_i)}{k_i}\intQ[i](w-k_i)_-^2\psi_i|(\psi_i)_t|+\Gamma_1\intQ[i](w-k_i)_-^2|\nabla\psi_i|^2+\Gamma_1\Big(\int_{t_0-\theta\rd_i^2}^{t_0}\big|A_{k_i,\rd_i}^-(t)\big|^\frac{q_1(q_2-1)}{q_2(q_1-1)}\Big)^\frac{q_1-1}{q_1}.
\end{align*}
Multiplying this inequality by $\theta_{k_i}^{-1}=\Phi'\big(\Phi^{-1}(k_i)\big)$ and recalling that $k_i\leq k_0= \mu^-+k\leq 2M$ for all $i\in\N_0$ and Lemma~\ref{lem:phi} ensure the existence of $C_1=C_1(\Phi,M)>0$ such that
\begin{align*}
\Phi'\big(\Phi^{-1}(k_i)\big)\Phi^{-1}(k_i)\leq C_1 \Phi\big(\Phi^{-1}(k_i)\big)=C_1 k_i.
\end{align*}
we find $\Gamma_2=\Gamma_2(\Phi,M,q_1,q_2,a,\hat{b})>0$ satisfying
\begin{align*}
J_i\leq\Gamma_2\intQ[i](w-k_i)_-^2\psi_i|(\psi_i)_t|+\frac{\Gamma_2}{\theta_{k_0}}\intQ[i](w-k_i)_-^2|\nabla\psi_i|^2+\Gamma_2\Big(\int_{t_0-\theta\rd_i^2}^{t_0}\!\!\big|A_{k_i,\rd_i}^-(t)\big|^\frac{q_1(q_2-1)}{q_2(q_1-1)}\Big)^\frac{q_1-1}{q_1}
\end{align*}
for all $i\in\N_0$, where we additionally used $\theta_{2M}\leq\theta_{k_0}\leq\theta_{k_i}$ for all $i\in\N_0$ entailed by the monotonicity properties of $\Phi'$ and $\Phi^{-1}$.

For the case $\theta=\theta_{\mu^-}$ we note -- again by monotonicity of $\Phi'$ and $\Phi^{-1}$ -- that necessarily $\mu^-\geq k>0$, so that the special case of Lemma~\ref{lem:Lambda} becomes applicable, yielding
\begin{align*}
\Lambda^-_i(w):=\Lambda^-(w,k_i)\leq\frac{1}{2\Phi'(\Phi^{-1}(\mu^-))}(w-k_i)_-^2\quad\text{ and }\quad\Lambda_i^-(w)\geq \frac{1}{2\Phi'(\Phi^{-1}(k_i))}(w-k_i)_-^2
\end{align*}
a.e. in $Q_i$. Hence, this time we find from employing Lemma~\ref{lem:local-energy}, multiplying by $\theta_{k_i}^{-1}=\Phi'\big(\Phi^{-1}(k_i)\big)$ and again using $k_i\leq 2M$, that there is $\Gamma_3=\Gamma_3(\Phi,M,q_1,q_2,a,\hat{b})>0$ satisfying
\begin{align*}
J_i\leq\Gamma_3\frac{\theta_{\mu^-}}{\theta_{k_i}}\intQ[i](w-k_i)_-^2\psi_i|(\psi_i)_t|
&+\frac{\Gamma_3}{\theta_{k_0}}\intQ[i](w-k_i)_-^2|\nabla\psi_i|^2+\Gamma_3\Big(\int_{t_0-\theta\rd_i^2}^{t_0}\!\!\big|A_{k_i,\rd_i}^-(t)\big|^\frac{q_1(q_2-1)}{q_2(q_1-1)}\Big)^\frac{q_1-1}{q_1}
\end{align*}
for all $i\in\N_0$. Since $\mu^-\geq k$ we have $k_i\leq k_0=\mu^-+k\leq 2\mu^-$, so that $\theta_{k_i}\geq \theta_{2\mu^-}$. Accordingly, drawing on the fact that $\Phi^{-1}$ is subadditive due to its assumed concavity and $\Phi^{-1}(0)=0$, we obtain from Lemma~\ref{lem:phi} that there is $C_2=C_2(\Phi,M)>0$ satisfying
\begin{align*}
\frac{\theta_{\mu^-}}{\theta_{k_i}}&\leq \frac{\theta_{\mu^-}}{\theta_{2\mu^-}}=\frac{\Phi'(\Phi^{-1}(2\mu^-))}{\Phi'(\Phi^{-1}(\mu^-))}\leq\frac{\Phi'(2\Phi^{-1}(\mu^-))}{\Phi'(\Phi^{-1}(\mu^-))}\leq C_2 \quad\text{for all }i\in\N_0.
\end{align*}
Hence, letting $\Gamma_4=\Gamma_4(\Phi,M,q_1,q_2,a,\hat{b}):=\max\{\Gamma_2,\Gamma_3C_2,\Gamma_3\}>0$ we find that
\begin{align*}
J_i\leq\Gamma_4\intQ[i](w-k_i)_-^2\psi_i|(\psi_i)_t|
&+\frac{\Gamma_4}{\theta_{k_0}}\intQ[i](w-k_i)_-^2|\nabla\psi_i|^2+\Gamma_4\Big(\int_{t_0-\theta\rd_i^2}^{t_0}\!\!\big|A_{k_i,\rd_i}^-(t)\big|^\frac{q_1(q_2-1)}{q_2(q_1-1)}\Big)^\frac{q_1-1}{q_1}
\end{align*}
holds for all $i\in\N_0$ irrespective of a specific case for $\theta$. In light of the fact that by the definition of $k_i$ we have $(w-k_i)_-\leq (w-k_0)_-\leq(\hat{\mu}_--\mu^--k)_-\leq k$ a.e. in $Q_{i}$, we derive from the estimates provided by \eqref{eq:first-alt-1-psi-est} that by writing $A_i^-:=A_{k_i,\rd_i}^-$, $\widehat{q_1}:=\frac{2q_1(1+\kappa)}{q_1-1}$ and $\widehat{q_2}:=\frac{2q_2(1+\kappa)}{q_2-1}$ we obtain
\begin{align}\label{eq:first-alt-1-eq0}
J_i\leq \Gamma_4\frac{2^{2i+4}}{\theta\rd^2}k^2\big|A_{i}^-\big|+\Gamma_4\frac{2^{2i+6}}{\theta_{k_0}\rd^2}k^2\big|A_{i}^-\big|+\Gamma_4\Big(\int_{t_0-\theta\rd_i^2}^{t_0}\big|A_{k_i,\rd_i}^-(t)\big|^\nfrac{\widehat{q_1}}{\widehat{q_2}}\Big)^\nfrac{2(1+\kappa)}{\widehat{q_1}}
\end{align}
is valid for all $i\in\N_0$. Now, we introduce a change in variable by setting $\bar{t}:=\frac{t-t_0}{\theta}$ and correspondingly denote $\bar{Q}_i:=\mK_{\rd_i}\times(-\rd_i^2,0]$, $\bar{A}_i^-:=\bar{Q}_i\cap\big\{(x,\bar{t})\in\Omega\times(-\frac{t_0}{\theta},\frac{T-t_0}{\theta}):\big(w(x,\bar{t})-k_i\big)_->0\big\}$ and $\bar{A}^-_{i}(\bar{t}):=A_{k_i,\rd_i}^-(\theta\bar{t}+t_0)$. This way we have $|A_i^-|=\theta|\bar{A}_i^-|$ and we conclude from \eqref{eq:first-alt-1-eq0} that
\begin{align}\label{eq:first-alt-1-eq1}
&\sup_{-\rd_i^2\leq\,\bar{t}\,\leq 0}\int_{\mK_{\rd_i}}\big(w(\cdot,\bar{t})-k_i\big)_-^2\psi_i^2(\cdot,\bar{t})+\frac{\theta}{\theta_{k_i}}\iint\limits_{\bar{Q}_i}\big|\nabla\big((w-k_i)_-\psi_i\big)\big|^2\nonumber\\
&\leq\Gamma_4\frac{2^{2i+4}}{\rd^2}k^2 \big|\bar{A}_i^-\big|+\Gamma_4\frac{\theta}{\theta_{k_0}}\frac{2^{2i+6}}{\rd^2}k^2\big|\bar{A}_i^-\big|+\Gamma_4\theta^{\nfrac{2(1+\kappa)}{\widehat{q_1}}}\Big(\int_{-\rd_i^2}^0\big|\bar{A}_{i}^-(\bar{t})\big|^\nfrac{\widehat{q_1}}{\widehat{q_2}}\intd{\bar{t}}\Big)^\nfrac{2(1+\kappa)}{\widehat{q_1}}
\end{align}
for all $i\in\N_0$. 
If $\mu^-\leq k$, we have $\frac{k}{2}\leq k_i\leq k_0\leq 2k$ and the subadditivity of $\Phi^{-1}$ then implies that
\begin{align*}
\Phi'\big(\Phi^{-1}(k_0)\big)\leq\Phi'\big(\Phi^{-1}(2k)\big)\leq\Phi'\big(2\Phi^{-1}(k)\big)\quad\text{and}\quad \Phi'\big(\Phi^{-1}(k)\big)\leq\Phi'\big(2\Phi^{-1}(\tfrac{k}{2})\big).
\end{align*}
Accordingly, in this case we find from $\theta=\theta_k$, $k\leq \mu^+\leq M$ and \eqref{eq:phi-prop-3} of Lemma~\ref{lem:phi} that there is some $C_3=C_3(M)>1$ for which
\begin{align*}
\frac{\theta}{\theta_{k_0}}=\frac{\Phi'\big(\Phi^{-1}(k_0)\big)}{\Phi'\big(\Phi^{-1}(k)\big)}\leq\frac{\Phi'\big(2\Phi^{-1}(k)\big)}{\Phi'\big(\Phi^{-1}(k)\big)}\leq C_3\quad\text{and}\quad \frac{1}{C_3}\leq\frac{\Phi'\big(\Phi^{-1}(\tfrac{k}{2})\big)}{\Phi'\big(2\Phi^{-1}(\tfrac{k}{2})\big)}\leq \frac{\Phi'\big(\Phi^{-1}(k_i)\big)}{\Phi'\big(\Phi^{-1}(k)\big)}=\frac{\theta}{\theta_{k_i}}
\end{align*}
for all $i\in\N_0$. Similarly, when $\mu^->k$ we have $\mu^-< k_i \leq k_0< 2\mu^-$ and thus
\begin{align*}
\Phi'\big(\Phi^{-1}(k_0)\big)\leq\Phi'\big(2\Phi^{-1}(\mu^-)\big)\quad\text{and}\quad\Phi'\big(\Phi^{-1}(\mu^-)\big)\leq\Phi'\big(\Phi^{-1}(k_i)\big),
\end{align*}
which together with Lemma~\ref{lem:phi} yields $C_4=C_4(M)> 1$ such that
\begin{align*}
\frac{\theta}{\theta_{k_0}}\leq \frac{\Phi'\big(2\Phi^{-1}(\mu^-)\big)}{\Phi'\big(\Phi^{-1}(\mu^-)\big)}\leq C_4\quad\text{and}\quad \frac{1}{C_4}\leq\frac{\Phi'\big(\Phi^{-1}(\mu^-)\big)}{\Phi'\big(\Phi^{-1}(\mu^-)\big)}\leq\frac{\theta}{\theta_{k_i}}\quad\text{for all }i\in\N_0.
\end{align*}
Hence, we obtain from \eqref{eq:first-alt-1-eq1} that in both cases there is $\Gamma_5=\Gamma_5(\Phi,M,q_1,q_2,a,\hat{b})>0$ such that
\begin{align}\label{eq:first-alt-1-eq1.5}
\bar{J}_i&:=\sup_{-\rd_i^2\leq\,\bar{t}\,\leq 0}\int_{\mK_{\rd_i}}\big(w(\cdot,\bar{t})-k_i\big)_-^2\psi_i^2(\cdot,\bar{t})+\iint\limits_{\bar{Q}_i}\big|\nabla\big((w-k_i)_-\psi_i\big)\big|^2\nonumber\\
&\leq \Gamma_5\frac{2^{2i}k^2}{\rd^2}\big|\bar{A}_i^-\big|+\Gamma_5\theta^{\frac{q_1-1}{q_1}}\Big(\int_{-\rd_i^2}^0\big|\bar{A}_{i}^-(\bar{t})\big|^\nfrac{\widehat{q_1}}{\widehat{q_2}}\intd{\bar{t}}\Big)^\nfrac{2(1+\kappa)}{\widehat{q_1}}\quad\text{for all }i\in\N_0.
\end{align}
In order to relate this to the smaller cylinder in step $i+1$, we estimate $\bar{J}_i$ from below using Corollary~\ref{cor:sob}. Since $\psi_i=1$ on $\bar{Q}_{i}$ and $\psi_i\geq0$ on $\bar{Q}_{i}\setminus \bar{Q}_{i+1}$, we have
\begin{align*}
\iint\limits_{\bar{Q}_{i+1}}(w-k_i)_-^2\leq \iint\limits_{\bar{Q}_i}(w-k_i)_-^2\psi_i^2
\end{align*}
for all $i\in\N_0$.
From Corollary \ref{cor:sob} we thereby find $C_5>0$ and $C_6=C_6(M)>0$ such that
\begin{align}\label{eq:first-alt-1-eq2}
\iint\limits_{\bar{Q}_{i+1}}(w-k_i)_-^2\nonumber\leq\ & C_5\big|\bar{A}_{i}^-\big|^\frac{2}{\dimN+2}\Big(\esssup_{\bar{t}\in(-\rd_i^2,0)}\|(w-k_i)_-\psi_i\|_{\LSpn{2}{\mK_{\rd_i}}}+\Big(\iint\limits_{\bar{Q}_{i}}\big|\nabla\big((w-k_i)_-\psi_i\big)\big|^2\Big)^\frac{1}{2}\Big)^2\nonumber\\
\leq\ & C_6\big|\bar{A}_{i}^-\big|^\frac{2}{\dimN+2}\Big(\esssup_{-\rd_i^2\leq \bar{t}\leq 0}\int_{\mK_{\rd_i}}(w-k_i)_-^2\psi_i^2+\iint\limits_{\bar{Q}_i}\big|\nabla\big((w-k_i)_-\psi_i\big)\big|^2\Big)\nonumber\\
\leq\ &C_6\big|\bar{A}_{i}^-\big|^\frac{2}{\dimN+2} \bar{J}_i\quad\text{for all }i\in\N_0.
\end{align}
Moreover, note that $(w-k_i)_-=-w+k_i\geq k_i-k_{i+1}$ whenever $w\leq k_{i+1}<k_i$, and thus $(w-k_i)_-^2\geq (k_i-k_{i+1})^2$ a.e. on $\bar{A}_{i+1}^-$. Accordingly, since $k_i-k_{i+1}=\frac{k}{2^{i+2}}$,
\begin{align*}
\iint\limits_{\bar{Q}_{i+1}}(w-k_i)_-^2\geq (k_i-k_{i+1})^2\big|\bar{A}_{i+1}^-\big|=\frac{k^2}{2^{2i+4}} \big|\bar{A}_{i+1}^-\big|,
\end{align*}
which upon combination with \eqref{eq:first-alt-1-eq1.5} and \eqref{eq:first-alt-1-eq2} implies that there is $\Gamma_6=\Gamma_6(\Phi,M,q_1,q_2,a,\hat{b})>0$ such that
\begin{align*}
\big|\bar{A}_{i+1}^-\big|&\leq\frac{C_6 2^{2i+4}}{k^2}\big|\bar{A}_{i}^-\big|^\frac{2}{\dimN+2} \bar{J}_i\\
&\leq \Gamma_6\frac{2^{4i}}{\rd^2}\big|\bar{A}_{i}^-\big|^{1+\frac{2}{\dimN+2}}+\Gamma_6\frac{2^{2i}}{k^2}\theta^{\frac{q_1-1}{q_1}}\big|\bar{A}_{i}^-\big|^\frac{2}{\dimN+2}\Big(\int_{-\rd_i^2}^{0}\big|\bar{A}_{i}^-(\bar{t})\big|^\nfrac{\widehat{q_1}}{\widehat{q_2}}\intd\bar{t}\Big)^\nfrac{2(1+\kappa)}{\widehat{q_1}}
\end{align*}
holds for all $i\in\N_0$. Dividing this by $\rd^{\dimN+2}$ and introducing the quantities $$\bar{Y}_i:=\frac{\big|\bar{A}_{i}^-\big|}{\rd^{\dimN+2}}\ \text{ as well as }\ \bar{Z}_i:=\frac{1}{\rd^N}\Big(\int_{-\rd_i^2}^{0}\big|\bar{A}_{i}^-(\bar{t})\big|^\nfrac{\widehat{q_1}}{\widehat{q_2}}\intd\bar{t}\Big)^\nfrac{2}{\widehat{q_1}},\quad i\in\N_0,$$ we therefore have
\begin{align}\label{eq:first-alt-1-Y}
\bar{Y}_{i+1}&\leq \Gamma_6 2^{4i}\Big(\bar{Y}_i^{1+\frac{2}{\dimN+2}}+\frac{\rd^{\dimN\kappa}}{k^{2}}\theta^{\frac{q_1-1}{q_1}}\bar{Y}_i^\frac{2}{\dimN+2}\bar{Z}_i^{1+\kappa} \Big)\quad\text{for all }i\in\N_0.
\end{align}
Regarding $\bar{Z}_{i+1}$, we note that $(w-k_i)_-^2\geq (k_i-k_{i+1})^2$ whenever $w<k_{i+1}$, and therefore we may estimate for $i\in\N_0$
\begin{align*}
&\rd^\dimN(k_i-k_{i+1})^2 \bar{Z}_{i+1}\\
=\ &(k_i-k_{i+1})^2\Big(\int_{-\rd_{i+1}^2}^{0}\big|\bar{A}_{i+1}^-(\bar{t})\big|^\nfrac{\widehat{q_1}}{\widehat{q_2}}\intd\bar{t}\Big)^\nfrac{2}{\widehat{q_1}}\\
\leq\ &\Big(\int_{-\rd_{i+1}^2}^{0}\Big(\int_{\mK_{\rd_{i+1}}}\big(w(x,\bar{t})-k_i\big)_-^{\widehat{q_2}}\Big)^\nfrac{\widehat{q_1}}{\widehat{q_2}}\Big)^\nfrac{2}{\widehat{q_1}}\\
\leq\ &\Big(\int_{-\rd_{i}^2}^{0}\Big(\int_{\mK_{\rd_{i}}}\big(w(x,\bar{t})-k_i\big)_-^{\widehat{q_2}}\psi_i^{\widehat{q_2}}(\cdot,\bar{t})\Big)^\nfrac{\widehat{q_1}}{\widehat{q_2}}\Big)^\nfrac{2}{\widehat{q_1}}\\
=\ &\big\|\big((w-k_i)_-\psi_i\big)(x,\bar{t})\big\|_{\LSpn{\widehat{q_1}}{(-\rd_i^2,0);\LSpn{\widehat{q_2}}{\mK_{\rd_i}}}}^2.
\end{align*}
Here, our assumptions on $q_1,q_2$ and $\kappa$ ensure that the embedding Lemma~\ref{lem:lemmaX} with $p=2$ applies, yielding $C_7>0$ satisfying
\begin{align*}
\rd^\dimN(k_i-k_{i+1})^2 \bar{Z}_{i+1}\leq C_7 \bar{J}_i\quad\text{for }i\in\N_0,
\end{align*}
which in view of $k^2=2^{2i+4}(k_i-k_{i+1})^2$ and \eqref{eq:first-alt-1-eq1.5} entails that with some $\Gamma_7=\Gamma_7(\Phi,M,q_1,q_2,a,\hat{b})>0$,
\begin{align}\label{eq:first-alt-1-Z}
\bar{Z}_{i+1}&\leq \frac{2^{2i+4}C_7}{\rd^\dimN k^2}\Big(\Gamma_5 \frac{2^{2i}k^2}{\rd^2}\rd^{\dimN+2}\bar{Y}_i+\Gamma_5\rd^{\dimN+\dimN\kappa}\theta^{\frac{q_1-1}{q_1}} \bar{Z}_i^{1+\kappa}\Big)\leq 2^{4i}\Gamma_7 \Big(\bar{Y}_{i}+\frac{\rd^{\dimN\kappa}}{k^{2}}\theta^{\frac{q_1-1}{q_1}}\bar{Z}_{i}^{1+\kappa}\Big)
\end{align}
for $i\in\N_0$. We conclude from \eqref{eq:first-alt-1-Y} and \eqref{eq:first-alt-1-Z} that with $\Gamma_8=\Gamma_8(\Phi,M,q_1,q_2,a,\hat{b}):=\max\big\{\Gamma_6,\Gamma_7 \big\}+1>1$ we have
\begin{align*}
\bar{Y}_{i+1}\leq \Gamma_8 2^{4i}\Big(\bar{Y}_i^{1+\frac{2}{\dimN+2}}+\frac{\rd^{\dimN\kappa}}{k^{2}}\theta^{\frac{q_1-1}{q_1}}\bar{Y}_i^\frac{2}{\dimN+2} \bar{Z}_i^{1+\kappa}\Big)\quad\text{and}\quad \bar{Z}_{i+1}\leq \Gamma_8 2^{4i}\Big(\bar{Y}_i+\frac{\rd^{\dimN\kappa}}{k^{2}}\theta^{\frac{q_1-1}{q_1}}\bar{Z}_i^{1+\kappa}\Big)
\end{align*}
for $i\in\N_0$. Assuming $k^{2}\theta^{-\frac{q_1-1}{q_1}}>\rd^{\dimN\kappa}$, we obtain recursive inequalities precisely of the form requested by Lemma~\ref{lem:fast-geo-conv} and are left with checking the condition on $\bar{Y}_0+\bar{Z}_0^{1+\kappa}$. Recall $\rd_0=\rd$ and $k_0=\mu^-+k$. Then, first assuming $\widehat{q_1}\geq\widehat{q_2}$, we use
$$\frac{2(1+\kappa)}{\widehat{q_1}}\left(\frac{\widehat{q_1}}{\widehat{q_2}}-1\right)=2(1+\kappa)\left(\frac{1}{\widehat{q_2}}-\frac{1}{\widehat{q_1}}\right)=\frac{1}{q_1}-\frac{1}{q_2}$$ and that $|\bar{A}_0^-(\bar{t})|\leq |K_\rd|$ for all $\bar{t}\in(-\tfrac{t_0}{\theta},\frac{T-t_0}{\theta})$ to estimate
\begin{align*}
\bar{Y}_0+\bar{Z}_0^{1+\kappa}&=\frac{\big|\bar{A}_0^-\big|}{\rd^{\dimN+2}}+\frac{1}{\rd^{\dimN(1+\kappa)}}\Big(\int_{-\rd^2}^{0}\big|\bar{A}_{0}^-(\bar{t})\big|^\nfrac{\widehat{q_1}}{\widehat{q_2}}\intd{\bar{t}}\Big)^\nfrac{2(1+\kappa)}{\widehat{q_1}}\\
&\leq \frac{\big|\bar{A}_{\mu_-+k,\rd}^-\big|}{\rd^{\dimN+2}}+\frac{|\mK_\rd|^{\frac{1}{q_1}-\frac{1}{q_2}}}{\rd^{\dimN(1+\kappa)}}\big|\bar{A}^-_{\mu^-+k,\rd}\big|^\nfrac{2(1+\kappa)}{\widehat{q_1}}\\
&\leq \frac{\big|\bar{A}_{\mu_-+k,\rd}^-\big|}{\rd^{\dimN+2}}+\omega_\dimN^{\frac{1}{q_1}-\frac{1}{q_2}}\rd^{\frac{\dimN}{q_1}-\frac{\dimN}{q_2}-\dimN-\dimN\kappa}\big|\bar{A}^-_{\mu^-+k,r}\big|^\frac{q_1-1}{q_1},
\end{align*}
where $\omega_\dimN:=2^\dimN$ denotes the volume of the $\dimN$-cube of edge length $2$. 
Here, we find by the assumption \eqref{eq:first-alt-1-ass1} that
$$\theta|\bar{A}^-_{\mu^-+k,\rd}|=|A_{\mu^-+k,\rd}^-|\leq\sig_0|Q_\rd^{\theta}|= \omega_\dimN\sig_0\theta\rd^{\dimN+2},$$ and, since in light of \eqref{eq:cond-kappa} $2-\frac{2}{q_1}-\frac{\dimN}{q_2}=\dimN\kappa$, we have
\begin{align}\label{eq:first-alt-eq3}
\bar{Y}_0+\bar{Z}_0^{1+\kappa}\leq \omega_\dimN \sig_0 +\rd^{2-\frac{2}{q_1}-\frac{\dimN}{q_2}-\dimN\kappa}\omega_\dimN^{\frac{q_2-1}{q_2}} \sig_0^{\frac{q_1-1}{q_1}}=\omega_\dimN\sig_0+\omega_\dimN^{\frac{q_2-1}{q_2}}\sig_0^{\frac{q_1-1}{q_1}}.
\end{align}
On the other hand, if $\widehat{q_1}<\widehat{q_2}$ we draw on Jensen's inequality applied to the concave function $s\mapsto s^{\nfrac{\widehat{q_1}}{\widehat{q_2}}},$ $s\geq 0$ to estimate
\begin{align*}
\int_{-r^2}^0\big|\bar{A}_0^-(\bar{t})\big|^{\nfrac{\widehat{q_1}}{\widehat{q_2}}}\intd\bar{t}\leq \rd^{2(1-\nfrac{\widehat{q_1}}{\widehat{q_2}})}\Big(\int_{-r^2}^0\big|\bar{A}_0^-(\bar{t})\big|\intd\bar{t}\Big)^\nfrac{\widehat{q_1}}{\widehat{q_2}} =\rd^{2(1-\nfrac{\widehat{q_1}}{\widehat{q_2}})}\big|\bar{A}_0^-\big|^{\nfrac{\widehat{q_1}}{\widehat{q_2}}}
\end{align*}
 and then use $$\frac{4(1+\kappa)}{\widehat{q_1}}\left(1-\frac{\widehat{q_1}}{\widehat{q_2}}\right)=4(1+\kappa)\left(\frac{1}{\widehat{q_1}}-\frac{1}{\widehat{q_2}}\right)=\frac{2}{q_2}-\frac{2}{q_1}$$ to obtain 
\begin{align}\label{eq:first-alt-eq3.5}
\bar{Y}_0+\bar{Z}_0^{1+\kappa}&\leq\frac{\big|\bar{A}_{\mu_-+k,\rd}^-\big|}{\rd^{\dimN+2}}+\rd^{\frac{2}{q_2}-\frac{2}{q_1}-\dimN-\dimN\kappa}\big|\bar{A}^-_{\mu^-+k,\rd}\big|^\frac{q_2-1}{q_2}\nonumber\\
&\leq\omega_\dimN\sig_0+\rd^{\frac{(\dimN+2)(q_2-1)}{q_2}+\frac{2}{q_2}-\frac{2}{q_1}-\dimN-\dimN\kappa}\omega_\dimN^{\frac{q_2-1}{q_2}}\sig_0^{\frac{q_2-1}{q_2}}=\omega_\dimN\sig_0+\omega_\dimN^\frac{q_2-1}{q_2}\sig_0^{\frac{q_2-1}{q_2}}
\end{align} 
similarly as before. Hence, we conclude from \eqref{eq:first-alt-eq3} and \eqref{eq:first-alt-eq3.5} that in both of the cases $\widehat{q_1}\geq\widehat{q_2}$ and $\widehat{q_1}<\widehat{q_2}$, taking $\sig_0=\sig_0(\Phi,M,q_1,q_2,a,\hat{b})$ so small that we have $\sig_0<\frac{1}{2}$ and that with $q:=\min\big\{\frac{q_1-1}{q_1},\frac{q_2-1}{q_2}\big\}\in(0,1)$ and $\sigma=\min\{\kappa,\frac{2}{\dimN+2}\}=\frac{2}{\dimN+2}\in(0,1)$ we also have
\begin{align*}
\omega_\dimN \sig_0+\omega_\dimN^{\frac{q_2-1}{q_2}} \sig_0^{q}\leq \big(2\Gamma_8\big)^{-\frac{1+\kappa}{\sigma}}2^{-\frac{4(1+\kappa)}{\sigma^2}},
\end{align*}
entails
\begin{align*}
\bar{Y}_0+\bar{Z}_0^{1+\kappa}\leq \big(2\Gamma_8\big)^{-\frac{1+\kappa}{\sigma}}2^{-\frac{4(1+\kappa)}{\sigma^2}}.
\end{align*}
Therefore, Lemma~\ref{lem:fast-geo-conv} becomes applicable to say $\bar{Y}_i\to0$ as $i\to\infty$ and reverting to the original variable $t$ we arrive at $$0=\frac{\big|\bar{A}_{k_\infty,\rd_\infty}^-\big|}{\rd^{\dimN+2}}=\frac{|A_{k_\infty,\rd_\infty}^-|}{\theta\rd^{\dimN+2}}=\frac{\omega_\dimN|A_{k_\infty,\rd_\infty}^-|}{|Q_{\rd}^{\theta}|}.$$ Recalling the values for $k_\infty$ and $\rd_\infty$ specified in \eqref{eq:first-alt-1-rinfty}, we obtain \eqref{eq:first-alt-concl}.
\end{bew}

In contrast to the first alternative presented in \cite[Proposition 4.1]{HwangZhang19} for the special case of $\Phi(s)=s^m$, $m\geq 1$, we got rid of the additional assumption that either $\mu^-\leq k$ or $m=1$. This allows for a sleeker iteration procedure in Section~\ref{sec4:iter} and, in particular, facilitates easier tracking of dependencies.

\subsection{The second alternative}\label{sec4:alt2}
We will now consider the second alternative in which the set where $w$ is significantly larger than the essential infimum is too large. For this purpose, we assume that the assumption \eqref{eq:first-alt-1-ass1} in Lemma~\ref{lem:first-alt-1} is false, i.e. $$\big|\big\{(x,t)\in Q_\rd^{\theta}:(w-\mu^--k)_->0\big\}\big|>\sig_0 |Q_\rd^{\theta}|,$$ then
\begin{align*}
\big|\big\{(x,t)\in Q_\rd^{\theta} : w\geq \mu^-+k\big\}\big|&=
\big|\big\{(x,t)\in Q_\rd^{\theta} :(w-\mu^--k)_-=0\big\}\big|\\&=|Q_\rd^{\theta}|-\big|\big\{(x,t)\in Q_\rd^{\theta} :(w-\mu^--k)_->0\big\}\big|\leq\big(1-\sig_0\big)|Q_\rd^{\theta}|.
\end{align*}
If we moreover suppose that $\essosc_{Q_{2\rd}^{\vartheta}}w\geq 2k$ -- and treat the case where the essential oscillation is small separately later -- we then have $\mu^+-k\geq \mu^-+k$ in $Q_{2r}^{\theta}\subseteq Q_{2r}^{\vartheta}$, so that
\begin{align*}
\big|\big\{(x,t)\in Q_\rd^{\theta} : (w-\mu^++k)_+>0\big\}\big|&=\big|\big\{(x,t)\in Q_\rd^{\theta} : w> \mu^+-k\big\}\big|\\
&\leq \big|\big\{(x,t)\in Q_\rd^{\theta} : w\geq \mu^-+k\big\}\big|\leq \big(1-\sig_0\big)|Q_\rd^{\theta}|.
\end{align*}
Accordingly, if \eqref{eq:first-alt-1-ass1} is false, then $w$ is close to $\mu^+$ throughout a large portion of $Q_\rd^\theta$. The next lemma shows, that then also at some fixed time level the corresponding superlevel set of $w$ is large when compared with the box.

\begin{lemma}\label{lem:alt2-time-fixed-level-set-est}
Let $M>0$ and denote by $\sig_0=\sig_0(\Phi,M,q_1,q_2,a,\hat{b})\in(0,\tfrac12)$ the number from Lemma~\ref{lem:first-alt-1}. Let $w$ be a nonnegative bounded weak solution of \eqref{eq:trans-gen-eq} satisfying \eqref{eq:linfty-bound-w} and let $x_0\in\Omega$, $0<t_0<T$, $0<\vartheta$ and $0<\rd<\min\big\{\frac{1}{\sqrt{\dimN}}\dist(x_0,\romega),\sqrt{\frac{t_0}{4\vartheta}},1\big\}$. Set $\mu^-:=\essinf_{Q_{2\rd}^{\vartheta}}w<\esssup_{Q_{2\rd}^{\vartheta}}w=:\mu^+$, let $k\in(0,M]$ and introduce $\theta:=\min\big\{\theta_k,\theta_{\mu^-}\big\}$.  If $\vartheta\geq \theta$ and if $$\big|\big\{(x,t)\in Q^{\theta}_\rd:(w-\mu^++k)_+>0\big\}\big|\leq \big(1-\sig_0\big)|Q^{\theta}_\rd|$$ then for any $l\in(0,k]$ and $\gamma\in\big(0,\sig_0\big)$ there is $\tau\in\big[t_0-\theta\rd^2,t_0-\gamma\theta\rd^2\big]$ such that 
\begin{align*}
\big|A^+_{\mu^+-l,\rd}(\tau)\big|\leq\frac{1-\sig_0 }{1-\gamma}|\mK_\rd|.
\end{align*}
\end{lemma}

\begin{bew}
First note that since $l\leq k$, we have $A^+_{\mu^+-l,\rd}(t)\subseteq A^+_{\mu^+-k,\rd}(t)$ for each $t\in[t_0-\theta\rd^2,t_0]$, so that it clearly suffices to verify the claim for $l=k$.\\
Now, assuming to the contrary that with $\gamma\in(0,\sig_0)\subset(0,\frac12)$ we had
$$\big|A_{\mu^+-k,\rd}^+(\tau)\big|>\frac{1-\sig_0}{1-\gamma}|\mK_\rd|\quad\text{for all }\tau\in\big[t_0-\theta\rd^2,t_0-\gamma\theta\rd^2\big],$$
then
\begin{align*}
\big|\big\{(x,t)\in Q^{\theta}_\rd:(w-\mu^++k)_+>0\big\}\big|&=\int_{t_0-\theta\rd^2}^{t_0}\big|A_{\mu^+-k,\rd}^+(\tau)\big|\intd \tau\\
&\geq\int_{t_0-\theta\rd^2}^{t_0-\gamma\theta\rd^2}\big|A^+_{\mu^+-k,\rd}(\tau)\big|\intd\tau\\
&>\frac{1-\sig_0 }{1-\gamma}(1-\gamma)\theta\rd^2|\mK_\rd|\\
&=\big(1-\sig_0 \big)|Q^{\theta}_\rd|,
\end{align*}
which clearly contradicts the assumption $\big|\big\{(x,t)\in Q^{\theta}_\rd:(w-\mu^++k)_+>0\big\}\big|\leq \big(1-\sig_0 \big)|Q^{\theta}_\rd|$.
\end{bew}

The next lemma now extends the quantitative information on the superset from the previous lemma for later times by utilizing the logarithmic estimate from Lemma~\ref{lem:log-est}.
\begin{lemma}\label{lem:expansion-in-time}
Let $M>0$ and denote by $\sig_0=\sig_0(\Phi,M,q_1,q_2,a,\hat{b})\in(0,\tfrac12)$ the number from Lemma~\ref{lem:first-alt-1}. Then, there is $\delta_0=\delta_0(\Phi,M,q_1,q_2,a,\hat{b})\in(0,\frac{1}{2})$ such that for all nonnegative bounded weak solutions $w$ of \eqref{eq:trans-gen-eq} satisfying \eqref{eq:linfty-bound-w} the following holds: Given $x_0\in\Omega$, $0<t_0<T$, $0<\vartheta$ and $0<\rd<\min\big\{\frac{1}{\sqrt{\dimN}}\dist(x_0,\romega),\sqrt{\frac{t_0}{4\vartheta}},1\big\}$ set $\mu^-:=\essinf_{Q_{2\rd}^{\vartheta}}w<\esssup_{Q_{2\rd}^{\vartheta}}w=:\mu^+$. Let $k\in(0,\mu^+]$ and introduce $\theta:=\min\big\{\theta_k,\theta_{\mu^-}\big\}$. If
\begin{align}\label{eq:exp-in-time}
\theta\leq\vartheta,\quad \text{and}\quad \mu^+-\mu^-\leq 4k,\quad\text{and}\quad\big|\big\{(x,t)\in Q^{\theta}_\rd:(w-\mu^++k)_+>0\big\}\big|\leq \big(1-\sig_0\big)|Q^{\theta}_\rd|,
\end{align}
then either
\begin{align*}
\delta_0^2 k^2\theta^{-\frac{q_1-1}{q_1}}\leq \rd^{\dimN\kappa}
\end{align*}
or
\begin{align}\label{eq:exp-in-time1}
\big|A^+_{\mu^+-\delta_0 k,\rd}(t)\big|<\Big(1-\frac{\sig_0^2}{4}\Big)|\mK_\rd|\quad\text{for all }t\in\big(t_0-\tfrac{\sig_0}{2}\theta\rd^2,t_0\big].
\end{align}
\end{lemma}

\begin{bew}
With $\sig_0$ taken from Lemma~\ref{lem:first-alt-1}, we set $\gamma:=\frac{\sig_0}{2}$ and denote by $C_1=C_1(\Phi,M,q_1,q_2,a,\hat{b})>0$ the constant provided by Lemma~\ref{lem:log-est}. In view of the monotonicity of $\Phi^{-1}$, we find from $\mu^-<\mu^+\leq M$ and $k\leq \mu^+\leq M$ that $\Phi^{-1}(\mu^-)\leq \Phi^{-1}(M)$ and $\Phi^{-1}(k)\leq \Phi^{-1}(M)$. Thus, according to Lemma~\ref{lem:phi}, we can pick $C_2=C_2(\Phi,M)>0$ such that
\begin{align}\label{eq:exp-in-time-eq0}
\frac{\Phi'\big(5\Phi^{-1}(k)\big)}{\Phi'\big(\Phi^{-1}(k)\big)}\leq C_2\quad\text{and}\quad\frac{\Phi'\big(5\Phi^{-1}(\mu^-)\big)}{\Phi'\big(\Phi^{-1}(\mu^-)\big)}\leq C_2.
\end{align}
Next, we define $\lambda:=\frac{\gamma^2}{\dimN}\in(0,\tfrac{1}{16})$ and pick a cutoff function $\zeta=\zeta_\rd\in\CSp{\infty}{\mK_\rd}$ satisfying $\zeta_{\vert_{\partial \mK_\rd}}=0$, $0\leq \zeta\leq 1$ in $\mK_\rd$ as well as $\zeta=1$ in $\mK_{(1-\lambda)\rd}$ and $|\nabla\zeta|\leq\frac{4}{\lambda\rd}$. With these choices, we then take $j=j(\Phi,M,q_1,q_2,a,\hat{b})\in\N$ so large that
\begin{align*}
\frac{j^2}{(j-1)^2}\cdot\frac{1-2\gamma}{1-\gamma}&\leq 1-4\gamma^2,\\
\frac{j}{(j-1)^2\ln(2)}\cdot\frac{96C_2}{\lambda^2}&\leq \gamma^2,\\
C_1\theta_M^{-1}\omega_\dimN^{-\frac{1}{q_2}}\Big(\frac{1+\ln(2)(M+1)j}{\ln^2(2)(j-1)^2}\Big)&\leq \gamma^2,
\end{align*}
which is possible due to $\gamma<\frac{1}{4}$, and let $\delta_0=\delta_0(\Phi,M,q_1,q_2,a,\hat{b})=2^{-j}\in(0,\frac{1}{2})$. Moreover, we notice that the assumptions of the current lemma allow for an application of Lemma~\ref{lem:alt2-time-fixed-level-set-est}, yielding a time $\tau\in\big[t_0-\theta\rd^2,t_0-\gamma\theta\rd^2\big]$ for which
\begin{align}\label{eq:exp-in-time-eq1}
\big|A^+_{\mu^+-k,\rd}(\tau)\big|\leq\frac{1-2\gamma}{1-\gamma}|\mK_\rd|.
\end{align}
Having these preparations at hand, we now consider the cylinder $Q:=Q_{\rd}(x_0,\tau,t)$ for arbitrary $t\in(\tau,t_0)$. In view of $(\tau,t]\subseteq(t_0-\theta\rd^2,t_0]$ we have $Q\subseteq Q_{\rd}^{\theta}\subseteq Q_{\rd}^{\vartheta}$ and hence also $\esssup_{Q} w\leq \mu^+$. Thus, we may draw on Lemma~\ref{lem:log-est} for $\xi=\xi_{\mu^+,k,\delta_0}$, as defined in \eqref{eq:xi-def}, inside the cylinder $Q=Q_{\rd}(x_0,\tau,t)$ to find that
\begin{align*}
\int_{\mK_\rd}\!\!\xi^2\big(w(\cdot,t)\big)\zeta^2\leq\int_{\mK_\rd}\!\!\xi^2\big(w(\cdot,\tau)\big)\zeta^2&+96\Phi'\big(\Phi^{-1}(\mu^+)\big)\frac{|t-\tau|}{\lambda^2 \rd^2}\ln\Big(\frac{1}{\delta_0}\Big)|\mK_\rd|\\
&+C_1\Big(\frac{1+\ln(\frac{1}{\delta_0})}{\delta_0^2 k^2}+\frac{\ln(\frac{1}{\delta_0})}{\delta_0 k}\Big)\Big(\int_{\tau}^t\big|A_{\mu^+-k,\rd}^+(s)\big|^{\frac{q_1(q_2-1)}{q_2(q_1-1)}}\intd s\Big)^{\frac{q_1-1}{q_1}},
\end{align*}
where again we made use of the monotonicity of $\Phi'$ and $\Phi^{-1}$. Then, since $t-\tau\leq t_0-t_0+\theta\rd^2=\theta\rd^2$ and $\big|A^+_{\mu^+-k,\rd}\big|\leq |\mK_\rd|$ for all $s\in(\tau,t)$ we obtain
\begin{align}\label{eq:exp-in-time-eq2}
\int_{\mK_\rd}\!\!\!\xi^2\big(w(\cdot,t)\big)\zeta^2\leq\!\int_{\mK_\rd}\!\!\!\xi^2\big(w(\cdot,\tau)\big)\zeta^2&+\frac{96}{\lambda^2}\frac{\theta}{\theta_{\mu^+}}\ln\!\Big(\frac{1}{\delta_0}\Big)|\mK_\rd|+C_1\Big(\frac{1+\ln(\frac{1}{\delta_0})}{\delta_0^2 k^2}+\frac{\ln(\frac{1}{\delta_0})}{\delta_0 k}\Big)\big(\theta\rd^2\big)^\frac{q_1-1}{q_1}|\mK_\rd|^{\frac{q_2-1}{q_2}}
\end{align}
for $t\in(\tau,t_0)$. Now, if $\mu^+-\mu^-\leq 4k$ and $\mu^-\leq k$, we have $\mu^+\leq 5k$ and in view of the subadditivity of $\Phi^{-1}$ also $\Phi^{-1}(\mu^+)\leq 5\Phi^{-1}(k)$. Similarly, if $\mu^+-\mu^-\leq 4k$ and $k<\mu^-$, we have $\mu^+\leq 5\mu^-$ yielding $\Phi^{-1}(\mu^+)\leq 5\Phi^{-1}(\mu^-)$, whence we find from \eqref{eq:exp-in-time-eq0} that in both cases $\frac{\theta}{\theta_{\mu^+}}\leq C_2$. Therefore, since $\xi(s)\equiv0$ for $s\in[0,\mu^+-(1-\delta_0)k]$ with $\xi(s)\leq\ln(\frac{1}{\delta_0})$ on $[0,\mu^+]$, we infer from a combination of \eqref{eq:exp-in-time-eq2}, \eqref{eq:exp-in-time-eq1} and \eqref{eq:cond-kappa} that
\begin{align}\label{eq:exp-in-time-eq3}
&\int_{\mK_\rd}\!\!\!\xi^2\big(w(\cdot,t)\big)\zeta^2\nonumber\\
\leq\;&\ln^2\!\Big(\frac{1}{\delta_0}\Big)\big|A^+_{\mu^+-k,\rd}(\tau)\big|+96C_2\ln\!\Big(\frac{1}{\delta_0}\Big)\frac{|\mK_\rd|}{\lambda^2}+C_1\Big(\frac{1+\ln(\frac{1}{\delta_0})}{\delta_0^2 k^2}+\frac{\ln(\frac{1}{\delta_0})}{\delta_0 k}\Big)\big(\theta\rd^2\big)^\frac{q_1-1}{q_1}|\mK_\rd|^{\frac{q_2-1}{q_2}}\nonumber\\
\leq\;&\ln^2\!\Big(\frac{1}{\delta_0}\Big)\frac{1-2\gamma}{1-\gamma}|\mK_\rd|+96C_2\ln\Big(\frac{1}{\delta_0}\Big)\frac{|\mK_\rd|}{\lambda^2}+C_1\Big(\frac{1+\ln(\frac{1}{\delta_0})}{\delta_0^2 k^2}+\frac{\ln(\frac{1}{\delta_0})}{\delta_0 k}\Big)\theta^\frac{q_1-1}{q_1}\omega_\dimN^{-\frac{1}{q_2}}\rd^{2-\frac{2}{q_1}-\frac{\dimN}{q_2}}|\mK_\rd|\nonumber\\
=\;&\left(\ln^2\!\Big(\frac{1}{\delta_0}\Big)\frac{1-2\gamma}{1-\gamma}+96C_2\ln\Big(\frac{1}{\delta_0}\Big)\frac{1}{\lambda^2}+C_1\Big(\frac{1+\ln(\frac{1}{\delta_0})}{\delta_0^2 k^2}+\frac{\ln(\frac{1}{\delta_0})}{\delta_0 k}\Big)\theta^\frac{q_1-1}{q_1}\omega_\dimN^{-\frac{1}{q_2}}\rd^{\dimN\kappa}\right)|\mK_\rd|
\end{align}
for $t\in(\tau,t_0)$. Moreover, noticing that from the construction of $\xi$ (see \eqref{eq:xi-def}) we have $\xi\big(w(x,t)\big)\geq \ln\big(\frac{1}{2\delta_0}\big)$ a.e. on $\{w>\mu_+-\delta_0k\}$, we may estimate the left hand side from below against
\begin{align}\label{eq:exp-in-time-eq4}
\int_{\mK_\rd}\xi^2\big(w(\cdot,t)\big)\zeta^2\geq \ln^2\Big(\frac{1}{2\delta_0}\Big)\big|\big\{x\in \mK_{(1-\lambda)\rd}\,:\, w(x,t)>\mu^+-\delta_0 k\big\}\big|\quad\text{for }t\in(\tau,t_0),
\end{align}
because $\zeta=1$ in $\mK_{(1-\lambda)\rd}$ and because $\ln\big(\frac{1}{2\delta_0}\big)$ is positive due to the restriction $\delta_0<\frac12$. Next, we notice that by Bernoulli's inequality we have
\begin{align*}
\big|\mK_{\rd}\setminus \mK_{(1-\lambda)\rd}\big|=|\mK_\rd|-|\mK_{(1-\lambda)\rd}|=\big(1-(1-\lambda)^\dimN\big)|\mK_\rd|\leq \dimN\lambda|\mK_\rd|
\end{align*}
and, accordingly,
\begin{align*}
\big|A_{\mu^+-\delta_0 k,\rd}^+(t)\big|&=\big|\big\{x\in \mK_{\rd}\,:\, w(x,t)>\mu^+-\delta_0 k\big\}\big|\\
&\leq\big|\big\{x\in \mK_{(1-\lambda)\rd}\,:\, w(x,t)>\mu^+-\delta_0 k\big\}\big|+\big|\mK_{\rd}\setminus \mK_{(1-\lambda)\rd}\big|\\
&\leq\big|\big\{x\in \mK_{(1-\lambda)\rd}\,:\, w(x,t)>\mu^+-\delta_0 k\big\}\big|+\dimN\lambda|\mK_\rd|.
\end{align*}
Plugging \eqref{eq:exp-in-time-eq3} and \eqref{eq:exp-in-time-eq4} into this, we arrive at
\begin{align}\label{eq:exp-in-time-eq5}
&\big|A_{\mu^+-\delta_0 k,\rd}^+(t)\big|\\
\leq\ &\left(\frac{\ln^2(\frac{1}{\delta_0})(1-2\gamma)}{\ln^2(\frac{1}{2\delta_0})(1-\gamma)}+\frac{96C_2\ln(\frac{1}{\delta_0})}{\ln^2(\frac{1}{2\delta_0})\lambda^2}+\frac{C_1\theta_M^{-1}}{\ln^2(\frac{1}{2\delta_0})}\Big(\frac{1+\ln(\frac{1}{\delta_0})}{\delta_0^2 k^2}+\frac{\ln(\frac{1}{\delta_0})}{\delta_0 k}\Big)\theta^\frac{q_1-1}{q_1}\omega_\dimN^{-\frac{1}{q_2}}\rd^{\dimN\kappa}+\dimN\lambda\right)|\mK_\rd|\nonumber
\end{align}
for all $t\in(\tau,t_0]$. Clearly, assuming $\delta_0^2k^2\theta^{-\frac{q_1-1}{q_1}}>\rd^{\dimN\kappa}$ to hold, we find that the facts that $\tau \leq t_0-\tfrac{\sig_0}{2}\theta\rd^2$ and $k\leq M$ as well as our choices for $\lambda,\delta_0$ and $j$ finally entail
\begin{align*}
\big|A_{\mu^+-\delta_0 k,\rd}^+(t)\big|\leq\big(1-4\gamma^2+\gamma^2+\gamma^2+\gamma^2\big)|\mK_\rd|=\big(1-\gamma^2\big)|\mK_\rd|=\Big(1-\frac{\sig_0^2}{4}\Big)|\mK_\rd|
\end{align*}
for all $t\in\big(t_0-\tfrac{\sig_0}{2}\theta\rd^2,t_0\big]$ as claimed.
\end{bew}

As an immediate consequence we also obtain the following. 
\begin{corollary}\label{cor:lvl-set-excl}
Let $M>0$ and denote by $\sig_0=\sig_0(\Phi,M,q_1,q_2,a,\hat{b})\in(0,\frac{1}{2})$ and $\delta_0=\delta_0(\Phi,M,q_1,q_2,a,\hat{b})\in(0,\frac12)$ the numbers from Lemma~\ref{lem:first-alt-1} and Lemma~\ref{lem:expansion-in-time}, respectively. Let $w$ be a nonnegative bounded weak solution of \eqref{eq:trans-gen-eq} satisfying \eqref{eq:linfty-bound-w}. Given $x_0\in\Omega$, $0<t_0<T$, $0<\vartheta$ and $0<\rd<\min\big\{\frac{1}{\sqrt{\dimN}}\dist(x_0,\romega),\sqrt{\frac{t_0}{4\vartheta}},1\big\}$, set $\mu^-:=\essinf_{Q_{2\rd}^{\vartheta}}w<\esssup_{Q_{2\rd}^{\vartheta}}w=:\mu^+$, $k\in(0,\mu^+]$ and $\theta:=\min\big\{\theta_k,\theta_{\mu^-}\big\}$. Assume that \eqref{eq:exp-in-time} and \eqref{eq:exp-in-time1} hold. Then for all $\delta\in(0,\delta_0)$ and all $t\in\big(t_0-\tfrac{\sig_0}{2}\theta\rd^2,t_0\big]$
\begin{align*}
\big|\mK_\rd\setminus A^+_{\mu^+-\delta k,\rd}(t)\big|\geq \frac{\sig_0^2}{4} |\mK_\rd|.
\end{align*}
\end{corollary}

\begin{bew}
By definition, $A^+_{\mu^+-\delta k,\rd}(t)\subseteq \mK_\rd$ for all $\delta\in(0,\delta_0)$ and $t\in(0,T)$, so that \[\big|\mK_\rd\setminus A^+_{\mu^+-\delta k,\rd}(t)\big|=|\mK_\rd|-\big|A^+_{\mu^+-\delta k,\rd}(t)\big|\geq |\mK_\rd|-\Big(1-\frac{\sig_0^2}{4}\Big)|\mK_\rd|=\frac{\sig_0^2}{4}|\mK_\rd|\] 
for all $t\in\big(t_0-\tfrac{\sig_0}{2}\theta\rd^2,t_0\big]$ is an evident consequence of Lemma~\ref{lem:expansion-in-time}.
\end{bew}

The corollary above when combined with De Giorgi's inequality (see Lemma~\ref{lem:DeGiorgi}) now provides first quantifiable information on the superlevel set also in time direction.
\begin{lemma}\label{lem:expansion-in-space}

Let $M>0$ and denote by $\sig_0=\sig_0(\Phi,M,q_1,q_2,a,\hat{b})\in(0,\frac{1}{2})$ and $\delta_0=\delta_0(\Phi,M,q_1,q_2,a,\hat{b})\in(0,\frac12)$ the numbers from Lemma~\ref{lem:first-alt-1} and Lemma~\ref{lem:expansion-in-time}, respectively. For any $\sig\in(0,1)$ one can find $\dels=\dels(\Phi,M,q_1,q_2,a,\hat{b},\sig)\in(0,\tfrac{\delta_0}{2}]$ such that for all nonnegative bounded weak solutions $w$ of \eqref{eq:trans-gen-eq} satisfying \eqref{eq:linfty-bound-w} the following holds: Given $x_0\in\Omega$, $0<t_0<T$, $0<\vartheta$ and $0<\rd<\min\big\{\frac{1}{\sqrt{\dimN}}\dist(x_0,\romega),\sqrt{\frac{t_0}{4\vartheta}},1\big\}$ set $\mu^-:=\essinf_{Q_{2\rd}^{\vartheta}}w<\esssup_{Q_{2\rd}^{\vartheta}}w=:\mu^+$, $k\in(0,\mu^+]$ and $\theta:=\min\big\{\theta_k,\theta_{\mu^-}\big\}$. If
\begin{align*}
\theta\leq \vartheta,\quad\text{and}\quad 2k<\mu^+-\mu^-<4k,\quad\text{and}\quad\big|\big\{(x,t)\in Q_\rd^{\theta}\,:\,(w-\mu^++k)_+>0\big\}\big|\leq (1-\sig_0)\big|Q_\rd^{\theta}\big|,
\end{align*}
then either
\begin{align}\label{eq:lem-exp-in-sp-dels}
\dels^2k^2\theta^{-\frac{q_1-1}{q_1}}\leq \rd^{\dimN\kappa}
\end{align}
or
\begin{align*}
\Big|\Big\{(x,t)\in Q_{\rd}^{\frac{\sig_0}{2}\theta}\,:\,(w-\mu^++\dels k)_+>0\Big\}\Big|\leq \sig\big|Q_{\rd}^{\frac{\sig_0}{2}\theta}\big|.
\end{align*}
\end{lemma}

\begin{bew}
With $\delta_0\in(0,\frac{1}{2})$ provided by Lemma~\ref{lem:expansion-in-time}, we claim that the assertion holds if we take $\dels=\delta_0 2^{-j_\star}$ with sufficiently large $j_\star=j_\star(\Phi,M,q_1,q_2,a,\hat{b},\sig)\in\N$ to be determined later. For $j\in\N_0$ we consider the numbers $\delta_j:=\frac{\delta_0}{2^j}$. Obviously, $\delta_j\leq \delta_0$ for all $j\in\N_0$ and $\delta_j-\delta_{j+1}=\delta_{j+1}=\frac{\delta_0}{2^{j+1}}$ for $j\in\N_0$. To shorten the notation, we set $\gamma:=\frac{\sig_0}{2}$ and introduce $A_j^+(t):=A_{\mu^+-\delta_j k,\rd}^+(t)$ for $t\in(t_0-\gamma\theta\rd^2,t_0]$. Employing the De Giorgi inequality from Lemma~\ref{lem:DeGiorgi}, with $L=\mu^+-\delta_{j+1} k$ and $K=\mu^+-\delta_j k$, to the function $x\mapsto w(x,t)$ entails that there is $C_1>0$ such that
\begin{align*}
\delta_{j+1}k\big|A_{j+1}^+(t)\big|=(\delta_j-\delta_{j+1})k\big|A_{j+1}^+(t)\big|\leq \frac{C_1\rd^{\dimN+1}}{\big|\mK_\rd\setminus A_j^+(t)\big|}\int_{A_j^+(t)\setminus A_{j+1}^+(t)}\big|\nabla w(\cdot,t)\big|
\end{align*}
holds for all $t\in(t_0-\gamma\theta\rd^2,t_0]$ and $j\in\N_0$. The denominator on the right can be estimated using Corollary~\ref{cor:lvl-set-excl}, so that integrating the resulting inequality over $(t_0-\gamma\theta\rd^2,t_0]$ and squaring yields
\begin{align*}
\delta_{j+1}^2 k^2\Big(\int_{t_0-\gamma\theta\rd^2}^{t_0}\big|A_{j+1}^+(t)\big|\Big)^2&\leq \Big(\frac{C_1\rd^{\dimN+1}}{\gamma^2|\mK_\rd|}\Big)^2\Big(\int_{t_0-\gamma\theta\rd^2}^{t_0}\!\int_{A_j^+(t)\setminus A_{j+1}^+(t)}\big|\nabla w(\cdot,t)\big|\Big)^2\\
&=\frac{C_1^2\rd^2}{\gamma^4\omega_\dimN^2}\Big(\int_{t_0-\gamma\theta\rd^2}^{t_0}\!\int_{A_j^+(t)\setminus A_{j+1}^+(t)}\big|\nabla w(\cdot,t)\big|\Big)^2
\end{align*}
for all $j\in\N_0$. Invoking Hölder's inequality to split the last integral, we then arrive at the estimate
\begin{align}\label{eq:exp-in-sp-eq1}
&\delta_{j+1}^2 k^2 \Big(\int_{t_0-\gamma\theta\rd^2}^{t_0}\big|A_{j+1}^+(t)\big|\Big)^2\nonumber\\
\leq\ &\frac{C_1^2\rd^2}{\gamma^4\omega_\dimN^2}\Big(\int_{t_0-\gamma\theta\rd^2}^{t_0}\!\int_{A_j^+(t)\setminus A_{j+1}^+(t)}\big|\nabla w(\cdot,t)\big|^2\Big)\Big(\int_{t_0-\gamma\theta\rd^2}^{t_0}\big|A^+_j(t)\setminus A_{j+1}^+(t)\big|\Big)
\end{align}
for all $j\in\N_0$. To control the gradient integral, we first note that with $\psi\in C_0^\infty(\OmT)$ satisfying $\psi\equiv 1$ on $Q_{\rd}^{\gamma\theta}$, $\psi\equiv 0$ on $\Omega_T\setminus Q_{2\rd}^{\gamma\theta}$, $|\nabla\psi|\leq\frac{4}{\rd}$ and $\psi_t\leq \frac{2}{3\gamma\theta\rd^2}$ we obtain 
\begin{align*}
\int_{t_0-\gamma\theta\rd^2}^{t_0}\!\int_{A_j^+(t)\setminus A_{j+1}^+(t)}\big|\nabla w(\cdot,t)\big|^2&\leq \iint_{Q_{\rd}^{\gamma\theta}}\Big|\nabla\Big(\big(w(\cdot,t)-\mu^++\delta_j k\big)_+\psi\Big)\Big|^2\\
&\leq \iint_{Q_{2\rd}^{\gamma\theta}}\Big|\nabla\Big(\big(w(\cdot,t)-\mu^++\delta_j k\big)_+\psi\Big)\Big|^2\quad\text{for all }j\in\N_0,
\end{align*}
so that drawing on Lemma~\ref{lem:local-energy} for $\Lambda^+(w)=\Lambda^+(w,\mu^+-\delta_j k)$ in the cylinder $Q^{\gamma\theta}_{2\rd}$ with $\psi$ as described above provides $C_2=C_2(\Phi,M,q_1,q_2,a,\hat{b})>0$ such that
\begin{align}\label{eq:exp-in-sp-eq2}
&\int_{t_0-\gamma\theta\rd^2}^{t_0}\!\int_{A_j^+(t)\setminus A_{j+1}^+(t)}\big|\nabla w(\cdot,t)\big|^2\nonumber\\
\leq\ &2\iint_{Q_{2\rd}^{\gamma\theta}}\Lambda^+(w)\big(\psi^2\big)_t+3\iint_{Q_{2\rd}^{\gamma\theta}}(w-\mu^++\delta_j k)_+^2|\nabla\psi|^2\nonumber\\&\hspace*{3.2cm}+C_2\Big(\int_{t_0-4\gamma\theta\rd^2}^{t_0}\!\Big|\big\{x\in \mK_{2\rd}:\big(w(x,t)-\mu^++\delta_j k\big)_+>0\big\}\Big|^\frac{q_1(q_2-1)}{q_2(q_1-1)}\Big)^\frac{q_1-1}{q_1}
\end{align}
for all $j\in\N_0$. Since $\gamma\theta<\theta\leq\vartheta$ we have $Q_{2\rd}^{\gamma\theta}\subseteq Q_{2\rd}^{\vartheta}$ and hence $w\leq \mu^+$ in $Q_{2\rd}^{\gamma\theta}$, so that
\begin{align*}
(w-\mu^++\delta_j k)_+^2\leq \delta_j^2 k^2\quad\text{a.e. on }Q_{2\rd}^{\gamma\theta}.
\end{align*}
Moreover noticing that by the assumption $2k\leq \mu^+-\mu^-$ we have $\mu^+-\delta_j k>\mu^+-k\geq \mu^-+k$, we may estimate
\begin{align*}
\frac{\theta_{\mu^+-\delta_j k}}{\theta_k}=\frac{\Phi'\big(\Phi^{-1}(k)\big)}{\Phi'\big(\Phi^{-1}(\mu^+-\delta_j k)\big)}\leq \frac{\Phi'\big(\Phi^{-1}(k)\big)}{\Phi'\big(\Phi^{-1}(\mu^-+k)\big)}\leq 1\quad\text{and}\quad \frac{\theta_{\mu^+-\delta_j k}}{\theta_{\mu^-}}\leq\frac{\Phi'\big(\Phi^{-1}(\mu^-)\big)}{\Phi'\big(\Phi^{-1}(\mu^-+k)\big)}\leq 1,
\end{align*}
which entails that in both cases for $\theta$ we have $\theta_{\mu^+-\delta_j k}\theta^{-1}\leq 1$. Thus, in view of the upper bound for $\Lambda^+$ from Lemma~\ref{lem:Lambda} we obtain
\begin{align*}
\frac{\Lambda^+(w)}{\theta}&\leq \frac{\theta_{\mu^+-\delta_j k}}{2\theta}(w-\mu^++\delta_j k)_+^2\leq (w-\mu^++\delta_j k)_+^2\leq \delta_j^2 k^2\quad\text{a.e. on }Q_{2\rd}^{\gamma\theta}.
\end{align*}
Making use of these estimates, the properties of $\psi$ and \eqref{eq:cond-kappa}, we can further refine \eqref{eq:exp-in-sp-eq2} to
\begin{align}\label{eq:exp-in-sp-eq3}
\int_{t_0-\gamma\theta\rd^2}^{t_0}\!\int_{A_j^+(t)\setminus A_{j+1}^+(t)}\!\!\big|\nabla w(\cdot,t)\big|^2&\leq\frac{8\delta_j^2 k^2}{3\gamma\rd^2}\big|Q_{2\rd}^{\gamma\theta}\big|
+\frac{48\delta_j^2 k^2}{\rd^2}\big|Q_{2\rd}^{\gamma\theta}\big|+C_2|\mK_{2\rd}|^{\frac{q_2-1}{q_2}}\big|4\gamma\theta\rd^2\big|^\frac{q_1-1}{q_1}\nonumber\\
&=\frac{\delta_j^2 k^2}{\rd^2}\Big(\frac{8}{3\gamma}+48\Big)\big|Q_{2\rd}^{\gamma\theta}\big|+C_2 2^{\dimN+2-\frac{\dimN}{q_2}-\frac{2}{q_1}}\omega_\dimN^{1-\frac{1}{q_2}}\rd^{\dimN+\dimN\kappa}\gamma^{1-\frac{1}{q_1}}\theta^{\frac{q_1-1}{q_1}}\nonumber\\
&=\frac{\delta_j^2 k^2}{\rd^2}\Big(\frac{8}{3\gamma}+48+\frac{C_2\rd^{\dimN\kappa}\theta^{\frac{q_1-1}{q_1}}}{\delta_j^2 k^2}2^{\dimN\kappa-2}\omega_\dimN^{-\frac{1}{q_2}}\theta^{-2}\gamma^{-\frac{1}{q_1}}\Big)\big|Q_{2\rd}^{\gamma\theta}\big|
\end{align}
for all $j\in\N_0$. Now, assuming \eqref{eq:lem-exp-in-sp-dels} to be false we have $\rd^{\dimN\kappa}<\dels^2 k^2\theta^{-\frac{q_1-1}{q_1}}$, so that we conclude from  the fact that due to the monotonicity of $\Phi'$ and $\Phi^{-1}$ and $\mu^-\leq M$ and $k\leq M$ we have $\theta^{-2}\leq\theta_M^{-2}$ and \eqref{eq:exp-in-sp-eq3}, that there is $C_3=C_3(\Phi,M,q_1,q_2,a,\hat{b})>0$ such that
\begin{align*}
\int_{t_0-\gamma\theta\rd^2}^{t_0}\!\int_{A_j^+(t)\setminus A_{j+1}^+(t)}\!\!\big|\nabla w(\cdot,t)\big|^2&\leq C_3\frac{\delta_j^2 k^2}{\rd^2}\Big(1+\frac{\dels^2}{\delta_j^2}\Big)\big|Q_{2\rd}^{\gamma\theta}\big|
\end{align*}
for all $j\in\N_0$. Plugging this back into \eqref{eq:exp-in-sp-eq1} and reordering yields
\begin{align*}
\Big(\int_{t_0-\gamma\theta\rd^2}^{t_0}\big|A_{j+1}^+(t)\big|\Big)^2\leq \frac{C_1^2\rd^2}{\gamma^4\omega_\dimN^2}\cdot C_3\frac{\delta_j^2 k^2}{\delta_{j+1}^2 k^2 \rd^2}\Big(1+\frac{\dels^2}{\delta_j^2}\Big)\big|Q_{2\rd}^{\gamma\theta}\big|\Big(\int_{t_0-\gamma\theta\rd^2}^{t_0}\big|A_j^+(t)\setminus A_{j+1}^+(t)\big|\Big)
\end{align*}
for all $j\in\N_0$. Then, since by definition we have $\delta_{j}\geq\dels$ for all $j\leq j_\star$ and $\frac{\delta_j}{\delta_{j+1}}=2$ for all $j\leq j_\star-1$, we find that there is $C_4=C_4(\Phi,M,q_1,q_2,a,\hat{b})>0$ such that  
\begin{align}\label{eq:exp-in-sp-eq4}
\Big(\int_{t_0-\gamma\theta\rd^2}^{t_0}\big|A_{j+1}^+(t)\big|\Big)^2\leq C_4\big|Q_{2\rd}^{\gamma\theta}\big|\Big(\int_{t_0-\gamma\theta\rd^2}^{t_0}\big|A_j^+(t)\setminus A_{j+1}^+(t)\big|\Big)\quad\text{for all }0\leq j\leq j_\star -1
\end{align}
and a summation of \eqref{eq:exp-in-sp-eq4} from $j=0$ to $j=j_\star-1$ entails
\begin{align*}
&j_\star\big|\big\{(x,t)\in Q_{\rd}^{\gamma\theta}\,:\, (w-\mu^++\dels k)_+>0\big\}\big|^2\\
=\ & j_\star\Big(\int_{t_0-\gamma\theta\rd^2}^{t_0}\big|A_{j_\star}^+(t)\big|\Big)^2\\
\leq\ &C_4\big|Q_{2\rd}^{\gamma\theta}\big|\cdot\sum_{j=0}^{j_\star-1}\int_{t_0-\gamma\theta\rd^2}^{t_0}\big|A_j^+(t)\setminus A_{j+1}^+(t)\big|\\
\leq\ &C_4|Q_{2\rd}^{\gamma\theta}\big|\cdot\big|\big\{(x,t)\in Q_{\rd}^{\gamma\theta}\,:\, (w-\mu^++\delta_0 k)_+>0\big\}\big|\\
\leq\ &C_4|Q_{2\rd}^{\gamma\theta}\big|\cdot\big|Q_{\rd}^{\gamma\theta}\big|
\end{align*}
in view of the nested inclusions $A_{j+1}^+\subseteq A_j^+(t)$ for all $j\leq j_\star-1$ and $t\in(t_0-\gamma\theta\rd^2,t_0]$. Rewriting $|Q_{2\rd}^{\gamma\theta}|=2^{\dimN+2}|Q_{\rd}^{\gamma\theta}|$ and dividing by $j_\star$ we obtain
\begin{align*}
\big|\big\{(x,t)\in Q_{\rd}^{\gamma\theta}:(w-\mu^++\dels k)_+>0\big\}\big|^2\leq \frac{C_4 2^{\dimN+2}}{j_\star}|Q_{\rd}^{\gamma\theta}|^2.
\end{align*}
Accordingly, our claim is proven if we take $j_\star=j_\star(\Phi,M,q_1,q_2,a,\hat{b},\sig)\in\N$ so large that $\sqrt{\frac{C_4 2^{\dimN+2}}{j_\star}}\leq \sig.$
\end{bew}

Finally, we can now state the complete second alternative.
\begin{lemma}\label{lem:second-alt}
Let $M>0$ and denote by $\sig_0=\sig_0(\Phi,M,q_1,q_2,a,\hat{b})\in(0,\frac{1}{2})$ the number from Lemma~\ref{lem:first-alt-1}. There is $\deld=\deld(\Phi,M,q_1,q_2,a,\hat{b})\in(0,\frac{1}{4})$ such that for all nonnegative bounded weak solutions $w$ of \eqref{eq:trans-gen-eq} satisfying \eqref{eq:linfty-bound-w} the following holds: Given $x_0\in\Omega$, $0<t_0<T$, $0<\vartheta$ and $0<\rd<\min\big\{\frac{1}{\sqrt{\dimN}}\dist(x_0,\romega),\sqrt{\frac{t_0}{4\vartheta}},1\big\}$, set $\mu^-:=\essinf_{Q_{2\rd}^{\vartheta}}w<\esssup_{Q_{2\rd}^{\vartheta}}w=:\mu^+$, $k\in(0,\mu^+]$ and $\theta:=\min\big\{\theta_k,\theta_{\mu^-}\big\}$. If
\begin{align*}
\theta\leq\vartheta ,\quad\text{and}\quad 2k\leq \mu^+-\mu^-\leq 4k,\quad\text{and}\quad\big|\big\{(x,t)\in Q_\rd^{\theta}:(w-\mu^++k)_+>0\big\}\big|\leq (1-\sig_0)\big|Q_\rd^{\theta}\big|,
\end{align*}
then either
\begin{align}\label{eq:sec-alt-deld}
\deld^2k^2\theta^{-\frac{q_1-1}{q_1}}\leq \rd^{\dimN\kappa}
\end{align}
or
\begin{align}\label{eq:sec-alt-zero}
\Big|\Big\{(x,t)\in Q_{\frac{\rd}{2}}^{\frac{\sig_0}{2}\theta}:(w-\mu^++\tfrac{\deld}{2} k)_+>0\Big\}\Big|=0.
\end{align}
\end{lemma}

\begin{bew}
We repeat the iteration procedure of Lemma~\ref{lem:first-alt-1} in an adjusted manner. (See also \cite[Proposition 4.5]{ChungHwangKangKim17}.) 
We claim that the assertion holds if we pick $\sig_1\in(0,1)$ sufficiently small with our choice only depending on $\Phi,M,a$ and $\hat{b}$, and then let $\deld=\deld(\Phi,M,q_1,q_2,a,\hat{b})=\dels(\Phi,M,q_1,q_2,a,\hat{b},\sig_1)\!\in(0,\frac{1}{4})$, where $\dels$ is the number provided by Lemma~\ref{lem:expansion-in-space}. We set $\gamma:=\frac{\sig_0}{2}$, $k_0:=\mu^+-\deld k$, $\rd_0=\rd$ and define
\begin{align*}
k_i:=\mu^+-\frac{\deld k}{2}-\frac{\deld k}{2^{i+1}},\quad\text{and}\quad \rd_i=\frac{\rd}{2}+\frac{\rd}{2^{i+1}}.
\end{align*}
Accordingly, for all $i\in\N_0$ we have $k_i\in[\mu^+-\deld k,\mu^+-\frac{\deld k}{2})$, $\rd_i\in(\frac{\rd}{2},\rd]$, $k_i\leq k_{i+1}$ and, moreover, $\rd_\infty:=\lim_{i\to\infty} \rd_i=\frac{\rd}{2}$ and $k_\infty:=\lim_{i\to\infty} k_i=\mu^+-\frac{\deld k}{2}$. Additionally, we see that the assumption $\mu^+-\mu^-\geq 2k$ implies that $k_0>\mu^+-k\geq \mu^-+k$. For $i\in\N_0$ we denote  $Q_i:=Q_{\rd_i}^{\gamma\theta}=\mK_{\rd_i}(x_0)\times(t_0-\gamma\theta\rd_i^2,t_0]$ and fix cutoff functions similar to those in Lemma~\ref{lem:first-alt-1}. In fact, we take $\psi_i\in C_0^\infty\big(\OmT\big)$ with $\supp \psi_i\subseteq \overline{Q_i}$ such that
\begin{align*}
\psi_i\equiv 1\ \text{in}\ Q_{i+1},\qquad \psi_i\equiv 0\text{ in }\partial_p Q_i,\qquad |\nabla\psi_i|\leq \frac{2^{i+3}}{\rd	}\qquad\text{and}\qquad |(\psi_i)_t|\leq\frac{2^{2(i+2)}}{\gamma\theta\rd^2}.
\end{align*}
This way, $Q_i\subseteq Q_{\rd_i}^{\theta}\subseteq Q_{2\rd}^{\vartheta}$ for all $i\in\N_0$ and hence $w\leq \mu^+$ in $Q_i$ for all $i\in\N_0$. We proceed by employing Lemma~\ref{lem:local-energy} for $\Lambda^+(w)=\Lambda^+(w,k_i)$ in $Q_i$ and obtain $\Gamma_1=\Gamma_1(\Phi,M,q_1,q_2,a,\hat{b})>0$ satisfying
\begin{align*}
&\sup_{t_0-\gamma\theta\rd_i^2\leq t \leq t_0}2\int_{\mK_{\rd_i}}\Lambda^+\big(w(\cdot,t)\big)\psi_i^2(\cdot,t)+\iint_{Q_i}\big|\nabla\big((w-k_i)_+\psi_i\big)\big|^2\\
\leq\ &4\iint_{Q_{i}}\Lambda^+(w)\psi_i|(\psi_i)_t|+3\iint_{Q_i}(w-k_i)_+^2|\nabla\psi_i|^2+\Gamma_1\Big(\int_{t_0-\gamma\theta\rd_i^2}^{t_0}\!\big|A_{k_i,\rd_i}^+(t)\big|^\frac{q_1(q_2-1)}{q_2(q_1-1)}\Big)^\frac{q_1-1}{q_1}
\end{align*}
for all $i\in\N_0$. Recalling that in view of Lemma~\ref{lem:Lambda} we have
\begin{align*}
\frac{\theta_{\mu^+}}{2}(w-k_i)_+^2=\frac{(w-k_i)_+^2}{2\Phi'\big(\Phi^{-1}(\mu^+)\big)}\leq \Lambda^+(w)\leq \frac{(w-k_i)_+^2}{2\Phi'\big(\Phi^{-1}(k_i)\big)}=\frac{\theta_{k_i}}{2}(w-k_i)_-^2\quad\text{in }Q_i
\end{align*}
yields
\begin{align*}
&J_i:=\sup_{t_0-\gamma\theta\rd_i^2\leq t \leq t_0}\int_{\mK_{\rd_i}}\big(w(\cdot,t)-k_i)_+^2\psi_i^2(\cdot,t)+\frac{1}{\theta_{\mu^+}}\iint_{Q_i}\big|\nabla\big((w-k_i)_+\psi_i\big)\big|^2\nonumber\\
&\qquad\leq \frac{2^{2i+5}}{\gamma\rd^2}\frac{\theta_{k_i}}{\theta_{\mu^+}\theta}\iint_{Q_{i}}(w-k_i)_+^2+\frac{2^{2i+8}}{\theta_{\mu^+}\rd^2}\iint_{Q_i}(w-k_i)_+^2+\frac{\Gamma_1}{\theta_{\mu^+}}\Big(\int_{t_0-\gamma\theta\rd_i^2}^{t_0}\!\big|A_{k_i,\rd_i}^+(t)\big|^\frac{q_1(q_2-1)}{q_2(q_1-1)}\Big)^\frac{q_1-1}{q_1}
\end{align*}
for all $i\in\N_0$ in view of the monotonicity of $\Phi'$, $\Phi^{-1}$ and the properties of $\psi_i$. As $\frac{\theta}{\theta_{\mu^+}}\geq 1$ due to $\mu^+\geq k$ and $\mu^+\geq \mu^-$, we obtain from repeating the change of variables $\bar{t}=\frac{t-t_0}{\theta}$ as in Lemma~\ref{lem:first-alt-1} and denoting $\bar{Q}_i:=\mK_{\rd_i}\times(-\gamma\rd_i^2,0]$ and $\bar{A}^+_i(\bar{t}):=A_{k_i,\rd_i}^+(\theta\bar{t}+t_0)$ that
\begin{align}\label{eq:sec-alt-eq1}
&\bar{J}_i:=\sup_{-\gamma\rd_i^2\leq \bar{t} \leq 0}\int_{\mK_{\rd_i}}\big(w(\cdot,\bar{t})-k_i)_+^2\psi_i^2(\cdot,\bar{t})+\iint_{\bar{Q}_i}\big|\nabla\big((w-k_i)_+\psi_i\big)\big|^2\nonumber\\
&\qquad\leq \frac{2^{2i+5}}{\gamma\rd^2}\frac{\theta_{k_i}}{\theta_{\mu^+}}\iint_{\bar{Q}_{i}}(w-k_i)_+^2+\frac{2^{2i+8}}{\rd^2}\frac{\theta}{\theta_{\mu^+}}\iint_{\bar{Q}_i}(w-k_i)_+^2+\Gamma_1\frac{\theta^{\frac{q_1-1}{q_1}}}{\theta\theta_{\mu^+}}\Big(\int_{-\gamma\rd_i^2}^{0}\!\big|\bar{A}_{i}^+(\bar{t})\big|^\frac{q_1(q_2-1)}{q_2(q_1-1)}\intd\bar{t}\Big)^\frac{q_1-1}{q_1}
\end{align}
for all $i\in\N_0$. Here, we note that our choice for $k_i$ and the assumption $\mu^+-\mu^-\geq 2k$ entail that
\begin{align*}
k_i\geq k_0=\mu^+-\deld k> \mu^+- k\geq\mu^-+k\geq\max\{\mu^-,k\}\quad\text{for all }i\in\N_0,
\end{align*}
so that in both instances for $\theta$ we have $\theta_{k_i}\leq \theta$ due to the monotonicity properties of $\Phi'$ and $\Phi^{-1}$. Moreover, the condition $\mu^+-\mu^-\leq 4k$ entails that $\mu^+\leq 5\max\{k,\mu^-\}\leq 5M$, so that in either case we may conclude from Lemma~\ref{lem:phi}, again utilizing the monotonicity of $\Phi'$ and $\Phi^{-1}$ as well as the subadditivity of $\Phi^{-1}$, that there is $C_1=C_1(\Phi,M)>0$ such that
\begin{align*}
\frac{\theta_{k_i}}{\theta_{\mu^+}}\leq\frac{\theta}{\theta_{\mu^+}}\leq \frac{\Phi'\big(\Phi^{-1}(5\max\{k,\mu^-\})\big)}{\Phi'\big(\Phi^{-1}(\max\{k,\mu^-\})\big)}\leq \frac{\Phi'\big(5\Phi^{-1}(\max\{k,\mu^-\})\big)}{\Phi'\big(\Phi^{-1}(\max\{k,\mu^-\})\big)}\leq C_1\quad\text{for all }i\in\N_0.
\end{align*}
Plugging this back into \eqref{eq:sec-alt-eq1}, estimating $\theta^{-1}\leq\theta_{\mu^+}^{-1}\leq\theta_M^{-1}$ as well as $(w-k_i)_+\leq \deld k$ for all $i\in\N_0$ and introducing the notation $\bar{A}_i^+:=\bar{Q}_i\cap\big\{(x,t)\in\Omega\times(-\frac{t_0}{\theta},\frac{T-t_0}{\theta}):\big(w(x,\bar{t})-k_i\big)_+>0\big\}$, we find $\Gamma_2=\Gamma_2(\Phi,M,q_1,q_2,a,\hat{b})>0$ fulfilling
\begin{align}\label{eq:sec-alt-eq2}
\bar{J}_i\leq \Gamma_2\frac{2^{2i}}{\rd^2}\Big(\frac{1}{\gamma}+1\Big)\deld^2k^2\big|\bar{A}_i^+\big|+\Gamma_2\theta^{\frac{q_1-1}{q_1}}\Big(\int_{-\gamma\rd_i^2}^{0}\big|\bar{A}_{i}^+(\bar{t})\big|^\frac{q_1(q_2-1)}{q_2(q_1-1)}\intd{\bar{t}}\Big)^\frac{q_1-1}{q_1}
\end{align}
for all $i\in\N_0$. Repeating steps similar to those in Lemma~\ref{lem:first-alt-1}, while drawing on Corollary~\ref{cor:sob} and the fact that $(w-k_i)_+\geq k_{i+1}-k_i=\frac{\deld k}{2^{i+2}}$ a.e. in $\bar{Q}_{i+1}$, we moreover see that there is $C_2>0$ such that
\begin{align}\label{eq:sec-alt-eq3}
\big|\bar{A}_{i+1}^+\big|\leq \frac{C_2 2^{2i+4}}{\deld^2 k^2}\big|\bar{A}_i^+\big|^\frac{2}{\dimN+2}\bar{J}_i\quad\text{for }i\in\N_0.
\end{align}
Setting $\widehat{q}_j:=\frac{2q_j(1+\kappa)}{q_j-1}$ for $j\in\{1,2\}$ and $\kappa$ as in \eqref{eq:cond-kappa} as well as
\begin{align*}
\bar{Y}_i:=\frac{\big|\bar{A}_i^+\big|}{\rd^{\dimN+2}}\quad\text{as well as}\quad\bar{Z}_i:=\frac{1}{\rd^\dimN}\Big(\int_{-\gamma\rd_i^2}^0\big|\bar{A}_{i}^+(\bar{t})\big|^\nfrac{\widehat{q_1}}{\widehat{q_2}}\intd{\bar{t}}\Big)^\nfrac{2}{\widehat{q_1}},\quad i\in\N_0,
\end{align*}
we conclude from  \eqref{eq:sec-alt-eq2} and \eqref{eq:sec-alt-eq3} that there is $\Gamma_3=\Gamma_3(\Phi,M,q_1,q_2,a,\hat{b})>0$ satisfying
\begin{align*}
\bar{Y}_{i+1}\leq\Gamma_3 2^{4i}\Big(\frac{1}{\gamma}+1\Big)\bar{Y}_i^{1+\frac{2}{\dimN+2}}+\Gamma_3 2^{4i}\frac{\rd^{\dimN\kappa}\theta^{\frac{q_1-1}{q_1}}}{\deld^2 k^2}\bar{Y}_i^{\frac{2}{\dimN+2}}\bar{Z}_{i}^{1+\kappa}\quad\text{for all }i\in\N_0.
\end{align*}
As the choice for $\gamma=\frac{\sig_0}{2}$ only depends on $\Phi,M,q_1,q_2,a$ and $\hat{b}$, we find that if \eqref{eq:sec-alt-deld} fails there is $\Gamma_4=\Gamma_4(\Phi,M,q_1,q_2,a,\hat{b})>0$ such that
\begin{align*}
\bar{Y}_{i+1}\leq\Gamma_4 2^{4i}\Big(\bar{Y}_i^{1+\frac{2}{\dimN+2}}+\bar{Y}_i^{\frac{2}{\dimN+2}}\bar{Z}_i^{1+\kappa}\Big)
\end{align*}
Mirroring the steps from Lemma~\ref{lem:first-alt-1} for $\bar{Z}_i$ we also infer from Lemma~\ref{lem:lemmaX} that if \eqref{eq:sec-alt-deld} fails there are $C_3>0$ and $\Gamma_5=\Gamma_5(\Phi,M,q_1,q_2,a,\hat{b})>0$ such that
\begin{align*}
\bar{Z}_{i+1}\leq\frac{2^{2i+4}C_3}{\rd^\dimN\deld^2 k^2}\bar{J}_i\leq 2^{4i}\Gamma_5 \bar{Y}_i+2^{4i}\Gamma_5\frac{\rd^{\dimN\kappa}\theta^{\frac{q_1-1}{q_1}}}{\deld^2 k^2}\bar{Z}_i^{1+\kappa}\leq 2^{4i}\Gamma_5\Big(\bar{Y}_i+\bar{Z}_i^{1+\kappa}\Big)\quad\text{for all }i\in\N_0,
\end{align*}
so that once more we are left with checking the condition on $\bar{Y}_0+\bar{Z}_0^{1+\kappa}$ from Lemma~\ref{lem:fast-geo-conv}. One can easily check (compare the steps leading to \eqref{eq:first-alt-eq3} and \eqref{eq:first-alt-eq3.5}) that according to the change in variable and \eqref{eq:cond-kappa} we have
\begin{align}\label{eq:sec-alt-eq4}
\bar{Y}_0+\bar{Z}_0^{1+\kappa}&\leq \frac{\big|\bar{A}_{0}^+\big|}{\rd^{\dimN+2}}+\omega_\dimN^{\max\{\frac{1}{q_1}-\frac{1}{q_2},0\}}\rd^{\frac{\dimN+2}{\min\{q_1,q_2\}}-\dimN-2}\big|\bar{A}_{0}^+\big|^{q}\nonumber\\&\leq\frac{\big|A_{0}^+\big|}{\theta\rd^{\dimN+2}}+\omega_\dimN^{\max\{\frac{1}{q_1}-\frac{1}{q_2},0\}}\rd^{\frac{\dimN+2}{\min\{q_1,q_2\}}-\dimN-2}\theta^{-q}\big|A_{0}^+\big|^q,
\end{align}
where we wrote $q=\min\big\{\frac{q_1-1}{q_1},\frac{q_2-2}{q_2}\big\}$ and $A_0^+=\big\{(x,t)\in Q_\rd^{\gamma\theta}:(w-\mu^++\deld k)_+>0\}$. Since we assume that \eqref{eq:sec-alt-deld} fails, we find from our definition of $\deld=\dels$ and Lemma~\ref{lem:expansion-in-space} that 
\begin{align*}
\big|A_{0}^+\big|=\Big|\Big\{(x,t)\in Q_{\rd}^{\gamma\theta}:(w-\mu^++\deld k)_+>0\Big\}\Big|\leq \sig_1\big|Q_{\rd}^{\gamma\theta}\big|=\sig_1\gamma\theta\omega_\dimN\rd^{\dimN+2},
\end{align*}
so that \eqref{eq:sec-alt-eq4} turns into
\begin{align*}
\bar{Y}_0+\bar{Z}_0^{1+\kappa}\leq \sig_1\gamma\omega_\dimN+\big(\sig_1\gamma\omega_\dimN\big)^{\frac{q_2-1}{q_2}}\rd^{2-\frac{2}{q_2}}\leq\sig_1\gamma\omega_\dimN+\big(\sig_1\gamma\omega_\dimN\big)^{\frac{q_2-1}{q_2}}
\end{align*}
because of $\rd<1$ and $q_2>1$. As $\gamma=\frac{\sig_0}{2}$ and $\sig_0$ only depends $\Phi,M,q_1,q_2,a,\hat{b}$, we can therefore pick $\sig_1=\sig_1(\Phi,M,q_1,q_2,a,\hat{b})>0$ so small that with $\sigma=\min\{\kappa,\frac{2}{\dimN+2}\}$ and $\Gamma_6=\Gamma_6(\Phi,M,q_1,q_2,a,\hat{b}):=\max\{\Gamma_4,\Gamma_5\}$ we have
\begin{align*}
\sig_1\gamma\omega_\dimN+\big(\sig_1\gamma\omega_\dimN\big)^{\frac{q_2-1}{q_2}}\leq \big(2\Gamma_6\big)^{-\frac{1+\kappa}{\sigma}}2^{-\frac{4(1+\kappa)}{\sigma^2}},
\end{align*}
thereby making Lemma~\ref{lem:fast-geo-conv} applicable to say $\bar{Y}_i\to 0$ as $i\to\infty$. This immediately yields
\begin{align*}
\big|\big\{(x,t)\in Q_{\frac{\rd}{2}}^{\gamma\theta}:(w-\mu^++\tfrac{\deld}{2}k)_+>0\big\}\big|&=\theta\big|\big\{(x,\bar{t})\in \bar{Q}_{\infty}:(w-\mu^++\tfrac{\deld}{2}k)_+>0\big\}\big|=\theta\big|\bar{A}^+_\infty\big|=0
\end{align*}
and hence \eqref{eq:sec-alt-zero}.
\end{bew}

\subsection{Constructing a sequence of nested cylinders}\label{sec4:iter}
With the alternatives prepared, we can now arrange the details of the iteration procedure. Starting with a cylinder of appropriate size we will iteratively construct a sequence of nested cylinders $Q_j$ and non-increasing positive numbers $k_j$ related to the oscillation, such that the oscillation of $w$ inside this sequence of cylinders decays with a fixed rate. For technical reasons concerning the scaling factors $\theta_j$, which are chosen according to $k_j$ and the essential infimum $\mu^-_{j}$ inside the current cylinder $Q_j$, the iteration will be split into two parts. For the first part we will assume that $\mu^-_{j}\leq k_j$, i.e. that $\theta=\theta_{k_j}$ in the previous alternatives Lemma~\ref{lem:first-alt-1} and Lemma~\ref{lem:second-alt}. If we would have $\mu^-_{j}=0$ for all $j\in\N$ we could just follows this iteration to completion. However, if $\essinf_{Q_{j_i}} w$ is strictly larger than zero for some $j_i\in\N$ the condition certainly gets violated at some point of our iteration, due to $k_j$ being decreasing. If this is the case we will use this first occurrence of $\mu^-_{j_i}>k_{j_i}$ as starting point for the second iteration procedure, where then we choose $\theta=\theta_{\mu^-_{j}}$ in the alternatives. Since the cylinders are nested, the essential infimum is non-decreasing with each step, so by the decreasing property of the $k_j$ we will always remain in this setting and can, in cases where the first iteration stops, then follow this second iteration to completion.

\begin{lemma}\label{lem:iteration-1}
Let $w$ be a nonnegative bounded weak solution of \eqref{eq:trans-gen-eq} with $M>0$ such that \eqref{eq:linfty-bound-w} holds. There are $\lambda=\lambda(\Phi,M,q_1,q_2,a,\hat{b})\in(0,\tfrac12)$, $\eta=\eta(\Phi,M,q_1,q_2,a,\hat{b})\in(\tfrac12,1)$, $\beta=\beta(q_1,q_2)\in(0,1)$, $C=C(\Phi,M,q_1,q_2,a,\hat{b})\geq 1$ with the following property: Given $x_0\in\Omega$, $0<t_0<T$, $0<k_0\leq\frac{M}{2}$ and $0<\rd_0<\min\big\{\frac{1}{\sqrt{\dimN}}\dist(x_0,\romega),\sqrt{\frac{t_0}{4\theta_{k_0}}},1\big\}$ such that $\essosc_{Q_{\rd_0}^{\theta_{k_0}}}w\leq 4k_0$ one can find a decreasing sequence $(\rd_j)_{j\in\N}$, a non-increasing sequence $(k_j)_{j\in\N}$, and an increasing sequence $(\theta_j)_{j\in\N}$ such that with $Q_j:=Q_{\rd_j}^{\theta_j}$ the properties
\begin{align}\label{eq:iteration-1}
Q_j\subseteq Q_{j-1},\ \; \essosc_{Q_{j}} w\leq 4k_j,\ \; \frac{1}{2}k_{j-1}\leq k_j\leq \max\big\{\eta k_{j-1}, C\rd_{j-1}^\beta\big\}\ \;\text{and}\ \; \lambda\rd_{j-1}\leq\rd_j\leq\frac{1}{2}\rd_{j-1}
\end{align}
either hold for all $j\in\N$ or there is $\js\geq0$ such that \eqref{eq:iteration-1} holds for all $j\leq \js$ and $\mu^{-}_{\js}:=\essinf_{Q_{\js}} w>k_{\js}$.
\end{lemma}

\begin{bew}
With $r_0$, $k_0$ and $\theta_0:=\theta_{k_0}=\big(\Phi'(\Phi^{-1}(k_0))\big)^{-1}$ provided we set $Q_0:=Q_{\rd_0}^{\theta_{0}}$, $\mu^-_{0}:=\essinf_{Q_0} w$ and $\mu^+_{0}:=\esssup_{Q_0}w$. If $\mu^-_{0}>k_0$ we set $j_\star=j_0$ and have nothing to prove. Assume now $\mu^-_{0}\leq k_0$. For $j>0$ we denote by $\mu^-_{j}:=\essinf_{Q_j} w$ and $\mu^+_{j}:=\esssup_{Q_j}w$ the infimum and supremum, respectively, in the current cylinder $Q_j$. If at some point $\mu^-_{j}>k_j$, we are finished. If not, we will generate $\rd_{j+1}$ and $k_{j+1}$ depending on which of the different cases we find ourselves in. The next scaling factor $\theta_{j+1}$ will always be chosen in dependence of $k_{j+1}$, i.e. $\theta_{j+1}:=\theta_{k_{j+1}}$ in all cases, so that in fact $\theta_j=\theta_{k_j}$ for all $j\leq j_\star$. Before giving details for the different cases, however, let us fix $C_1=C_1(\Phi,2,M)>1$ such that
\begin{align}\label{eq:it1-c1}
\Phi'\big(\Phi^{-1}(s)\big)=\Phi'\big(2\cdot\tfrac12\Phi^{-1}(s)\big)\leq C_1\Phi'\big(\tfrac12\Phi^{-1}(s)\big)\quad\text{for all }s\in[0,M],
\end{align}
which is possible in view of Lemma~\ref{lem:phi}. Moreover, let the numbers $\sig_0=\sig_0(\Phi,M,q_1,q_2,a,\hat{b})\in(0,\tfrac12)$ and $\deld=\deld(\Phi,M,q_1,q_2,a,\hat{b})\in(0,\tfrac14)$ be given by Lemma~\ref{lem:first-alt-1} and Lemma~\ref{lem:second-alt}, respectively. We claim that the lemma holds for the choices
\begin{equation}\label{eq:it1-const}
\begin{split}
\lambda&:=\frac{\sqrt{\sig_0}}{4\sqrt{2C_1}}\in(0,\tfrac18),\qquad\eta:=1-\tfrac{\deld}{8}\in(\tfrac{31}{32},1),\qquad \beta:=\frac{\dimN\kappa}{2+\frac{q_1-1}{q_1}}<1\\\text{ and } \  C&:=\max\Big\{1,\Big(\big(\Phi^{-1}(M)\big)^{\frac{q_1-1}{q_1}}\deld^{-2}\Big)^\frac{1}{2+\frac{q_1-1}{q_1}}\Big\}\geq 1.
\end{split}
\end{equation}
Let us now turn to the construction of the next cylinder.
\smallskip

\emph{Case 1}: $\essosc_{Q_j} w< 2k_j$.\smallskip\\
With $C_1>1$ as provided above, we set $r_{j+1}=\frac{r_j}{2\sqrt{C_1}}$, $k_{j+1}=\frac{k_j}{2}$ and $\theta_{j+1}:=\theta_{k_{j+1}}$. Clearly, the conditions for $k_{j+1}$ and $r_{j+1}$ in \eqref{eq:iteration-1} are satisfied in view of $\eta>\tfrac12$ and $\lambda<\tfrac1{2\sqrt{C_1}}$. Concerning the first condition, we note that \eqref{eq:it1-c1}, the monotonicity of $\Phi'$ and the subadditivity of $\Phi^{-1}$ imply that $$\Phi'\big(\Phi^{-1}(k_j)\big)\leq C_1\Phi'\big(\tfrac{1}{2}\Phi^{-1}(k_j)\big)=C_1\Phi'\big(\tfrac{1}{2}\Phi^{-1}(2\cdot\tfrac{k_j}{2})\big)\leq C_1\Phi'\big(\Phi^{-1}(\tfrac{k_j}{2})\big).$$ Thus, we have
\begin{align*}
\rd_{j+1}^2\theta_{{j+1}}=\frac{\rd_j^2}{4C_1\Phi'\big(\Phi^{-1}(\frac{k_j}{2})\big)}\leq \frac{\rd_j^2}{4\Phi'\big(\Phi^{-1}(k_j)\big)}=\frac{\rd_j^2}{4}\theta_{k_j}\leq\rd_j^2\theta_{j},
\end{align*}
so that $Q_{j+1}\subseteq Q_{j}$ follows from our choice $\rd_{j+1}\leq \frac{\rd_j}{2}$. We are left with verifying $\essosc_{Q_{j+1}}w\leq 4 k_{j+1}$. This is an evident consequence of the nesting of cylinders, our case assumption and the particular choice of $k_{j+1}$. In fact, due to $Q_{j+1}\subseteq Q_j$ we have
\begin{align*}
\essosc_{Q_{j+1}} w\leq \essosc_{Q_j} w<2k_j=4 k_{j+1}.
\end{align*}\smallskip

\emph{Case 2}: $2k_j\leq \essosc_{Q_j} w\leq 4k_j$.\smallskip\\
We further distinguish cases indicated by the alternatives in Lemma~\ref{lem:first-alt-1} and Lemma~\ref{lem:second-alt}, which will be employed for $\vartheta=\theta_j$, $\rd=\frac{\rd_j}{2}$ and $k=k_j$. Note that due to the standing assumption that $k_j\geq \mu^-_{j}$ we have $\theta=\theta_{k_j}$ in these lemmas. \smallskip

\emph{Case 2a}: $\deld^2 k_j^2\theta_{k_j}^{-\frac{q_1-1}{q_1}}\leq \big(\tfrac{\rd_j}{2}\big)^{\dimN\kappa}$, i.e. \eqref{eq:sec-alt-deld} is true for $k=k_j$, $\theta=\theta_{k_j}$ and $\rd=\tfrac{\rd_j}{2}$.\smallskip\\
In this case we set $\rd_{j+1}:=\frac{\rd_j}{2}$, $k_{j+1}:=k_j$ and $\theta_{j+1}:=\theta_{k_{j+1}}$. These choices entail that
\begin{align*}
\rd_{j+1}^2\theta_{k_{j+1}}=\frac{\rd_j^2}{4}\theta_{k_j}\leq \rd_j^2\theta_{k_j}\quad\text{and}\quad \rd_{j+1}\leq \rd_{j},
\end{align*}
and hence $Q_{j+1}\subseteq Q_{j}$. Due to this and the choice $k_{j+1}=k_j$ we also immediately obtain $$\essosc_{Q_{j+1}} w\leq\essosc_{Q_j} w\leq 4{k_j}=4 k_{j+1}.$$ Fulfillment of $\frac12 k_{j}\leq k_{j+1}$ and $\lambda r_{j}\leq r_{j+1}\leq \frac{1}{2}\rd_{j}$ are also evident by the given choices due to $\lambda<\frac{1}{2}$. We are left with checking the upper bound for $k_{j+1}$. In view of the mean value theorem, the monotonicity of $\Phi'$ and $\Phi(0)=0$ we find that
$\Phi(s)\leq s\Phi'(s)$ for all $s>0$, so that with $s=\Phi^{-1}(k_j)>0$ we get
\begin{align*}
k_j^2\theta_{k_j}^{-\frac{q_1-1}{q_1}}\big(\Phi^{-1}(k_j)\big)^{\frac{q_1-1}{q_1}}=k_j^2\big(\Phi'\big(\Phi^{-1}(k_j)\big)\big)^{\frac{q_1-1}{q_1}}\big(\Phi^{-1}(k_j)\big)^{\frac{q_1-1}{q_1}}\geq k_j^{2+\frac{q_1-1}{q_1}}.
\end{align*}
In particular, utilizing \eqref{eq:sec-alt-deld} and monotonicity of $\Phi^{-1}$ we have
\begin{align*}
\deld^2 k_j^{2+\frac{q_1-1}{q_1}}<\big(\Phi^{-1}(k_j)\big)^{\frac{q_1-1}{q_1}}\Big(\frac{\rd_j}{2}\Big)^{\dimN\kappa}\leq \big(\Phi^{-1}(M)\big)^{\frac{q_1-1}{q_1}}\rd_j^{\dimN\kappa}.
\end{align*}
Recalling that we took $C=\max\Big\{1,\Big(\big(\Phi^{-1}(M)\big)^{\frac{q_1-1}{q_1}}\deld^{-2}\Big)^\frac{1}{2+\frac{q_1-1}{q_1}}\Big\}\geq 1$ and $\beta=\frac{\dimN\kappa}{2+\frac{q_1-1}{q_1}}<1$, we therefore obtain
\begin{align*}
k_{j+1}=k_j\leq \Big(\big(\Phi^{-1}(M)\big)^{\frac{q_1-1}{q_1}}\deld^{-2}\rd_j^{\dimN\kappa}\Big)^\frac{1}{2+\frac{q_1-1}{q_1}}\leq C\rd_j^\beta.
\end{align*}\smallskip

\emph{Case 2b}: $\deld^2 k_j^2\theta_{k_j}^{-\frac{q_1-1}{q_1}}> \big(\tfrac{\rd_j}{2}\big)^{\dimN\kappa}$ and $\big|\big\{(x,t)\in Q_{\frac{\rd_j}{2}}^{\theta_{k_j}}:(w-\mu^-_{j}-k_j)_->0\big\}\big|\leq \sig_0\big|Q_{\frac{\rd_j}{2}}^{\theta_{k_j}}\big|$, i.e \eqref{eq:sec-alt-deld} is false and \eqref{eq:first-alt-1-ass1} is true for $r=\frac{\rd_j}{2}$, $k=k_j$ and $\theta=\theta_{k_j}$.\smallskip\\
Note that $\deld^2 k_j^2\theta_{k_j}^{-\frac{q_1-1}{q_1}}>\big(\frac{\rd_j}{2}\big)^{\dimN\kappa}$ also entails that \eqref{eq:first-alt-1-r} is false, since $\deld<\frac{1}{4}$. Accordingly, we conclude from Lemma~\ref{lem:first-alt-1} that
\begin{align*}
\Big|\Big\{(x,t)\in Q_{\frac{\rd_j}{4}}^{\theta_{k_j}}:\big(w-\mu^-_{j}-\tfrac{k_j}{2}\big)_->0\Big\}\Big|=0,
\end{align*}
implying that
\begin{align}\label{eq:it1-2b-inf}
{\textstyle\essinf_{Q_{\frac{\rd_j}{4}}^{\theta_{k_j}}}} w>\mu^-_{j}+\frac{k_j}{2}.
\end{align}
Similar to Case 1, we find from \eqref{eq:it1-c1}, monotonicity of $\Phi'$ and subadditivity of $\Phi^{-1}$ that
\begin{align*}
\Phi'\big(\Phi^{-1}(k_j)\big)\leq C_1\Phi'\big(\tfrac12\Phi^{-1}(k_j)\big)< C_1\Phi'\big(\tfrac12\Phi^{-1}(\tfrac74 k_j)\big)\leq C_1\Phi'\big(\Phi^{-1}(\tfrac78 k_j)\big).
\end{align*}
Thus, choosing $r_{j+1}:=\frac{\rd_j}{4\sqrt{C_1}}$, $k_{j+1}:=\frac{7}{8}k_j$ and $\theta_{j+1}:=\theta_{k_{j+1}}$ with $C_1>1$ as in \eqref{eq:it1-c1}, we obtain that
\begin{align*}
\rd_{j+1}^2\theta_{k_{j+1}}=\frac{\rd_j^2}{16 C_1\Phi'\big(\Phi^{-1}(\tfrac78 k_j)\big)}\leq\frac{\rd_j^2}{16\Phi'\big(\Phi^{-1}(k_j)\big)}\leq \Big(\frac{\rd_j}{4}\Big)^2\theta_{k_j}\quad\text{and}\quad \rd_{j+1}\leq\frac{\rd_j}{4},
\end{align*}
so that in fact $Q_{j+1}\subseteq Q_{\frac{\rd_j}{4}}^{\theta_{k_j}}\subseteq Q_j$. Hence, we may use \eqref{eq:it1-2b-inf} and $\essosc_{Q_j} w\leq 4 k_{j}$ to estimate
\begin{align*}
\essosc_{Q_{j+1}} w\leq \textstyle{\essosc_{Q_{\frac{\rd_j}{4}}^{\theta_{k_j}}}} w\leq \mu^+_{j}-\mu^-_{j}-\frac{k_j}{2}\leq 4k_j-\frac{k_j}{2}=\frac{7}{2}k_j=4 k_{j+1}.
\end{align*}
Since the properties for $k_{j+1}$ and $\rd_{j+1}$ are clearly satisfied because $\eta>\frac{31}{32}$ and $\lambda<\frac{1}{4\sqrt{C_1}}$, we are finished with this case.\smallskip

\emph{Case 2c}: $\deld^2 k_j^2\theta_{k_j}^{-\frac{q_1-1}{q_1}}> \big(\tfrac{\rd_j}{2}\big)^{\dimN\kappa}$ and $\big|\big\{(x,t)\in Q_{\frac{\rd_j}{2}}^{\theta_{k_j}}:(w-\mu^+_{j}+k_j)_+>0\big\}\big|\leq (1-\sig_0)\big|Q_{\frac{\rd_j}{2}}^{\theta_{k_j}}\big|$, i.e. \eqref{eq:sec-alt-deld} and \eqref{eq:first-alt-1-ass1} are both false. (Recall the reasoning at the start of Section~\ref{sec4:alt2}.)\smallskip\\
Here, writing $\gamma:=\frac{\sig_0}{2}<\frac14$, we can draw on Lemma~\ref{lem:second-alt} to obtain $$\Big|\Big\{(x,t)\in Q_{\frac{\rd_j}{4}}^{\gamma\theta_{k_j}}:\big(w-\mu^+_{j}+\tfrac{\deld}{2}k_j\big)_+>0\Big\}\Big|=0,$$
which also clearly entails
\begin{align}\label{eq:it1-2c-sup}
{\textstyle\esssup_{Q_{\frac{\rd_j}{4}}^{\gamma\theta_{k_j}}}} w<\mu^+_{j}-\frac{\deld}{2}k_j.
\end{align}
We choose $\rd_{j+1}:=\lambda\rd_j$, $k_{j+1}:=\eta k_j$ and $\theta_{j+1}:=\theta_{k_{j+1}}$. The properties for $k_j$ and $\rd_j$ are obviously satisfied. Additionally, since $2\eta>1$ we may follow the same reasoning as before to find that
\begin{align*}
\Phi'\big(\Phi^{-1}(k_j)\big)\leq C_1\Phi'\big(\tfrac{1}{2}\Phi^{-1}(k_j)\big)\leq C_1\Phi'\big(\tfrac{1}{2}\Phi^{-1}(2\eta k_j)\big)\leq C_1\Phi'\big(\Phi^{-1}(\eta k_j)\big).
\end{align*}
Therefore, in this case
\begin{align*}
\rd_{j+1}^2\theta_{k_{j+1}}=\frac{\gamma\rd_j^2}{16 C_1\Phi'\big(\Phi^{-1}(\eta k_j)\big)}\leq\frac{\gamma\rd_j^2}{16}\theta_{k_j}=\Big(\frac{\rd_j}{4}\Big)^2\gamma\theta_{k_j}\quad\text{and}\quad \rd_{j+1}\leq \frac{\rd_j}{4}
\end{align*}
and hence $Q_{j+1}\subseteq Q_{\frac{\rd_j}{4}}^{\gamma\theta_{k_j}}\subseteq Q_j$. This, together with \eqref{eq:it1-2c-sup} and $\essosc_{Q_j}w\leq 4k_j$, yields the estimate
\begin{align*}
\essosc_{Q_{j+1}} w\leq \textstyle{\essosc_{Q_{\frac{\rd_j}{4}}^{\gamma\theta_{k_j}}}} w\leq \mu^+_{j}-\frac{\deld}{2}k_j-\mu^-_j\leq 4k_j-\frac{\deld}{2}k_j=4\big(1-\frac{\deld}{8}\big)k_j=4k_{j+1},
\end{align*}
finishing the proof in this case and also the proof of the lemma as a whole.\hfill\qedhere
\end{bew}

The second iteration now starts with the assumption that $k_0<\mu_0^-$. Accordingly, in Lemma~\ref{lem:first-alt-1} and Lemma~\ref{lem:second-alt} we always choose $\theta=\theta_{\mu^-_j}$ during the iteration. In the proof, which mimics the reasoning above, we will mostly make sure that the switch in the scaling factors $\theta_j$ does not impede the nesting of the cylinders.

\begin{lemma}\label{lem:iteration-2}
Let $w$ be a nonnegative bounded weak solution of \eqref{eq:trans-gen-eq} with $M>0$ such that \eqref{eq:linfty-bound-w} holds. With $\lambda\in(0,\tfrac12)$, $\eta\in(\tfrac12,1)$, $\beta\in(0,1)$, $C=C(\Phi,M,q_1,q_2,a,\hat{b})\geq 1$ provided by Lemma~\ref{lem:iteration-1} the following holds: Given $x_0\in\Omega$, $0<t_0<T$, $0<k_0\leq\frac{M}{2}$ and $0<\rd_0<\min\big\{\frac{1}{\sqrt{\dimN}}\dist(x_0,\romega),\sqrt{\frac{t_0}{4\theta_{k_0}}},1\big\}$ set $\theta_0:={\theta_{k_0}}$ and $Q_0:=Q_{\rd_0}^{\theta_0}$. If $\essosc_{Q_{0}}w\leq 4k_0$ and if $\mu^-_{0}:=\essinf_{Q_0} w>k_0$ one can find non-increasing sequences $(\rd_j)_{j\in\N}$, $(k_j)_{j\in\N}$, $(\theta_j)_{j\in\N}$ such that with $Q_j:=Q_{\rd_j}^{\theta_j}$ the properties
\begin{equation}\label{eq:iteration-2}
Q_{j}\subseteq Q_{j-1},\ \;\essosc_{Q_{j}} w\leq 4k_j,\ \;\frac{1}{2}k_{j-1}\leq k_j\leq \max\big\{\eta k_{j-1}, C\rd_{j-1}^\beta\big\}\ \;\text{and}\ \;\lambda\rd_{j-1}\leq\rd_j\leq\frac{1}{2}\rd_{j-1}
\end{equation}
are satisfied for all $j\in\N$.
\end{lemma}

\begin{bew}
The proof is similar to the previous lemma. This time, using the fact that $(k_j)_{j\in\N}$ is non-increasing and that the $(Q_j)_{j\in\N}$ are nested, we always have $\mu^-_{j}:=\essinf_{Q_j} w> k_j$ inside the cylinder of the current iteration step. Accordingly, whenever we draw on Lemma~\ref{lem:first-alt-1} or Lemma~\ref{lem:second-alt} we are in the case $\theta=\theta_{\mu^-_{j}}$ and we will adjust our choice for the scaling factor of the next cylinder to $\theta_{j+1}=\theta_{\mu^-_{j}}$ instead of $\theta_{k_{j+1}}$ as taken in Lemma~\ref{lem:iteration-1}. Notice that by this choice, $\theta_{j+1}\leq\theta_j$ and $Q_{j+1}\subseteq Q_j$ are evident in most of the cases, as $(\rd_j)_{j\in\N}$ and $(\theta_j)_{j\in\N}$ are non-increasing, the latter due to $\mu^-_{j-1}\leq\mu^-_{j}$, so that $\rd_{j+1}^2\theta_{j+1}=\rd_{j+1}^2\theta_{\mu^-_{j}}\leq\rd_j^2\theta_{\mu^-_{j-1}}=\rd_j^2\theta_j$. Aside from the switch in scaling factors, the general steps are quite similar to the previous lemma, nevertheless, let us briefly verify the main steps.\medskip

Again we let the numbers $\sig_0=\sig_0(\Phi,M,q_1,q_2,a,\hat{b})\in(0,\tfrac12)$ and $\deld=\deld(\Phi,M,q_1,q_2,a,\hat{b})\in(0,\tfrac14)$ be given by Lemma~\ref{lem:first-alt-1} and Lemma~\ref{lem:second-alt}, respectively, and claim that the lemma holds if we take the same $\lambda=\lambda(\Phi,M,q_1,q_2,a,\hat{b})\in(0,\tfrac12)$, $\eta=\eta(\Phi,M,q_1,q_2,a,\hat{b})\in(\tfrac{31}{32},1)$, $\beta<1$ and $C=C(\Phi,M,q_1,q_2,a,\hat{b})\geq 1$ as in \eqref{eq:it1-const}. Now, assume \eqref{eq:iteration-2} holds for some $j\in\N$. For the choice for $\rd_{j+1}$, $k_{j+1}$ and $\theta_{j+1}$ we discern the following cases:\smallskip

\emph{Case 1}: $\essosc_{Q_j} w<2 k_j$.\smallskip\\
Pick $\rd_{j+1}:=\frac{\rd_j}{2}$, $k_{j+1}:=\frac{k_j}{2}$, $\theta_{j+1}:=\theta_{\mu^-_{j}}$. As discussed above $\theta_{j+1}\leq \theta_{j}$ and $Q_{j+1}\subseteq Q_j$ follows from $\rd_{j+1}\leq\rd_j$. Thus, $\essosc_{Q_{j+1}} w\leq\essosc_{Q_j} w<2k_j=4k_{j+1}$. The inequalities for $\rd_{j+1}$ and $k_{j+1}$ are obvious.\smallskip

\emph{Case 2}: $2k_j\leq \essosc_{Q_j} w\leq 4k_j$.\smallskip\\
Again we split this case in accordance with the Lemmas~\ref{lem:first-alt-1} and \ref{lem:second-alt}. The lemmas will be employed for $\vartheta=\theta_j$, $\rd=\frac{\rd_j}{2}$, $k=k_j$, which is possible due to $k_j<\mu^-_{j}$ and $\theta_{j+1}\leq\theta_{j}$ implying that here $\theta_j= \vartheta\geq \theta=\theta_{\mu^-_{j}}=\theta_{j+1}$.
\smallskip

\emph{Case 2a}: $\deld^2k_j^2\theta_{j+1}^{-\frac{q_1-1}{q_1}}\leq \big(\tfrac{\rd_j}{2}\big)^{\dimN\kappa}$.\smallskip\\
Set $\rd_{j+1}:=\frac{\rd_j}{2}$, $k_{j+1}:=k_j$, $\theta_{j+1}:=\theta_{\mu^-_j}$. The arguments for $\theta_{j+1}\leq\theta_j$, $Q_{j+1}\subseteq Q_j$ and $\essosc_{Q_{j+1}} w$ remain the same as in the previous case. That the properties for $\rd_{j+1}$ and the lower bound for $k_j$ are satisfied is also obvious. Let us check the upper bound for $k_{j+1}$. Since $k_j\leq k_{j-1}\leq\dots \leq k_0<\mu^-_{0}\leq \dots\leq\mu^-_{j}$ we have $\theta_{k_j}\geq\theta_{k_0}=\theta_{0}\geq \theta_{j+1}$ and therefore the case assumption also implies that $\deld^2 k_j^2\theta_{k_j}^{-\frac{q_1-1}{q_1}}\leq\big(\frac{\rd_j}{2} \big)^{\dimN\kappa}$ and we can copy the steps presented in Case~2a of Lemma~\ref{lem:iteration-1}.\smallskip

\emph{Case 2b}: $\deld^2k_j^2\theta_{j+1}^{-\frac{q_1-1}{q_1}}> \big(\tfrac{\rd_j}{2}\big)^{\dimN\kappa}$ and $\big|\big\{(x,t)\in Q_{\frac{\rd_j}{2}}^{\theta_{j+1}}:(w-\mu^-_{j}-k_j)_->0\big\}\big|\leq\sig_0\big|Q_{\frac{\rd_j}{2}}^{\theta_{j+1}}\big|$.\smallskip\\
We choose $\rd_{j+1}:=\frac{\rd_j}{4}$, $k_{j+1}:=\frac{7}{8}k_j$ and $\theta_{j+1}:=\theta_{\mu^-_{j}}$. From $\deld<\frac{1}{4}$ we find that the assumption in this case implies $k_j\theta_{j+1}^{-\frac{q_1-1}{q_1}}>\big(\frac{\rd_j}{2}\big)^{\dimN\kappa}$ and in view of Lemma~\ref{lem:first-alt-1} we therefore obtain
\begin{align*}
\Big|\Big\{(x,t)\in Q_{\frac{\rd_j}{4}}^{\theta_{j+1}}:\big(w-\mu^-_{j}-\tfrac{k_j}{2}\big)_->0\Big\}\Big|=0\quad\text{and thus }\quad{\textstyle\essinf_{Q_{\frac{\rd_j}{4}}^{\theta_{j+1}}}} w>\mu^-_{j}+\frac{k_j}{2},
\end{align*}
so that we may conclude analogously to Case 2b of Lemma~\ref{lem:iteration-1} if we show $Q_{j+1}\subseteq Q_{\frac{\rd_j}{4}}^{\theta_{j+1}}\subseteq Q_{j}$. The first inclusion is clear from the choice of $\rd_{j+1}$ since the cylinders have the same scaling factor, whereas the second inclusion relies on $\theta_{j+1}\leq\theta_{j}$.\smallskip

\emph{Case 2c}: $\deld^2k_j^2\theta_{j+1}^{-\frac{q_1-1}{q_1}}> \big(\tfrac{\rd_j}{2}\big)^{\dimN\kappa}$ and $\big|\big\{(x,t)\in Q_{\frac{\rd_j}{2}}^{\theta_{j+1}}:(w-\mu^+_{j}+k_j)_+>0\big\}\big|\leq(1-\sig_0)\big|Q_{\frac{\rd_j}{2}}^{\theta_{j+1}}\big|$.\smallskip\\
Pick $\rd_{j+1}:=\lambda\rd_j$, $k_{j+1}:=\eta k_j$, $\theta_{j+1}:=\theta_{\mu^-_{j}}$ and let $\gamma:=\frac{\varpi_0}{2}\in(0,\tfrac14)$. In view of Lemma~\ref{lem:second-alt} we have
\begin{align*}
\Big|\Big\{(x,t)\in Q_{\frac{\rd_j}{4}}^{\gamma\theta_{j+1}}:\big(w-\mu^+_{j}+\tfrac{\deld}{2}k_j\big)_+>0\Big\}\Big|=0\quad\text{and thus}\quad {\textstyle\esssup_{Q_{\frac{\rd_j}{4}}^{\gamma\theta_{j+1}}}} w<\mu^+_{j}-\frac{\deld}{2}k_j.
\end{align*}
Since $\lambda\leq \frac{\sqrt{\gamma}}{4}<\frac18$ we have
\begin{align*}
\rd_{j+1}^2\theta_{j+1}\leq\frac{\gamma\rd_j^2}{16}\theta_{j+1}=\Big(\frac{\rd_j}{4}\Big)^2\gamma\theta_{j+1}<\rd_j^2\theta_j,
\end{align*}
which together with $\rd_{j+1}\leq \frac{\rd_j}{4}\leq \rd_j$ implies $Q_{j+1}\subseteq Q_{\frac{\rd_j}{4}}^{\gamma\theta_{j+1}}\subseteq Q_j$. Thus,
\begin{align*}
\essosc_{Q_{j+1}}w\leq \textstyle{\essosc_{Q_{\frac{\rd_j}{4}}^{\gamma\theta_{j+1}}}}w\leq \mu^+_{j}-\tfrac{\deld}{2}k_j-\mu_j^-\leq 4k_j-\tfrac{\deld}{2}k_j=4\eta k_j=4k_{j+1},
\end{align*}
as intended, finishing the proof.\hfill\qedhere
\end{bew}

An immediate consequence of a combination of the two iterations above is the following decay estimate for the essential oscillation with respect to shrinking cylinders.

\begin{corollary}\label{cor:iteration}
Let $w$ be a nonnegative bounded weak solution of \eqref{eq:trans-gen-eq} with $M>0$ such that \eqref{eq:linfty-bound-w} holds and set $k_0:=\frac{M}{2}$. Then, there are constants $\lambda_\star=\lambda_\star(\Phi,M,q_1,q_2,a,\hat{b})\in(0,\tfrac12)$, $\eta_\star=\eta_\star(\Phi,M,q_1,q_2,a,\hat{b})\in(\tfrac12,1)$ and $C_\star=C_\star(\Phi,M,q_1,q_2,a,\hat{b})\geq 1$ such that for all $x_0\in\Omega$, $0<t_0<T$ and $0<\rd_0<\min\big\{\frac{1}{\sqrt{\dimN}}\dist(x_0,\romega),\sqrt{\frac{t_0}{4\theta_{k_0}}},1\big\}$ the estimate
\begin{align*}
\essosc_{Q_{\lambda_\star^j \rd_0}^{\theta_M}} w\leq C_\star\eta_\star^j\quad\text{for all }j\in\N
\end{align*}
holds.
\end{corollary}

\begin{bew} Let $Q_0:=Q_{\rd_0}^{\theta_{k_0}}$. Then $\essosc_{Q_0} w\leq 2M=4k_0$ and accordingly, $\rd_0$, $Q_0$ and $k_0$ fulfill the requirements of Lemma~\ref{lem:iteration-1}, providing $\lambda=\lambda(\Phi,M,q_1,q_2,a,\hat{b})\in(0,\tfrac12)$, $\eta=\eta(\Phi,M,q_1,q_2,a,\hat{b})\in(\tfrac12,1)$, $\beta=\beta(q_1,q_2)\in(0,1)$, $C=C(\Phi,M,q_1,q_2,a,\hat{b})\geq 1$ and sequences $(r_j)_{j\in\N}$, $(k_j)_{j\in\N}$, $(\theta_j)_{j\in\N}$ such that with $Q_j:=Q_{\rd_j}^{\theta_j}$ either \eqref{eq:iteration-1} holds for all $j\in\N$, or there is $j_\star\in\N$ such that \eqref{eq:iteration-1} holds for $j\leq j_\star$ and $\mu^-_{j_\star}:=\essinf_{Q_{j_\star}}w>k_{j_\star}$. If such a $j_\star$ exists, then $\rd_{j_\star}$, $Q_{j_\star}=Q_{\rd_{\js}}^{\theta_{\js}}$ and $k_{j_\star}$ satisfy the conditions of Lemma~\ref{lem:iteration-2} and we find that \eqref{eq:iteration-2} also holds for all $j>j_\star$. Accordingly, we may continue the sequence obtained for $j\leq j_\star$ and find that, irrespective of the existence of $j_\star$,
\begin{align*}
\lambda \rd_{j-1}\leq\rd_j\leq\frac{1}{2}\rd_{j-1},\quad k_j\leq \max\big\{\eta k_{j-1}, C\rd_{j-1}^{\beta}\big\},\quad\essosc_{Q_j} w\leq 4k_j
\end{align*}
hold for all $j\in\N$. From these recursive inequalities we conclude that
\begin{align*}
\lambda_\star^j\rd_0\leq \rd_j\leq 2^{-j}\rd_0,\qquad\text{and}\qquad\essosc_{Q_j} w\leq 4k_j\leq 4\max\big\{k_0\eta^j, C 2^{-(j-1)\beta}\rd_0^{\beta}\big\}\leq C_\star\eta_\star^j
\end{align*}
hold for all $j\in\N$, where we used $r_0\leq 1$ and $k_0\leq M$ and put $\lambda_\star=\lambda$, $\beta_\star=\beta$, $\eta_\star:=\max\{\eta,2^{-\beta}\}$ and $C_\star:=4\max\{M,2^\beta C\}$. Then, since $\theta_j\geq \theta_M$ and $\lambda_\star^j\rd_0\leq\rd_j$ for all $j\in\N$ we have $$Q_{\lambda_\star^j\rd_0}^{\theta_M}\subseteq Q_{\lambda_\star^j \rd_0}^{\theta_j}\subseteq  Q_{\rd_j}^{\theta_j}=Q_j\quad\text{for all }j\in\N,$$
and therefore \[\essosc_{Q_{\lambda_\star^j\rd_0}^{\theta_M}} w\leq \essosc_{Q_j} w\leq C_\star\eta_\star^j\quad\text{for all }j\in\N.\hfill\qedhere\]
\end{bew}

\setcounter{equation}{0}
\section{Hölder regularity in compact subsets \texorpdfstring{of $\OmT$}{}}\label{sec5}
Note that in light of Corollary~\ref{cor:iteration} the existence of $\esslim_{(x,t)\to(x_0,t_0)} w(x,t)$ for all $(x_0,t_0)\in\OmT$ is evident and, moreover, setting $\widetilde{w}(x_0,t_0):=\esslim_{(x,t)\to(x_0,t_0)} w(x,t)$ specifies a continuous representative of our equivalence class of weak solutions. For the remainder will identify $w$ with this continuous representative $\widetilde w$. We proceed by verifying that $w$ satisfies a local Hölder condition.
\begin{lemma}\label{lem:höl}
Let $w$ be a nonnegative bounded weak solution of \eqref{eq:trans-gen-eq} with $M>0$ such that \eqref{eq:linfty-bound-w} holds. Set $\Gamma:=\partial_p(\OmT)=\romega\times[0,T]\,\cup\,\Omega\times\{0\}$. Then there are constants $\alpha=\alpha(\Phi,M,q_1,q_2,a,\hat{b})\in(0,1)$ and $C=C(\Phi,M,q_1,q_2,a,\hat{b})>0$ such that for any compact $K\subset \OmT$
\begin{align*}
\big|w(x_0,t_0)-w(x_1,t_1)\big|\leq C\left(\frac{|x_0-x_1|+|t_0-t_1|^\frac{1}{2}}{d_p(K,\Gamma\big)}\right)^\alpha
\end{align*}
for all $(x_0,t_0)$, $(x_1,t_1)\in K$, where $d_p(K,\Gamma):=\inf_{(x,t)\in K}\min\big\{\frac{1}{\sqrt{\dimN}}\dist(x,\romega),\sqrt{t},1\big\}.$
\end{lemma}

\begin{bew}
Given any compact $K\subset\OmT$ we set $d_K:=d_p(K,\Gamma)$ and pick $R:=\frac{d_K}{2(1+\sqrt{\theta_{0}})}$ with $\theta_0:=\theta_{\frac{M}{2}}$. Given two points $(x_0,t_0),(x_1,t_1)\in K$, where without loss of generality we assume $t_0>t_1$, we fix the cylinder $Q_R^{\theta_M}:=Q_R(x_0,t_0-\theta_M R^2,t_0)$. Since $d_K\leq \min\big\{\frac{1}{\sqrt{\dimN}}\dist(x_0,\romega),\sqrt{t_0},1\big\}$ we in particular have $R\leq \min\big\{\frac{1}{\sqrt{\dimN}}\dist(x_0,\romega),\sqrt{\frac{t_0}{4\theta_0}},1\big\}$ and because of $\theta_{0}>\theta_M$ also $Q_R^{\theta_M}\subset Q_R^{\theta_0}\subset\OmT$. Denote by $\lambda_\star=\lambda_\star(\Phi,M,q_1,q_2,a,\hat{b})\in(0,\tfrac12)$, $\eta_\star=\eta_\star(\Phi,M,q_1,q_2,a,\hat{b})\in(\tfrac12,1)$ and $C_\star=C_\star(\Phi,M,q_1,q_2,a,\hat{b})\geq1$ the numbers from Corollary \ref{cor:iteration} and set $\alpha=\alpha(\Phi,M,q_1,q_2,a,\hat{b})=\frac{\log\eta_\star}{\log\lambda_\star}\in(0,1)$. Now consider the following cases:

If $(x_1,t_1)\not\in Q_{\lambda_\star R}^{\theta_M}$ then $\|x_0-x_1\|_\infty\geq \lambda_\star R$ or $\theta_M^{-\frac12}|t_0-t_1|^\frac{1}{2}\geq \lambda_\star R$. Accordingly, because of $\|x\|_\infty\leq|x|$, also $|x_0-x_1|+\theta_M^{-\frac12}|t_0-t_1|^\frac{1}{2}\geq \lambda_\star R$ and 
since $|w|\leq M$, we can thus estimate
\begin{align*}
\big|w(x_0,t_0)-w(x_1,t_1)\big|\leq 2M&\leq 2M\left(\frac{|x_0-x_1|+\theta_M^{-\frac12}|t_0-t_1|^\frac12}{\lambda_\star R}\right)^\alpha\\
&\leq \frac{2^{1+\alpha}M}{\lambda_\star^\alpha}\big(1+\sqrt{\theta_{0}}\big)^\alpha\big(1+\theta_M^{-\frac12}\big)^\alpha\left(\frac{|x_0-x_1|+|t_0-t_1|^\frac12}{d_K}\right)^\alpha\\
&\leq \frac{4M}{\lambda_\star}\big(1+\sqrt{\theta_0}\big)\big(1+\theta_M^{-\frac12}\big)\left(\frac{|x_0-x_1|+|t_0-t_1|^\frac12}{d_K}\right)^\alpha,
\end{align*}
where we made use of $\lambda_\star\leq 1$ and $\alpha\in(0,1)$.\smallskip

On the other hand, if $(x_1,t_1)\in Q_{\lambda_\star R}^{\theta_M}$, we can pick $j\geq 1$ such that $(x_1,t_1)\in Q_{\lambda_\star^j R}^{\theta_M}$ but $(x_1,t_1)\not\in Q_{\lambda_\star^{j+1} R}^{\theta_M}$, meaning that $\lambda_\star^{j+1}R\leq |x_0-x_1|+\theta_M^{-\frac12}|t_0-t_1|^\frac12$. Drawing on Corollary~\ref{cor:iteration} and our choice for $\alpha$ then provides the estimate 
\begin{align*}
\big|w(x_0,t_0)-w(x_1,t_1)\big|&\leq \essosc_{Q_{\lambda_\star^j R}^{\theta_M}} w\leq C_\star\eta_\star ^j=C_\star\lambda_\star^{j\alpha}=\frac{C_\star \big(\lambda_\star^{j+1}R\big)^\alpha}{\big(\lambda_\star R\big)^{\alpha}}\\
&\leq \frac{C_\star}{\lambda_\star^\alpha}\left(\frac{|x_0-x_1|+\theta_M^{-\frac12}|t_0-t_1|^\frac12}{R}\right)^\alpha\\
&\leq 2^\alpha\frac{C_\star}{\lambda_\star^\alpha}\big(1+\sqrt{\theta_0}\big)^\alpha\big(1+\theta_M^{-\frac12}\big)^\alpha\left(\frac{|x_0-x_1|+|t_0-t_1|^\frac12}{d_K}\right)^\alpha\\
&\leq\frac{2C_\star}{\lambda_\star}\big(1+\sqrt{\theta_0}\big)\big(1+\theta_M^{-\frac12}\big)\left(\frac{|x_0-x_1|+|t_0-t_1|^\frac12}{d_K}\right)^\alpha
\end{align*}
and we may conclude the lemma.
\end{bew}

\begin{proof}[\textbf{Proof of Theorem \ref{theo:1}}:]
The proof of Theorem~\ref{theo:1} is an immediate consequence of Lemma~\ref{lem:höl}.
\end{proof}

\setcounter{equation}{0}
\section{Regularity up to the boundary. Proof of Theorem~\ref{theo:2}}\label{sec6:boundary}
In the remainder of this paper we will present the necessary changes to treat the Hölder regularity of the associated initial-boundary-value problem to \eqref{eq:trans-gen-eq} as stated in \eqref{eq:bvp-w}. Let us briefly recall that in this setting we suppose that $\romega$ is of class $C^{1}$, that the initial data $w_0$ are Hölder continuous and that the boundary data $g\in \CSp{0}{\romega\times[0,T]}$ are assumed to admit an extension onto $\Omega$ for a.e. $t\in(0,T)$ denoted by $\hat{g}$, such that the conditions in \eqref{eq:cond-g} are satisfied, i.e.
\begin{align*}
\big|\hat{g}\big(x,t,u(x,t)\big)\big|&\leq g_0(x,t)\Phi\big(u(x,t)\big)\quad\text{a.e. in } \OmT,\nonumber\\
\big|\partial_u \hat{g}\big(x,t,u(x,t)\big)\big|&\leq g_0(x,t)\Phi'\big(u(x,t)\big)\quad\text{a.e. in } \OmT,\\
\text{and}\quad \big|\partial_{x_i} \hat{g}\big(x,t,u(x,t)\big)\big|&\leq g_0^2(x,t),\ i\in\{1,\dots,\dimN\}\ \text{a.e. in } \OmT,\nonumber
\end{align*}
hold with nonnegative $g_0\in\LSp{2q_1}{(0,T);\Lo[2q_2]}$ where $q_1,q_2>1$ satisfying \eqref{eq:cond-kappa} for some $\kappa\in(0,\frac{2}{\dimN})$.\smallskip

We will now assume that $w$ is a nonnegative bounded weak solution of \eqref{eq:bvp-w} on $\OmT$ for some $T>0$ in the sense of Definition~\ref{def:weaksol-bv} such that $w$ satisfies \eqref{eq:linfty-bound-w-2} for some $M>0$. The following line of reasoning is not too different from our previous arguments and adopts well-established methods (see e.g. \cite[Chapter II.4]{dibenedettoDegenerateParabolicEquations1993}). In particular, we will derive local energy estimates akin to the estimates presented in Lemma~\ref{lem:local-energy} and Lemma~\ref{lem:log-est}. This time, however, the cylinders will be spatially centered around a boundary point and, accordingly, the cutoff functions involved do not disappear on $\romega$, so that we will have to discuss how to treat the newly appearing boundary integral. We will see that the structure of the inequality can be mostly retained with only minor changes necessary. Since the arguments in Section~\ref{sec4} and Section~\ref{sec5} were solely dependent on these energy estimates -- and not in the slightest on the PDE in question -- we can then, for the most part, duplicate the arguments of the previous sections to arrive at the conclusion.

\subsection{Local energy estimates and regularity at the lateral boundary}
To shorten some of the upcoming notation we introduce for $x=(x_1,\dots,x_\dimN)\in\R^\dimN$ the dimension reduced vector $\cx$ via the notation $\cx:=(x_1,\dots,x_{\dimN-1})\in\R^{\dimN-1}$ and write $$\check{\mK}_\rd(x_0)=\big\{\cx\in\R^{\dimN-1}\ \big\vert\ \|\cx-\cx_0\|_{\infty}<\rd\big\}.$$
Since we assume $\romega$ to be of class $C^1$ we know (\cite{fiorenza2016}) that for fixed $x_0\in\romega$ there is a reorientation and relabeling of the coordinate system such that in a neighborhood around $x_0$ we can express $\bomega$ as the epigraph of $\gamma\in \CSp{1}{\R^{\dimN-1}}$. In fact, the proof of \cite[Theorem 3.2.1]{fiorenza2016} tells us that there is an orthogonal matrix $S=S(x_0)\in\R^{\dimN\times\dimN}$, which links the canonical coordinate vector $(x_1,\dots,x_\dimN)$ of a point $x$ to its new coordinates $(y_1,\dots,y_\dimN)$ by the relations 
\begin{align}\label{eq:coord-trafo}
x=x_0+Sy\quad\text{ and }\quad y=S^T(x-x_0).
\end{align}
The result \cite[Theorem 3.2.1]{fiorenza2016} further substantiates that in the new coordinate system the boundary point $x_0$ corresponds to the origin and $\gamma$ and all its first derivatives vanish at $\check{0}$, so that the hyperplane $\{y_\dimN=0\}$ is in fact tangent to the boundary at $x_0$ and the outward unit normal at $x_0$ points straight down. With respect to the $y$-coordinates $\mK_\rd(x_0)\cap\Omega$ can be thought of as a cube with just a deformed bottom. The new coordinate system is quite beneficial for the boundary regularity investigations since the orthogonal nature of $S$ makes the switch between the coordinate systems highly compatible with the PDE as a whole. In fact, drawing on \eqref{eq:coord-trafo}, one can verify from direct calculations that the structure of the PDE does not change. Let us consider the following: Given $w(x,t)$ we set $\tilde{w}(y,t):=w(x_0+Sy,t)$. From the fact that $y=S^T(x-x_0)$ entails $\frac{\partial y_j}{\partial x_i}=S_{ij}$ and the chain rule we find that $$\frac{\partial w(x,t)}{\partial x_i}=\sum_{j=1}^\dimN \frac{\partial \tilde{w}(y,t)}{\partial y_j}\frac{\partial y_j}{\partial x_i}=\sum_{j=1}^\dimN S_{ij}\frac{\partial \tilde{w}(y,t)}{\partial y_j},\quad i\in\{1,\dots,\dimN\},$$ so that in particular $\nabla_x w(x,t)=S\nabla_y\tilde{w}(y,t)$. Moreover, differentiating once more
$$\frac{\partial^2 w(x,t)}{\partial x_i^2}
=\sum_{j,k=1}^\dimN S_{ij}S_{ik}\frac{\partial^2\tilde{w}(y,t)}{\partial y_k\partial y_j},\quad i\in\{1,\dots,\dimN\}.$$
Recalling that $S$ is orthogonal, $$\sum_{i=1}^\dimN S_{ij}S_{ik}=\begin{cases}1\quad&\text{if }j=k,\\0&\text{if }j\neq k,\end{cases}$$ and thus
$$\Delta_x w(x,t)=\sum_{i=1}^\dimN\sum_{j,k=1}^\dimN S_{ij}S_{ik}\frac{\partial^2\tilde{w}(y,t)}{\partial y_k\partial y_j}=\sum_{j,k=1}^\dimN\frac{\partial^2\tilde{w}(y,t)}{\partial y_k \partial y_j}\sum_{i=1}^\dimN S_{ij}S_{ik}=\sum_{j=1}^\dimN\frac{\partial^2 \tilde{w}(y,t)}{\partial y_j^2}=\Delta_y \tilde{w}(y,t).$$
Similarly, setting $\tilde{a}(y,t):=S^T a(x_0+Sy,t)$ and using the chain rule one calculates that
\begin{align*}
\nabla_x\cdot\big(a(x,t)\Phi^{-1}(w(x,t))\big)&=\sum_{i=1}^\dimN \frac{\partial}{\partial x_i} \Big(a_i(x,t)\Phi^{-1}\big(w(x,t)\big)\Big)\\
&=\sum_{i,j=1}^\dimN \frac{\partial}{\partial y_j}\Big(\sum_{k=1}^\dimN S_{ik}\tilde{a}_k(y,t)\Phi^{-1}\big(\tilde{w}(y,t)\big)\Big)\frac{\partial y_j}{\partial x_i}\\
&=\sum_{j=1}^\dimN \frac{\partial}{\partial y_j}\Big(\sum_{i,k=1}^\dimN S_{ij}S_{ik}\tilde{a}_k(y,t)\Phi^{-1}\big(\tilde{w}(y,t)\big)\Big)\\
&=\sum_{j=1}^\dimN\frac{\partial}{\partial y_j}\Big(\tilde{a}_j(y,t)\Phi^{-1}\big(\tilde{w}(y,t)\big)\Big)=\nabla_y\cdot \big(\tilde{a}(y,t)\Phi^{-1}(\tilde{w}(y,t))\big).
\end{align*}
The coordinate transform also entails that the outward normals are linked by $\tilde{\nu}(y)=S^T\nu(x_0+Sy)$ and, using that $S$ preserves the inner product, we have
\begin{align*}
\big(\nabla_x w(x,t)+a(x,t)\Phi^{-1}(w(x,t))\big)\cdot\nu(x)&=\big(S\nabla_y \tilde{w}(y,t)+S\tilde{a}(y,t)\Phi^{-1}(\tilde{w}(y,t))\big)\cdot S\tilde{\nu}(y)\\&=\big(\nabla_y \tilde{w}(y,t)+\tilde{a}(y,t)\Phi^{-1}(\tilde{w}(y,t))\big)\cdot\tilde{\nu}(y).
\end{align*}
With $\tilde{g}(y,t,\xi)=g(x_0+Sy,t,\xi)$ and $\tilde{b}(y,t,\xi)=b(x_0+Sy,t,\xi)$ we conclude from the above that \eqref{eq:bvp-w} can be rewritten in the spatial $y$-coordinates without any changes in the structure. Moreover, expanding the arguments above to the bounds in \eqref{eq:cond-g}, we find that the bounds with respect to the original $x$-coordinates  transfer almost identically to the extension of $\tilde{g}$ given by $\tilde{\hat{g}}(y,t,\xi)=\hat{g}(x_0+Sy,t,\xi)$ by setting $\tilde{g}_0(y,t)=\sqrt{N}g_0^2(x_0+Sy,t)$. The additional factor of $\sqrt{N}$ therein is necessary for the bound on the spatial derivatives and is an upper bound on the column sum operator norm of $S$. Let us conclude the above remarks by summarizing that without loss of generality we may assume that for a fixed $x_0\in\romega$ the problem \eqref{eq:bvp-w} is already written in the coordinate system described in \cite[Theorem 3.2.1]{fiorenza2016}. In the lemma below, we will just further restrict the size of the neighborhood of $x_0$ in order to get a stricter control on the gradient of $\gamma$ which we require to extend Lemma~\ref{lem:deformed-DeGiorgi} to contorted cylinders (see Lemma~\ref{lem:deformed-DeGiorgi} below) in a straightforward manner.

\begin{lemma}\label{lem:coord-transf}
Suppose that $\romega$ is of class $C^{1}$. Let $x_0\in\romega$ and assume that the coordinate axes have been reoriented and relabeled such that $x_0$ is an element of the hyperplane $\{x_\dimN=0\}$ and that the outward unit normal $\nu(x_0)$ satisfies $\nu(x_0)=(0,\dots,0,-1)^T$. Then, there exist $\gamma\in\CSp{1}{\R^{\dimN-1}}$ and $\rd_0\in(0,1]$ such that \begin{align*}
\big|\nabla\gamma(\check{x})\big|\leq\frac{1}{2\sqrt{\dimN}}\quad\text{for all }\check{x}\in\check{\mK}_{\rd_0}(x_0).
\end{align*} and such that for all $\rd\in(0,\rd_0]$ 
\begin{align}\label{eq:dom-traf}
\begin{split}
\mK_\rd(x_0)\cap\romega&=\big\{x\in\R^\dimN\ \vert\  \check{x}\in\check{\mK}_\rd(x_0),\ x_\dimN=\gamma(\check{x})\big\}\\
\text{and}\qquad \mK_\rd(x_0)\cap\Omega&=\big\{x\in\R^\dimN\ \vert\ \check{x}\in\check{\mK}_\rd(x_0),\ \gamma(\check{x})< x_\dimN < r\big\}.
\end{split}
\end{align}
\end{lemma}
\begin{bew}
We note that the assumed coordinate system coincides with the system described in \cite[Theorem 3.2.1]{fiorenza2016} and that accordingly there are $\rho_0>0$ and $\gamma\in\CSp{1}{\R^{\dimN-1}}$ such that
$$B_{\rho_0}(x_0)\cap \Omega=\big\{x\in B_{\rho_0}(x_0)\,\vert\, x_\dimN>\gamma(x_1,\dots,x_{\dimN-1})\big\}.$$ Recall moreover that the partial derivatives of $\gamma$ vanish when evaluated at $\cx_0$. Accordingly, by continuity we can find $\rho_1\in(0,\rho_0)$ with the property that for all $x\in B_{\rho_1}(x_0)\cap\Omega$ we have $|\nabla\gamma(\cx)|\leq \frac{1}{2\sqrt{\dimN}}$. In light of the mean value theorem, we hence find that for all $0<\rd\leq \rd_0:=\min\{\frac{\rho_1}{\sqrt{\dimN}},1\}$ and each $\cx\in\check{\mK}_\rd(x_0)$ we have $|\gamma(\cx)|<\rd$, ensuring that $\romega$ does not leave $\mK_\rd(x_0)$ through the top or bottom of the cube with respect to the $x_\dimN$-direction. Consequently, for $\rd\in(0,\rd_0]$ we obtain \eqref{eq:dom-traf}.
\end{bew}

Similar to before we denote for $K\in\R$ the sublevel and superlevel sets of $w$ in $\mK_\rd(x_0)\cap\,\Omega$ at time $\tau$ by
\begin{align*}
B_{K,\rd,x_0}^{\pm}(\tau)&:=\Big\{ x\in \mK_\rd(x_0)\cap\,\Omega\ \big\vert\ \big(w(x,\tau)-K\big)_\pm>0\Big\},
\end{align*}
again dropping $x_0$ from the subscript if the center of the box is clear. (We also refer the reader to \cite[Proposition II.4.1]{dibenedettoDegenerateParabolicEquations1993} for related arguments.) 

\begin{lemma}\label{lem:loc-en-bdry}
Let $M>0$. Suppose that $\romega$ is of class $C^{1}$. Assume that $\hat{g}$ satisfies \eqref{eq:cond-g}. Let $w$ be a nonnegative bounded weak solution of \eqref{eq:bvp-w} satisfying \eqref{eq:linfty-bound-w-2}. There is $C=C(\Phi,M,q_1,q_2,a,\hat{b},g_0)>0$ with the following property: For all $x_0\in\romega$, all $0< T_0<T_1\leq T$ there is $\rd_0\in(0,1]$ such that for all $0<\rd\leq r_0$ and each $K\in(0,2M]$ the functions $\Lambda^\pm$ as defined in \eqref{eq:Lambda} and every smooth cutoff function $\psi$ inside the cylinders $Q_\rd=Q_\rd(x_0,T_0,T_1)$ fulfilling $0\leq \psi\leq 1$ and $\psi_{\vert\partial_p Q_\rd}=0$ satisfy
\begin{align}\label{eq:loc-en-bdry}
&\int\limits_{\mK_\rd\cap\,\Omega}\!\!\!\Lambda^\pm\big(w(\cdot,T_1)\big)\psi^2(\cdot,T_1)+\frac{1}{4}\intQO \big|\nabla\big((w-K)_\pm\psi\big)\big|^2\nonumber\\
\leq\ &\intQO\!\! \Lambda^\pm(w)\big(\psi^2\big)_t+2\intQO(w-K)_\pm^2|\nabla\psi|^2+C\Big(\int_{T_0}^{T_1}|B_{K,\rd}^\pm(t)|^{\frac{q_1(q_2-1)}{q_2(q_1-1)}}\Big)^\frac{q_1-1}{q_1}.
\end{align}
\end{lemma}

\begin{bew}
As explained at the start of the section we may without loss of generality assume that \eqref{eq:bvp-w} is written in the coordinate system fulfilling the conditions of Lemma~\ref{lem:coord-transf} and let $r_0\in(0,1]$ and $\gamma\in\CSp{1}{\R^{\dimN-1}}$ be provided by Lemma~\Ref{lem:coord-transf}. We fix an arbitrary $\rd\in(0,\rd_0]$ and pick any smooth cutoff function $\psi$ in $Q_\rd$ satisfying $0\leq\psi\leq 1$ and $\psi_{\vert\partial_p Q_\rd}=0$ and set $\Psi_h=\pm\big(\Phi\big([\Phi^{-1}(w)]_h\big)-K\big)_{\pm}\psi^2$. Noticing that $\Psi_h$ is an eligible test-function for \eqref{eq:weak-sol-bv-steklov}, we repeat the steps of Lemma~\ref{lem:local-energy} and obtain
\begin{align*}
0&=\!\intQO\!\!\big(\Lambda^\pm(w)\big)_t\psi^2\pm\!\intQO\!\!\big(\nabla w+a(x,t)\Phi^{-1}(w)\big)\cdot\nabla\big((w-K)_\pm\psi^2\big)\mp\!\intQO\!\! b\big(x,t,\Phi^{-1}(w)\big)(w-K)_\pm\psi^2\\&\hspace*{3.cm}\mp\int_{T_0}^{T_1}\!\int_{\mK_\rd\cap\romega} g\big(x,t,\Phi^{-1}(w)\big)(w-K)_\pm\psi^2=:I_1^\pm+I_2^\pm+I_3^\pm+I_4^\pm
\end{align*}
due to $\psi$ disappearing on $\partial_p Q_\rd$ but not on $\romega$. The first three integrals may be treated in exactly the same fashion as in Lemma~\ref{lem:local-energy} -- while replacing $A_{K,\rd}^\pm$ with $B_{K,\rd}^\pm$ -- using the well-signed term $$\intQO\!\big|\nabla\big((w-K)_\pm\psi\big)\big|^2,$$ so we only concern ourselves with a detailed treatment of $I_4^\pm$. In light of \eqref{eq:dom-traf} we can draw on $\gamma$ to rewrite the boundary integral in the following way.
\begin{align*}
|I_4^\pm|&=\Big|\int_{T_0}^{T_1}\!\!\int_{\mK_\rd\cap\romega} g\big(x,\tau,\Phi^{-1}(w)\big)(w-K)_\pm\psi^2\intd S\intd \tau\Big|\\
&=\Big|\int_{T_0}^{T_1}\!\!\int_{\check{\mK}_\rd} \Big(g\big(x,\tau,\Phi^{-1}(w)\big)\big(w-K)_\pm\psi^2\Big)\Big\vert_{x_\dimN=\gamma(\cx)}\sqrt{1+|\nabla\gamma(\cx)|^2}\intd \cx\intd \tau\Big|
\end{align*}
Making use of the assumed extension $\hat{g}$ of $g$ onto $\Omega$ for a.e. $t\in(0,T)$, the facts that $\psi=0$ on $\partial_p Q_\rd$ and $|\nabla\gamma|\leq 1$ on $\check{\mK}_\rd$ and \eqref{eq:dom-traf} once more, we can further estimate
\begin{align*}
|I_4^\pm|&\leq\sqrt{2}\int_{T_0}^{T_1}\!\!\int_{\check{\mK}_\rd}\Big|\int_{\gamma(\cx)}^\rd \frac{\partial}{\partial x_\dimN}\Big(\hat{g}\big(x,\tau,\Phi^{-1}(w)\big)(w-K)_\pm\psi^2\Big)\intd x_\dimN\Big|\intd\cx\intd\tau\\
&\leq \sqrt{2}\int_{T_0}^{T_1}\!\!\int_{\mK_\rd\cap\,\Omega}\Big|\frac{\partial}{\partial{x_\dimN}}\Big(\hat{g}\big(x,\tau,\Phi^{-1}(w)\big)(w-K)_\pm\psi^2\Big)\Big|\intd x\intd\tau.
\end{align*}
Hence, drawing on \eqref{eq:cond-g} and Young's inequality, we have
\begin{align*}
|I_4^\pm|&\leq \sqrt{2}\intQO\big|\partial_{x_\dimN}\hat{g}\big(x,\tau,\Phi^{-1}(w)\big)\big|(w-K)_\pm\psi^2+\sqrt{2}\intQO\big|\partial_{u}\hat{g}\big(x,\tau,\Phi^{-1}(w)\big)\big|\big|\nabla\Phi^{-1}(w)\big|(w-K)_\pm\psi^2\\&\quad+\sqrt{2}\intQO\big|\hat{g}\big(x,\tau,\Phi^{-1}(w)\big)\big|\big|\nabla(w-K)_\pm\big|\psi^2+2\sqrt{2}\intQO\big|\hat{g}\big(x,\tau,\Phi^{-1}(w)\big)\big|(w-K)_\pm\psi|\nabla\psi|\\
&\leq \sqrt{2}\intQO g_0^2(x,\tau)(w-K)_\pm\psi^2+\sqrt{2}\intQO \frac{g_0(x,\tau)\Phi'\big(\Phi^{-1}(w)\big)}{\Phi'\big(\Phi^{-1}(w)\big)}|\nabla (w-K)_\pm|(w-K)_\pm\psi^2\\
&\qquad \sqrt{2}\intQO g_0(x,\tau)w|\nabla(w-K)_\pm|\psi^2+2\sqrt{2}\intQO g_0(x,\tau) w(w-K)_\pm\psi|\nabla\psi|\\
&\leq \sqrt{2}\intQO g_0^2(x,\tau)(w-K)_\pm\psi^2+\frac{1}{4}\intQO\big|\nabla(w-K)_\pm\big|^2\psi^2+4\intQO g_0^2(x,\tau)(w-K)_\pm^2\psi^2\\
&\qquad+8\intQO g_0^2(x,\tau) w^2\psi^2\chi_{[(w-K)_\pm>0]}+\frac12\intQO(w-K)_\pm^2|\nabla\psi|^2,
\end{align*}
where $\chi_{[(w-K)_\pm>0]}$ denotes the characteristic function of the set $\{(w-K)_\pm>0\}$. Since $(w-K)_\pm\leq 2M$ and $w\leq M$, we can treat the integrals containing $g_0$ with the same methods we used for the integrals containing $a$ and $b$ in Lemma~\ref{lem:local-energy}. E.g., using two applications of Hölder
's inequality in the third integral on the right leads to
\begin{align*}
\intQO g_0^2(x,\tau)(w-K)_\pm^2\psi^2&\leq 4M^2\!\!\intQO g_0^2(x,\tau)\chi_{[(w-K)_\pm>0]}\leq\int_{T_0}^{T_1}\!\Big(\int_{B_{K,\rd}^\pm(\tau)} g_0^{2q_2}(x,\tau)\Big)^\frac{1}{q_2}\big|B_{K,\rd}^\pm(\tau)\big|^\frac{q_2-1}{q_2}\\&\leq \|g_0\|_{\LSp{2q_1}{(0,T);\Lo[2q_2]}}^2\Big(\int_{T_0}^{T_1}\big|B_{K,\rd}^\pm(\tau)\big|^\frac{q_1(q_2-1)}{q_2(q_1-1)}\Big)^\frac{q_1-1}{q_1}.
\end{align*}
Treating the remaining integrals in a similar fashion and gathering everything finally leads to \eqref{eq:loc-en-bdry}.
\end{bew}

Complementary reasoning also extends the logarithmic energy inequality to the boundary setting. Recall the notations
\begin{align*}
\xi(s):=\xi_{L,K,\delta}(s):=\Big(\ln\Big(\frac{K}{(1+\delta)K-(s-L+K)_+}\Big)\Big)_+,\quad s\in[0,L].
\end{align*}
and
\begin{align*}
\Psi(w):=\Psi_{\Phi,L,K,\delta,\rd}(w):=\Phi'\big(\Phi^{-1}(w)\big)\big(\xi^2\big)'\big(w(x,t)\big)\zeta^2(x),\quad (x,t)\in \OmT.
\end{align*}
from \eqref{eq:xi-def} and \eqref{eq:tpsi-def}, respectively.

\begin{lemma}\label{lem:log-est-bdry}
Let $M>0$. Suppose that $\romega$ is of class $C^{1}$. Assume that $\hat{g}$ satisfies \eqref{eq:cond-g}. Let $w$ be a nonnegative bounded weak solution of \eqref{eq:bvp-w} such that \eqref{eq:linfty-bound-w} holds. Then there is $C=C(\Phi,M,q_1,q_2,a,\hat{b},g_0)>0$ with the following property: For all $x_0\in\romega$, $0<T_0<T_1\leq T$ there is $\rd_0\in(0,1]$ such that for all cylinders $Q_\rd=Q_\rd(x_0,T_0,T_1)$ with $0<\rd\leq\rd_0$, all $L\in(0,M]$ satisfying $L\geq \esssup_{Q_\rd} w$, each $K\in(0,L)$, $\delta\in(0,1)$ and any cutoff function $\zeta$ inside $\mK_\rd$ fulfilling $\zeta_{\vert\partial \mK_\rd(x_0)}=0$ and $0\leq \zeta\leq 1$ in $\mK_\rd(x_0)$, the function $\xi=\xi_{L,K,\delta}$ provided by \eqref{eq:xi-def} satisfies
\begin{align}\label{eq:log-est-bdry}
\int_{\mK_\rd\cap\,\Omega}\!\!\xi^2\big(w(\cdot,T_1)\big)\zeta^2\nonumber&\leq\int_{\mK_\rd\cap\,\Omega}\xi^2\big(w(\cdot,T_0)\big)\zeta^2+12\ln\Big(\frac{1}{\delta}\Big)\iint\limits_{Q_\rd\cap\,\OmT}\Phi'\big(\Phi^{-1}(w)\big)|\nabla\zeta|^2\nonumber\\&\quad+C\Big(\frac{\ln\big(\frac1\delta\big)}{\delta K}+\frac{1+\ln\big(\frac{1}{\delta}\big)}{\delta^2 K^2}\Big)\Big(\int_{T_0}^{T_1}\!\big|B^+_{L-K,\rd}(t)\big|^{\frac{q_1(q_2-1)}{q_2(q_1-1)}}\Big)^\frac{q_1-1}{q_1}.
\end{align}
\end{lemma}

\begin{bew}
As in the previous lemma, we again assume without loss of generality that \eqref{eq:bvp-w} is written in a coordinate system satisfying the conditions of Lemma~\ref{lem:coord-transf} and that $\rd_0\in(0,1]$ and $\gamma\in\CSp{1}{\R^{\dimN-1}}$ are given accordingly. Following the testing procedures from Lemma~\ref{lem:log-est} and slightly adjusting the coefficients used in the applications of Young's inequality to keep a bit more of the well-signed terms stemming from the diffusion, we again only have to consider how to treat the newly arising boundary integral in 
\begin{align}\label{eq:log-est-bdry-eq1}
0&=\intQO\big(\xi^2(w)\zeta^2\big)_t+\intQO \big(\nabla w+ a(x,t)\Phi^{-1}(w)\big)\cdot\nabla\Psi(w)-\intQO b\big(x,t,\Phi^{-1}(w)\big)\Psi(w)\nonumber\\
&\hspace*{4cm}-\int_{T_0}^{T_1}\!\int_{\mK_\rd\cap\romega}g\big(x,t,\Phi^{-1}(w)\big)\Psi(w).
\end{align}
Recalling \eqref{eq:tpsi-grad}, we find that the well-signed terms from the diffusion are of the form
\begin{align*}
J_1:=\intQO\frac{\Phi''\big(\Phi^{-1}(w)\big)}{\Phi'\big(\Phi^{-1}(w)\big)}\xi(w)\xi'(w)\zeta^2|\nabla w|^2\quad\text{and}\quad J_2:=\intQO\Phi'\big(\Phi^{-1}(w)\big)\big(1+\xi(w)\big)\big(\xi'(w)\big)^2\zeta^2|\nabla w|^2.
\end{align*}
Repeating the re-expression procedure from the previous Lemma and using the properties of $\xi$, \eqref{eq:cond-g} and the calculations for $\nabla\Psi(w)$ in \eqref{eq:tpsi-grad}, we estimate
\begin{align*}
&\left|\int_{T_0}^{T_1}\!\!\int_{\mK_\rd\cap\romega}g\big(x,t,\Phi^{-1}(w)\big)\Psi(w)\intd S\intd t\right|\\
=\ &\left|\int_{T_0}^{T_1}\!\!\int_{\check{\mK}_\rd} \Big(g\big(x,\tau,\Phi^{-1}(w)\big)\big(w-K)_\pm\Psi(w)\Big)\Big\vert_{x_\dimN=\gamma(\cx)}\sqrt{1+|\nabla\gamma(\cx)|^2}\intd \cx\intd \tau\right|\\
\leq\ &\sqrt{2}\int_{T_0}^{T_1}\!\!\int_{\check{\mK}_\rd}\left|\int_{\gamma(\cx)}^\rd\frac{\partial}{\partial_{x_\dimN}}\Big(\hat{g}\big(x,\tau,\Phi^{-1}(w)\big)\Psi(w)\Big)\intd x_\dimN\right|\intd\cx\intd \tau\\
\leq\ &\sqrt{2}\intQO |\partial_{x_\dimN} \hat{g}|\big|\Psi(w)\big|+\sqrt{2}\intQO |\partial_u \hat{g}|\big|\Psi(w)\big|\frac{|\nabla w|}{\Phi'\big(\Phi^{-1}(w)\big)}+\sqrt{2}\intQO|\hat{g}|\big|\nabla\Psi(w)\big|\\
=\ &2\sqrt{2}\intQO |\partial_{x_\dimN} \hat{g}|\Phi'\big(\Phi^{-1}(w)\big)\xi(w)\xi'(w)\zeta^2+2\sqrt{2}\intQO |\partial_u \hat{g}| \xi(w)\xi'(w)\zeta^2|\nabla w|\\&+2\sqrt{2}\intQO|\hat{g}|\frac{\Phi''\big(\Phi^{-1}(w)\big)}{\Phi'\big(\Phi^{-1}(w)\big)}\xi(w)\xi'(w)\zeta^2|\nabla w|+2\sqrt{2}\intQO|\hat{g}|\Phi'\big(\Phi^{-1}(w)\big)\big(1+\xi(w)\big)\big(\xi'(w)\big)^2\zeta^2|\nabla w|\\&+4\sqrt{2}\intQO|\hat{g}|\Phi'\big(\Phi^{-1}(w)\big)\xi(w)\xi'(w)\zeta|\nabla\zeta|=I_1+\dots+I_5.
\end{align*}
Drawing on the estimates for $\xi$ and $\xi'$ in \eqref{eq:xi-prop}, the monotonicity properties of $\Phi'$ and $\Phi^{-1}$ and recalling that $\xi(w)=0$ on $\{w\leq L-K\}$ we obtain from Hölder's inequality that
\begin{align*}
I_1&\leq 2\sqrt{2}\Phi'\big(\Phi^{-1}(M)\big)\frac{\ln\big(\frac{1}{\delta}\big)}{\delta K}\intQO g_0^2\,\chi_{[w>L-K]}\\&\leq 2\sqrt{2}\Phi'\big(\Phi^{-1}(M)\big)\frac{\ln\big(\frac1\delta\big)}{\delta K}\|g_0\|_{\LSp{2q_1}{(0,T);\Lo[2q_2]}}^2\Big(\int_{T_0}^{T_1}\big|B^+_{L-K,\rd}(\tau)\big|^\frac{q_1(q_2-1)}{q_2(q_1-1)}\Big)^\frac{q_1-1}{q_1}.
\end{align*}
Relying on Young's inequality, we also have
\begin{align*}
I_2&\leq \tilde{\eta}_1\intQO \Phi'\big(\Phi^{-1}(w)\big)\xi(w)\big(\xi'(w)\big)^2\zeta^2|\nabla w|^2+\frac{2}{\tilde{\eta}_1}\intQO|g_0|^2\Phi'\big(\Phi^{-1}(w)\big)\xi(w)\zeta^2\\&\leq \tilde{\eta}_1 J_2+\frac{2\Phi'\big(\Phi^{-1}(M)\big)}{\tilde{\eta}_1}\frac{M\ln\big(\frac1\delta\big)}{\delta K}\|g_0\|_{\LSp{2q_1}{(0,T);\Lo[2q_2]}}^2 \Big(\int_{T_0}^{T_1}\!\!\big|B^+_{L-K,\rd}(\tau)\big|^\frac{q_1(q_2-1)}{q_2(q_1-1)}\Big)^\frac{q_1-1}{q_1},
\end{align*}
where we also made use of the fact that $\delta K\leq M$. Similarly, we find that
\begin{align*}
I_4&\leq \tilde{\eta}_2 J_2+\frac{2\Phi'\big(\Phi^{-1}(M)\big)}{\tilde{\eta}_2}\frac{1+\ln\big(\frac1\delta\big)}{\delta^2 K^2}\intQO|g_0|^2 w^2\chi_{[w>L-K]}\\
&\leq \tilde{\eta}_2 J_2+\!\frac{2\Phi'\big(\Phi^{-1}(M)\big)M^2}{\tilde{\eta}_2}\frac{1+\ln\big(\frac1\delta\big)}{\delta^2 K^2}\|g_0\|_{\LSp{2q_1}{(0,T);\Lo[2q_2]}}^2\Big(\int_{T_0}^{T_1}\!\!\big|B^+_{L-K,\rd}(\tau)\big|^\frac{q_1(q_2-1)}{q_2(q_1-1)}\Big)^\frac{q_1-1}{q_1}
\end{align*}
and
\begin{align*}
I_5&\leq 4\tilde{\eta}_3\intQO\Phi'\big(\Phi^{-1}(w)\big)\xi(w)|\nabla\zeta|^2+\frac{2}{\tilde{\eta}_3}\intQO|g_0|^2w^2\Phi'\big(\Phi^{-1}(w)\big)\xi(w)\big(\xi'(w)\big)^2\zeta^2\\
&\leq 4\tilde{\eta}_3\ln\Big(\frac1\delta\Big)\intQO\Phi'\big(\Phi^{-1}(w)\big)|\nabla\zeta|^2\\&\qquad+\frac{2\Phi'\big(\Phi^{-1}(M)\big)M^2}{\tilde{\eta}_3}\frac{\ln\big(\frac1\delta\big)}{\delta^2 K^2}\|g_0\|_{\LSp{2q_1}{(0,T);\Lo[2q_2]}}^2\Big(\int_{T_0}^{T_1}\!\!\big|B^+_{L-K,\rd}(\tau)\big|^\frac{q_1(q_2-1)}{q_2(q_1-1)}\Big)^\frac{q_1-1}{q_1}.
\end{align*}
For $I_3$ first note that due to the mean-value theorem, $\Phi^{-1}(0)=0$ and the monotonicity properties of $\Phi'$ and $\Phi^{-1}$ we have $$\frac{1}{\Phi'(\Phi^{-1}(2M))}\leq\frac{\Phi^{-1}(s)}{s}\quad \text{for all }s\in(0,2M),$$ so that with $C_1=C_1(\Phi,M):= M\Phi'\big(\Phi^{-1}(2M)\big)$ we certainly have $s^2\leq  C_1\Phi^{-1}(s)$ for $s\in[0,M]$. Moreover, recall that by Lemma~\ref{lem:phi} there is $C_2=C_2(\Phi,M)>0$ such that $$\Phi^{-1}(s)\frac{\Phi''\big(\Phi^{-1}(s)\big)}{\Phi'\big(\Phi^{-1}(s)\big)}\leq C_2\quad\text{for all }s\in[0,M].$$ Hence, we can treat $I_3$ with Young's inequality in the following way:
\begin{align*}
I_3&\leq \tilde{\eta}_4\intQO\frac{\Phi''\big(\Phi^{-1}(w)\big)}{\Phi'\big(\Phi^{-1}(w)\big)}\xi(w)\xi'(w)\zeta^2|\nabla w|^2+\frac{2}{\tilde{\eta}_4}\intQO|g_0|^2 w^2\frac{\Phi''\big(\Phi^{-1}(w)\big)}{\Phi'\big(\Phi^{-1}(w)\big)}\xi(w)\xi'(w)\zeta^2\\
&\leq \tilde{\eta}_4\intQO\frac{\Phi''\big(\Phi^{-1}(w)\big)}{\Phi'\big(\Phi^{-1}(w)\big)}\xi(w)\xi'(w)\zeta^2|\nabla w|^2+\frac{2C_1}{\tilde{\eta}_4}\!\intQO|g_0|^2 \Phi^{-1}(w)\frac{\Phi''\big(\Phi^{-1}(w)\big)}{\Phi'\big(\Phi^{-1}(w)\big)}\xi(w)\xi'(w)\zeta^2\\
&\leq \tilde{\eta}_4\intQO\frac{\Phi''\big(\Phi^{-1}(w)\big)}{\Phi'\big(\Phi^{-1}(w)\big)}\xi(w)\xi'(w)\zeta^2|\nabla w|^2+\frac{2C_1C_2}{\tilde{\eta}_4}\frac{\ln\big(\frac{1}{\delta}\big)}{\delta K}\!\intQO|g_0|^2 \chi_{[w>L-K]}\\
&\leq \tilde{\eta}_4 J_1+\frac{2C_1C_2}{\tilde{\eta}_4}\frac{\ln\big(\frac1\delta\big)}{\delta K}\|g_0\|_{\LSp{2q_1}{(0,T);\Lo[2q_2]}}^2\Big(\int_{T_0}^{T_1}\!\big|B^+_{L-K,\rd}(\tau)\big|^\frac{q_1(q_2-1)}{q_2(q_1-1)}\Big)^\frac{q_1-1}{q_1}.
\end{align*}
Gathering the estimates above and treating the non-boundary integrals in \eqref{eq:log-est-bdry-eq1} as in Lemma~\ref{lem:log-est} -- this time choosing $\eta_1=\frac14$, $\eta_2=\frac12$, $\eta_3=\frac14$, $\eta_4=\frac14$ -- we can pick $\tilde{\eta}_1=\frac14$, $\tilde{\eta}_2=\frac14$, $\tilde{\eta}_3=\frac12$, $\tilde{\eta}_4=1$ and finally conclude \eqref{eq:log-est-bdry}.
\end{bew}

To fully manage cubes with contorted bottom one additional modification has to be undertaken in Lemma~\ref{lem:DeGiorgi}. In deformed cubes as provided by Lemma~\ref{lem:coord-transf} the conclusion of the lemma remains valid due to the assumption that $\gamma\in\CSp{1}{\R^{\dimN-1}}$ satisfies $|\nabla\gamma|<\frac{1}{2\sqrt{\dimN}}$. As we will witness below this is basically a consequence of the fact that Lemma~\ref{lem:DeGiorgi} holds for convex sets as stated in \cite[Remark~2.1]{BenedettoGianazzaVespri-Harnack}. To employ the convex version in the proof of the following lemma, we will alter the shape of the given cube by an explicit coordinate transform. To better distinguish between deformed (possibly non-convex) cubes and flat cubes we will use the abbreviation $$\mK_\rd(x_0,f):=\big\{x\in\R^N\,\vert\,\cx\in\check{\mK}_\rd(x_0),\, f(\cx)<x_\dimN<\rd\big\}$$ for a cube with bottom described by $f\in\CSp{1}{\R^{\dimN-1}}$.
\begin{lemma}\label{lem:deformed-DeGiorgi}
Let $\dimN\geq1$ be an integer. There exists a constant $C=C(\dimN)>0$ with the following property: If $x_0\in\R^{\dimN}$, $\rd>0$ and $\gamma\in\CSp{1}{\R^{\dimN-1}}$ are such that $\gamma(\cx_0)=0$ and $|\nabla \gamma|\leq\frac{1}{2\sqrt{\dimN}}$ in $\check{\mK}_{\rd}(x_0)$, and if $L,K\in\R$ satisfy $L>K$ and $\varphi\in W^{1,1}(\mK_{\rd}(x_0,\gamma))$, then 
\begin{align*}
(L-K)\big|\big\{x\in \mK_\rd(x_0,\gamma)&\,:\,\varphi(x)>L\big\}\big|\\&\leq \frac{C\rd^{\dimN+1}}{\big|\big\{x\in \mK_\rd(x_0,\gamma)\,:\,\varphi(x)<K\big\}\big|}\int_{\big\{x\in \mK_\rd(x_0,\gamma)\,:\,K<\varphi(x)<L\big\}}\big|\nabla\varphi(x)\big|\intd x.
\end{align*}
\end{lemma}
\begin{bew}
From $\gamma(x_0)=0$ and $|\nabla\gamma|\leq\frac{1}{2\sqrt{\dimN}}$ in $\check{\mK}_{\rd}(x_0)$ we find in view of the mean value theorem that
\begin{align}\label{eq:deformed-DeGiorgi-eq1}
|\gamma(\cx)|=|\gamma(\cx)-\gamma(\cx_0)|\leq\max_{\check{\xi}\in\check{\mK}_{\rd}(x_0)}\big|\nabla\gamma(\check{\xi})\big||\cx-\cx_0|\leq \frac{\rd}{2} \quad\text{for all }\cx\in\check{\mK}_{\rd}(x_0).
\end{align}
We note that with $F:\mK_{\rd}(x_0,\gamma)\to\R^\dimN$ given by $F(x)=\big(x_1,\dots,x_{\dimN-1},\rd\frac{x_\dimN-\gamma(\cx)}{\rd-\gamma(\cx)}\big)^T$ the image of the deformed cube under $F$ is a flat half-cube, i.e. $F\big(\mK_\rd(x_0,\gamma)\big)=\mK_{\rd}(x_0,0)$. Moreover, $F$ is invertible on $\mK_\rd(x_0,\gamma)$ and the inverse $F^{-1}:\mK_\rd(x_0,0)\to\mK_\rd(x_0,\gamma)$ is explicitly given by $F^{-1}(y)=\big(y_1,\dots,y_{\dimN-1},\frac{\rd-\gamma(\check{y})}{\rd}y_\dimN+\gamma(\check{y})\big)^T$. Both $F$ and $F^{-1}$ are continuously differentiable and, recalling \eqref{eq:deformed-DeGiorgi-eq1}, the Jacobian determinants satisfy 
$$\frac12\leq \big|\det DF(x)\big|=\frac{\rd}{\rd-\gamma(\cx)}\leq 2\quad\text{for all }x\in\mK_\rd(x_0,\gamma)$$
and
$$\frac12\leq \big|\det DF^{-1}(y)\big|=\frac{\rd-\gamma(\check{y})}{\rd}\leq 2\quad\text{for all }y\in\mK_\rd(y_0,0),$$ respectively. With the coordinate transform given by $y=F(x)$ and the Jacobian determinant bounds we now find that
\begin{align}\label{eq:deformed-DeGiorgi-eq2}
\big|\big\{x\in\mK_\rd(x_0,\gamma):\varphi(x)>L\big\}\big| 
&=\int_{\{y\in\mK_\rd(x_0,0)\,:\,\varphi(F^{-1}(y))>L\}}\big|\det DF^{-1}(y)\big|\intd y\nonumber\\
&\leq 2\Big|\big\{y\in\mK_\rd(x_0,0):\varphi\big(F^{-1}(y)\big)>L\big\}\Big|
\end{align}
and, along similar lines,
\begin{align}\label{eq:deformed-DeGiorgi-eq3}
\Big|\big\{y\in\mK_\rd(x_0,0)\,:\,\varphi\big(F^{-1}(y)\big)<K\big\}\Big|\geq\frac{1}{2}\big|\big\{x\in\mK_\rd(x_0,\gamma)\,:\,\varphi(x)<K\big\}\big|.
\end{align}
Since $\mK_\rd(x_0,0)$ is convex, we may draw on Lemma~\ref{lem:DeGiorgi} (together with \cite[Remark~2.1]{BenedettoGianazzaVespri-Harnack}) to find $C_1=C_1(\dimN)>0$ such that
\begin{align*}
&(L-K)\Big|\big\{y\in\mK_\rd(x_0,0):\varphi\big(F^{-1}(y)\big)>L\big\}\Big|\\\leq\ &\frac{C_1\rd^{\dimN+1}}{\big|\big\{y\in\mK_\rd(x_0,0):\varphi(F^{-1}(y))<K\big\}\big|}\int_{\{y\in\mK_\rd(x_0,0)\,:\,L<\varphi(F^{-1}(y))<K\}}\big|\nabla\varphi\big(F^{-1}(y)\big)\big|\intd y,
\end{align*}
wherein substituting back and using the upper bound on $|\det DF(x)|$ entails
\begin{align}\label{eq:deformed-DeGiorgi-eq4}
&(L-K)\Big|\big\{y\in\mK_\rd(x_0,0):\varphi\big(F^{-1}(y)\big)>L\big\}\Big|\nonumber\\\leq\ &\frac{2C_1\rd^{\dimN+1}}{\big|\big\{y\in\mK_\rd(x_0,0)\,:\,\varphi(F^{-1}(y))<K\big\}\big|}\int_{\{x\in\mK_\rd(x_0,\gamma)\,:\,L<\varphi(x)<K\}}\big|\nabla\varphi(x)\big|\intd x.
\end{align}
Combining \eqref{eq:deformed-DeGiorgi-eq4} with the estimates in \eqref{eq:deformed-DeGiorgi-eq2} and \eqref{eq:deformed-DeGiorgi-eq3} we conclude upon choosing $C=8C_1$.
\end{bew}

\begin{remark}
With the augmented inequalities above, we can now proceed to obtain Hölder regularity up to the lateral boundary by replacing the cylinders appearing in the proofs of Sections~\ref{sec4} and \ref{sec5} with their intersection with $\OmT$ wherever necessary and drawing on the inequalities from Lemmas~\ref{lem:loc-en-bdry}, \ref{lem:log-est-bdry} and \ref{lem:deformed-DeGiorgi} instead of Lemma~\ref{lem:local-energy}, Lemma~\ref{lem:log-est} and Lemma~\ref{lem:DeGiorgi}, respectively. Except for the additional dependencies of $\lambda$, $\eta$ and $C$ on $g_0$ in \eqref{eq:it1-const} our iteration procedures in Section~\ref{sec4:iter} remain unchanged.
\end{remark}

\subsection{Local energy estimates near \texorpdfstring{$t=0$}{t=0}}
In contrast to the backwards-in-time cylinders used in the arguments before, we will make use of forward-in-time cylinders to extend the obtained regularity up to the starting time $t=0$. Accordingly, arguments relying on cutoff functions disappearing on $\partial_p Q_\rd$ will not provide any reasonable information for this cause. Instead we will restrict the allowed levels $k$ to get rid of the terms containing initial data. For this purpose the range of permissible levels is dependent on the size of the initial data in the corresponding box centered at $x_0$.\smallskip

Adapting the proofs of Lemma~\ref{lem:local-energy} (with a time-constant cutoff function) and Lemma~\ref{lem:loc-en-bdry} we can show the following. (See also \cite[Section II.4-(iii)]{dibenedettoDegenerateParabolicEquations1993}.)

\begin{lemma}\label{lem:loc-en-init}
Let $M>0$. Suppose that $\romega$ is of class $C^{1}$. Assume that $\hat{g}$ satisfies \eqref{eq:cond-g} and that $w_0\in\CSp{\beta_0}{\bomega}$ for some $\beta_0\in(0,1)$ is nonnegative. Let $w$ be a nonnegative bounded weak solution of \eqref{eq:bvp-w} such that \eqref{eq:linfty-bound-w} holds. There is $C=C(\Phi,M,q_1,q_2,a,\hat{b},g_0)>0$ with the following property: For all $x_0\in\bomega$ and all $0<T_1\leq T$ there is $R_0\in(0,1]$ such that for all $0<\rd\leq R_0$ and each $K\in(0,2M]$ satisfying
\begin{align}\label{eq:restrict-K}
\begin{cases}
K\geq \esssup_{\mK_\rd(x_0)}w_0\qquad&\text{for the function }\Lambda^+,\\
K\leq \essinf_{\mK_\rd(x_0)}w_0\qquad&\text{for the function }\Lambda^-,
\end{cases}
\end{align}
and every smooth cutoff function $\zeta=\zeta(x)$ inside $\mK_\rd(x_0)$ fulfilling $0\leq \zeta\leq 1$ and $\zeta_{\vert\partial\mK_\rd}=0$ the inequality
\begin{align*}
\int\limits_{\mK_\rd\cap\,\Omega}\!\!\!\Lambda^\pm\big(w(\cdot,T_1)\big)\zeta^2&+\frac{1}{4}\intQO\! \big|\nabla\big((w-K)_\pm\zeta\big)\big|^2\leq 2\!\!\intQO\!(w-K)_\pm^2|\nabla\zeta|^2+C\Big(\int_{0}^{T_1}\!\!|B_{K,\rd}^\pm(t)|^{\frac{q_1(q_2-1)}{q_2(q_1-1)}}\Big)^\frac{q_1-1}{q_1}
\end{align*}
holds for the cylinder $Q_\rd=Q_\rd(x_0,0,T_1)$. The functions $\Lambda^\pm$ are as defined in Lemma \ref{lem:Lambda}.
\end{lemma}

\begin{bew}
Keeping close to the reasoning of Lemmas~\ref{lem:local-energy} and \ref{lem:loc-en-bdry} we test \eqref{eq:bvp-w} against $\pm(w-K)_\pm\zeta^2$ and only note that in view of the level restrictions in \eqref{eq:restrict-K} we find that $K\geq \esssup_{\mK_\rd}w_0$ entails $(w_0-K)_+=0$ in $\mK_\rd\cap\,\Omega$ and that $K\leq \essinf_{\mK_\rd}w_0$ implies $(w_0-K)_-=0$ in $\mK_\rd\cap\,\Omega$, so that in both of the cases we have $\Lambda^\pm(w_0)=0$ in $\mK_\rd\cap\,\Omega$. The remaining terms can be estimated as before.
\end{bew}

Recall from \eqref{eq:xi-def}, that for $L>0$, $K\in(0,L]$ and $\delta\in(0,1)$ we defined
\begin{align*}
\xi(s):=\xi_{L,K,\delta}(s):=\Big(\ln\Big(\frac{K}{(1+\delta)K-(s-L+K)_+}\Big)\Big)_+,\quad s\in[0,L]
\end{align*}
and
\begin{align*}
\Psi(w):=\Psi_{\Phi,L,K,\delta,\rd}(w):=\Phi'\big(\Phi^{-1}(w)\big)\big(\xi^2\big)'\big(w(x,t)\big)\zeta^2(x),\quad(x,t)\in\OmT.
\end{align*}
Under the level restriction $L-K\geq\esssup_{\mK_\rd(x_0)}w_0$ we also obtain the following version of Lemmas~\ref{lem:log-est} and \ref{lem:log-est-bdry}.

\begin{lemma}\label{lem:log-est-init}
Let $M>0$. Suppose that $\romega$ is of class $C^{1}$. Assume that $\hat{g}$ satisfies \eqref{eq:cond-g} and that $w_0\in\CSp{\beta_0}{\bomega}$ for some $\beta_0\in(0,1)$ is nonnegative. Let $w$ be a nonnegative bounded weak solution of \eqref{eq:bvp-w} such that \eqref{eq:linfty-bound-w} holds. There is $C=C(\Phi,M,q_1,q_2,a,\hat{b},g_0)>0$ with the following property: For all $x_0\in\bomega$ and all $0<T_1\leq T$ there is $R_0\in(0,1]$ such that for all $0<\rd\leq R_0$, all $L\in(0,M]$ and $K\in(0,L]$ such that $$L\geq\esssup_{Q_r\cap\,\OmT} w\quad{ and }\quad L-K\geq \esssup_{\mK_\rd(x_0)}w_0,$$
each $\delta\in(0,1)$ and every smooth cutoff function $\zeta=\zeta(x)$ inside $\mK_\rd(x_0)$ fulfilling $0\leq \zeta\leq 1$ and $\zeta_{\vert\partial\mK_\rd}=0$  the function $\xi=\xi_{L,K,\delta}$ provided by \eqref{eq:xi-def} satisfies
\begin{align}\label{eq:log-est-init}
&\int_{\mK_\rd\cap\,\Omega}\!\!\xi^2\big(w(\cdot,T_1)\big)\zeta^2\\\leq\ &12\ln\Big(\frac{1}{\delta}\Big)\iint\limits_{Q_\rd\cap\,\OmT}\Phi'\big(\Phi^{-1}(w)\big)|\nabla\zeta|^2 +C\Big(\frac{1+\ln\big(\frac{1}{\delta}\big)}{\delta^2 K^2}+\frac{\ln\big(\frac1\delta\big)}{\delta K}\Big)\Big(\int_{0}^{T_1}\!\big|B^+_{L-K,\rd}(t)\big|^{\frac{q_1(q_2-1)}{q_2(q_1-1)}}\Big)^\frac{q_1-1}{q_1},\nonumber
\end{align}
where $Q_\rd$ is given by $Q_\rd=Q_\rd(x_0,0,T_1)$.
\end{lemma}

\begin{bew}
Note that because of $L-K\geq \esssup_{\mK_{\rd}(x_0)}w_0$ we have $\big(w_0(x)-L+K\big)_+= 0$ for all $x\in\mK_\rd\cap\,\Omega$ and accordingly $\xi(w_0)= 0$ throughout $\mK_\rd\cap\,\Omega$. Hence, closely following the proofs of Lemma~\ref{lem:log-est} and Lemma~\ref{lem:log-est-bdry} we may prove \eqref{eq:log-est-init}.
\end{bew}

Since the energy inequalities remain valid for the restricted levels even when considering cylinders lying at the bottom of $\OmT$, we can extend the two alternatives previously discussed in Sections~\ref{sec4:alt1} and \ref{sec4:alt2} by replacing the cylinders appearing in these results with their corresponding intersections with $\overline{\OmT}$. We let $\bOmT:=\bomega\times[0,T]$ and for $\sig<1$ set $$Q_{\rd}^{\sig,\theta}:=Q_{\rd}\big(x_0,(\sig-1)\theta\rd^2,\sig\theta\rd^2\big).$$ Notice that $\sig-1<0$, so that 
\begin{align}\label{eq:size-of-shortenend-cyl}
\big|Q_{\rd}^{\sig,\theta}\cap\bOmT\big|\leq\sig\omega_\dimN\theta\rd^{\dimN+2}.
\end{align}
With these notations we obtain a version of Lemma~\ref{lem:first-alt-1} under the restriction $\mu^-+k\leq\essinf_{\mK_\rd} w_0$ of the following form.

\begin{lemma}\label{lem:first-alt-init}
Let $M>0$. Suppose that $\romega$ is of class $C^{1}$. Assume that $\hat{g}$ satisfies \eqref{eq:cond-g} and that $w_0\in\CSp{\beta_0}{\bomega}$ for some $\beta_0\in(0,1)$ is nonnegative. There is $\sig_0=\sig_0(\Phi,M,q_1,q_2,a,\hat{b},g_0)\in(0,\tfrac12)$ such that for all nonnegative bounded weak solutions $w$ of \eqref{eq:bvp-w} satisfying \eqref{eq:linfty-bound-w} the following holds: Given $x_0\in\Omega$ one can find $R_0\in(0,1]$ such that for all $0<\rd\leq R_0$, all $\sig\in(0,1)$, each $0<\vartheta$, letting $\mu^-:=\essinf_{Q_{2\rd}^{\sig,\vartheta}}w<\esssup_{Q_{2\rd}^{\sig,\vartheta}}w=:\mu^+$ and $k\in(0,M]$ be such that $\mu^-+k\leq\essinf_{\mK_{2r}} w_0$, then if $\vartheta\geq\theta:=\min\big\{\theta_k,\theta_{\mu^-}\big\}$ and
\begin{align}\label{eq:first-alt-init-ass1}
\big|\big\{(x,t)\in Q_{\rd}^{\sig,\theta}\cap\bOmT\,:\,(w- \mu^--k)_->0\big\}\big|\leq \sig_0\omega_\dimN\theta\rd^{\dimN+2} 
\end{align}
then either
\begin{align*}
k^{2}\theta^{-\frac{q_1-1}{q_1}}\leq \rd^{\dimN\kappa}
\end{align*}
or
\begin{align*}
\Big|\Big\{(x,t)\in Q_{\frac{\rd}{2}}^{\sig,\theta}\cap\bOmT\,:\,\big(w-\mu^--\tfrac{k}{2}\big)_->0\Big\}\Big|=0.
\end{align*}
\end{lemma}

\begin{remark}
Assuming $\sig\leq\sig_0<\frac12$, the condition \eqref{eq:first-alt-init-ass1} is trivially satisfied according to \eqref{eq:size-of-shortenend-cyl}.
\end{remark}

Similarly, assuming the levels to satisfy the constriction $\mu^+-k\geq\esssup_{\mK_{2r}}w_0$, we can start with the logarithmic inequality of Lemma~\ref{lem:log-est-init} and arguments akin to Lemma~\ref{lem:expansion-in-time} to first find $\delta_0\in(0,\frac12)$ and $e_0\in(0,1)$ such that 
$$\big|B_{\mu^+-\delta_0 k,\rd}^+(t)\big|\leq e_0|\mK_\rd|\quad\text{for all }t\in(0,\sig_0\theta\rd^2].$$ Note that a smallness requirement on the level set as in \eqref{eq:exp-in-time} of Lemma~\ref{lem:expansion-in-time} is not required in this case, because in contrast to Lemma~\ref{lem:log-est} there is no evaluation at $T_0$ in Lemma~\ref{lem:log-est-init} due to the level restrictions. After this, we can continue by employing Lemma~\ref{lem:loc-en-init} for $\Lambda^+$ within a reasoning along the lines of Lemmas~\ref{lem:expansion-in-space} and \ref{lem:second-alt} enables us to establish a result of the following form.

\begin{lemma}\label{lem:second-alt-init}
Let $M>0$. Suppose that $\romega$ is of class $C^{1}$. Assume that $\hat{g}$ satisfies \eqref{eq:cond-g} and that $w_0\in\CSp{\beta_0}{\bomega}$ for some $\beta_0\in(0,1)$ is nonnegative. Denote by $\sig_0=\sig_0(\Phi,M,q_1,q_2,a,\hat{b},g_0)\in(0,\frac12)$ and $R_0\in(0,1]$ the numbers form Lemma~\ref{lem:first-alt-init}. There exists $\deld=\deld(\Phi,M,q_1,q_2,a,\hat{b},g_0)\in(0,\tfrac14)$ such that for all nonnegative bounded weak solutions $w$ of \eqref{eq:bvp-w} satisfying \eqref{eq:linfty-bound-w} the following holds: Given $x_0\in\Omega$ one can find $R_0\in(0,1]$ such that for all $0<\rd\leq R_0$, each $0<\vartheta$, letting $\mu^-:=\essinf_{Q_{2\rd}^{\sig_0,\vartheta}}w<\esssup_{Q_{2\rd}^{\sig_0,\vartheta}}w=:\mu^+$ and $k\in(0,M]$ be such that $\mu^+-k\geq\esssup_{\mK_{2r}} w_0$, then if $$\vartheta\geq\theta:=\min\big\{\theta_k,\theta_{\mu^-}\big\}\quad\text{and}\quad \mu^+-\mu^-\leq 4k$$ then either
\begin{align*}
\deld^2 k^{2}\theta^{-\frac{q_1-1}{q_1}}\leq \rd^{\dimN\kappa}
\end{align*}
or
\begin{align*}
\Big|\Big\{(x,t)\in Q_{\frac{\rd}{2}}^{\sig_0,\theta}\cap\bOmT\,:\,\big(w-\mu^++\tfrac{\deld}{2}k\big)_+>0\Big\}\Big|=0.
\end{align*}
\end{lemma}

With these adjusted alternatives we can include marginal modifications to the iteration procedure, in order to accommodate for cylinders lying at the bottom of $\Omega\times(0,T)$. For this, we will slightly restructure the cases considered in the induction step according to the level restrictions.
\begin{lemma}\label{lem:iteration-1-init}
Let $M>0$. Suppose that $\romega$ is of class $C^{1}$. Assume that $\hat{g}$ satisfies \eqref{eq:cond-g} and that $w_0\in\CSp{\beta_0}{\bomega}$ for some $\beta_0\in(0,1)$ is nonnegative. Let $w$ be a nonnegative bounded weak solution of \eqref{eq:bvp-w} such that \eqref{eq:linfty-bound-w} holds. Denote by $\sig_0=\sig_0(\Phi,M,q_1,q_2,a,\hat{b},g_0)\in(0,\frac{1}{2})$ and $R_0\in(0,1]$ the numbers from Lemma~\ref{lem:first-alt-init}. There are $\lambda=\lambda(\Phi,M)\in(0,\tfrac14)$, $\eta(\Phi,M,q_1,q_2,a,\hat{b},g_0)\in(\tfrac12,1)$, $\beta=\beta(q_1,q_2)\in(0,1)$, $C=C(\Phi,M,q_1,q_2,a,\hat{b},g_0)\geq 1$ with the following property: Given $x_0\in\bomega$, $\rd_0\leq R_0$ and $0<k_0\leq\frac{M}{2}$ such that 
$$\essosc_{Q_{\rd_0}^{\sig_0,\theta_{k_0}}\cap\,\bOmT}w\leq 4k_0,$$
 one can find non-increasing sequences $(\rd_j)_{j\in\N}$, $(k_j)_{j\in\N}$, and an increasing sequence $(\theta_j)_{j\in\N}$ such that with $\mK_j:=\mK_{\rd_j}(x_0)\cap\,\bomega$ and $Q_j:=Q_{\rd_j}^{\sig_0,\theta_j}\cap\,\bOmT$ the properties
\begin{align}\label{eq:iteration-1-init}
Q_j&\subseteq Q_{j-1},&\essosc_{Q_{j}} w&\leq \max\Big\{4k_j,8\essosc_{\mK_{j-1}}w_0\Big\},\nonumber\\
\frac{1}{2}k_{j-1}&\leq k_j\leq \max\big\{\eta k_{j-1}, C\rd_{j-1}^\beta\big\}\ \;\text{and}\ \; &\lambda\rd_{j-1}&\leq\rd_j\leq\frac{1}{2}\rd_{j-1}
\end{align}
either hold for all $j\in\N$ or there is $\js\geq0$ such that \eqref{eq:iteration-1-init} holds for all $j\leq \js$ and $\mu^{-}_{\js}:=\essinf_{Q_{\js}} w>k_{\js}$.
\end{lemma}

\begin{bew}
With $\rd_0$, $k_0$ and $\theta_{k_0}=\big(\Phi'(\Phi^{-1}(k_0))\big)^{-1}$ provided we set $Q_0:=Q_{\rd_0}^{\sig_0,\theta_{k_0}}\cap\bOmT$, $\mu^-_{0}:=\essinf_{Q_0} w$ and $\mu^+_{0}:=\esssup_{Q_0}w$. If $\mu^-_{0}>k_0$ we set $j_\star=j_0$ and have nothing to proof. Assume now $\mu^-_{0}\leq k_0$. For $j>0$ we denote by $\mu^-_{j}:=\essinf_{Q_j} w$ and $\mu^+_{j}:=\esssup_{Q_j}w$ the infimum and supremum of $w$ in the current cylinder $Q_j$ and by $\mu^-_{j,0}:=\essinf_{\mK_{j}} w_0$ and $\mu^+_{j,0}:=\esssup_{\mK_{j}}w_0$ the infimum and supremum of the initial datum in the corresponding spatial box centered at $x_0$. If at some point $\mu^-_{j}>k_j$ the argument is finished. If not, we will generate $\rd_{j+1}$, $k_{j+1}$, $\theta_{j+1}$ depending on which of the different cases we find ourselves in. Let $C_1=C_1(\Phi,2,M)>1$ again be taken from Lemma~\ref{lem:phi} such that
\begin{align*}
\Phi'\big(\Phi^{-1}(s)\big)\leq C_1\Phi'\big(\tfrac12\Phi^{-1}(s)\big)\quad\text{for all }s\in[0,M],
\end{align*}
and let $\deld=\deld(\Phi,M,q_1,q_2,a,\hat{b},g_0)\in(0,\tfrac14)$ be provided by Lemma~\ref{lem:second-alt-init}. We claim that the lemma holds for the choices
\begin{equation*}
\begin{split}
\lambda&:=\frac{1}{4\sqrt{C_1}}\in(0,\tfrac14),\qquad \eta:=\big(1-\tfrac{\deld}{32}\big)\in(\tfrac{31}{32},1),\qquad\beta:=\frac{\dimN\kappa}{2+\frac{q_1-1}{q_1}}<1\\ \text{and} \quad C&:=\max\Big\{1,\Big(\big(\Phi^{-1}(M)\big)^{\frac{q_1-1}{q_1}}\deld^{-2}\Big)^\frac{1}{2+\frac{q_1-1}{q_1}}\Big\}\geq 1.
\end{split}
\end{equation*}
We will now construct the next cylinder in accordance with the following cases. Since the properties of $\rd_{j+1}$, $k_{j+1}$ and the inclusion of $Q_{j+1}$ in $Q_j$ can be easily adapted from Lemma~\ref{lem:iteration-1}, we will mostly only present details for the inequality with the essential oscillation in \eqref{eq:iteration-1-init}.\smallskip

\emph{Case 1}: $\essosc_{Q_j} w<2k_j$.\smallskip\\
Set $\rd_{j+1}:=\frac{\rd_j}{2\sqrt{C_1}}$, $k_{j+1}:=\frac{k_j}{2}$, $\theta_{j+1}=\theta_{k_{j+1}}$. Then, due to the choice of $C_1$, we have $Q_{j+1}\subset Q_j$ and thus
\begin{align*}
\essosc_{Q_{j+1}} w\leq \essosc_{Q_j} w< 2k_j=4k_{j+1}\leq \max\Big\{4k_{j+1},8\essosc_{K_{j}}w_0\Big\}.
\end{align*}\vspace*{-0.25cm}

\emph{Case 2}: $2k_j\leq \essosc_{Q_j} w\leq 4k_j$.\smallskip\\
Here, we further separate the case according to the size difference between $k_j$ and the oscillation of $w_0$ in $\mK_j$.\smallskip

\emph{Case 2.1}: $k_j< 2\essosc_{\mK_j}w_0$.\smallskip\\
Set $\rd_{j+1}:=\frac{\rd_j}{2\sqrt{C_1}}$, $k_{j+1}:=\frac{k_j}{2}$, $\theta_{j+1}=\theta_{k_{j+1}}$. Again $Q_{j+1}\subset Q_j$ and hence
\begin{align*}
\essosc_{Q_{j+1}} w\leq \essosc_{Q_j} w \leq 4 k_j<8 \essosc_{\mK_j}w_0\leq\max\Big\{4k_{j+1},8\essosc_{\mK_j}w_0\Big\}.
\end{align*}\vspace*{-0.25cm}

\emph{Case 2.2}: $2\essosc_{\mK_j}w_0\leq k_j$.\smallskip\\
In this subcase, we also split the case with regards to the level restrictions of Lemma~\ref{lem:loc-en-init}.\smallskip

\emph{Case 2.2a}: $\mu_j^+-\frac{k_j}{4}\leq \mu_{j,0}^+$ and $\mu_j^-+\frac{k_j}{4}\geq\mu_{j,0}^-$.\smallskip\\
This case cannot occur. Notice that substracting the second inequality of Case 2.2a from the first would entail that
\begin{align*}
\essosc_{Q_j} w\leq\essosc_{K_j}w_0+\frac{k_j}{2}\stackrel{\text{Case 2.2}}{\leq} k_j,
\end{align*}
which contradicts the assumption of Case 2.\smallskip

\emph{Case 2.2b}: $\mu_j^-+\frac{k_j}{4}<\mu_{j,0}^-$.\smallskip\\
In this case, $K:=\mu_j^-+\frac{k_j}{4}<\mu_{j,0}^-$ fulfills \eqref{eq:restrict-K} for $\Lambda^-$ and 
Lemma~\ref{lem:first-alt-init} becomes applicable for the choices $\sig=\sig_0$, $k=\frac{k_j}{4}$ and $\rd=\frac{\rd_j}{2}$. Thus, either 
\begin{align}\label{eq:iteration-1-case-alt}
\frac{k_j^2}{16}\theta_j^{-\frac{2q_1-1}{q_2}}\leq\left(\frac{\rd_j}{2}\right)^{\dimN\kappa}\quad\text{or}\quad \Big|\Big\{(x,t)\in Q_{\frac{\rd_j}{4}}^{\sig_0,\theta_j}\cap\bOmT\,:\,\big(w-\mu_j^--\tfrac{k_j}{8}\big)_->0\Big\}\Big|=0.
\end{align}
If the first is true we pick $\rd_{j+1}=\frac{\rd_j}{2}$, $k_{j+1}=k_j$, $\theta_{j+1}=\theta_{k_{j+1}}$ and follow the reasoning of Lemma~\ref{lem:iteration-1} Case 2a. In particular, we easily see that then $$\essosc_{Q_{j+1}}w\leq 4k_j=4k_{j+1}\leq\max\Big\{4k_{j+1},8\essosc_{\mK_j}w_0\Big\}.$$
If the second part in \eqref{eq:iteration-1-case-alt} is true, we pick $\rd_{j+1}=\lambda\rd_j$, $k_{j+1}=\frac{31}{32}k_j$, $\theta_{j+1}=\theta_{k_{j+1}}$ and obtain $Q_{j+1}\subset Q_{\frac{\rd_j}{4}}^{\sig_0,\theta}$ and hence\vspace*{-0.3cm}
\begin{align*}
\essosc_{Q_{j+1}}w\leq \essosc_{Q_{\frac{\rd_j}{4}}^{\sig_0,\theta_j}} w\leq \mu^+_j-\mu^-_j-\frac{k_j}{8}\leq \frac{31}{8}k_j=4k_{j+1}\leq\max\Big\{4k_{j+1},8\essosc_{\mK_j}w_0\Big\}.
\end{align*}\vspace*{-0.25cm}

\emph{Case 2.2c}: $\mu_j^+-\frac{k_j}{4}>\mu_{j,0}^+$.\smallskip\\
Here, Lemma~\ref{lem:second-alt-init} becoms applicable for $k=\frac{k_j}{4}$ and $\rd=\frac{\rd_j}{2}$. Hence, either
\begin{align*}
\frac{\deld^2 k_j^2}{16}\theta_j^{-\frac{2q_1-1}{q_2}}\leq\left(\frac{\rd}{2}\right)^{\dimN\kappa}\quad\text{or}\quad \Big|\Big\{(x,t)\in Q_{\frac{\rd_j}{4}}^{\sig_0,\theta_j}\cap\bOmT\,:\,\big(w-\mu_j^++\tfrac{\deld}{8}k\big)_+>0\Big\}\Big|=0.
\end{align*}
If the first is true, again take $\rd_{j+1}=\frac{\rd_j}{2}$, $k_{j+1}=k_j$ and $\theta_{j+1}=\theta_{k_{j+1}}$ and follow Case 2a of Lemma~\ref{lem:iteration-1}.
To the contrary, if the second holds, pick $\rd_{j+1}=\lambda\rd_j$, $k_{j+1}=\eta k_j$ and $\theta_{j+1}=\theta_{k_{j+1}}$. Then, as in Case 2c of Lemma~\ref{lem:iteration-1}, we have $Q_{j+1}\subset Q_{\frac{\rd_j}{4}}^{\sig_0,\theta_j}$ and
\[
\essosc_{Q_{j+1}} w\leq \mu_j^+-\frac{\deld k_j}{8}-\mu^-_j\leq 4k_j-\frac{\deld k_j}{8}=4\eta k_j=4k_{j+1}\leq\max\Big\{4k_{j+1},8\essosc_{\mK_j}w_0\Big\}.
\qedhere\]
\end{bew}

In the same manner, one obtains a version of the second iteration, which then can be combined to a result akin to Corollary~\ref{cor:iteration}. Afterwards, the steps of Section~\ref{sec5} can be repeated to prove Theorem~\ref{theo:2}.

\begin{proof}[\textbf{Proof of Theorem \ref{theo:2}}:]
The proof can be completed by combining Lemma~\ref{lem:iteration-1-init} with arguments corresponding to Lemma~\ref{lem:iteration-2} and Corollary~\ref{cor:iteration} and the reasoning outlined in Section~\ref{sec5}.
\end{proof}

\section*{Acknowledgements}
The author would like to thank Frederic Heihoff (Paderborn University) and the anonymous referee for their valuable comments concerning the presentation of the argument.\\
The author acknowledges support of the {\em Deutsche~Forschungsgemeinschaft} (Project No.~462888149).


\footnotesize{
\setlength{\bibsep}{3pt plus 0.5ex}

}
\enlargethispage{1cm}
\vfill
\hfill\begin{minipage}[r]{0.3\textwidth}
Tobias Black\\
Institute of Mathematics\\
Paderborn University\\
33098 Paderborn, Germany\\
tblack@math.upb.de
\end{minipage}


\begin{thebibliography}{42}
\providecommand{\natexlab}[1]{#1}

\bibitem[Adler(1966)]{Adler66}{https://doi.org/10.1126/science.153.3737.708}
J.~Adler.
\newblock \emph{Chemotaxis in {{Bacteria}}}.
\newblock Science, {\bfseries 153} (1966), \penalty0 708--716.

\bibitem[Bellomo et~al.(2015)Bellomo, Bellouquid, Tao, and
  Winkler]{BBWT15}{https://doi.org/10.1142/S021820251550044X}
N.~Bellomo, A.~Bellouquid, Y.~Tao, and M.~Winkler.
\newblock \emph{Toward a mathematical theory of {K}eller-{S}egel models of
  pattern formation in biological tissues}.
\newblock Math. Models Methods Appl. Sci., {\bfseries 25} (2015), \penalty0
  1663--1763.

\bibitem[Black(2018)]{TB2017_nonlindiff}{https://doi.org/10.1137/17M1159488}
T.~Black.
\newblock \emph{Global {{Very Weak Solutions}} to a {{Chemotaxis-Fluid System}}
  with {{Nonlinear Diffusion}}}.
\newblock SIAM J. Math. Anal., {\bfseries 50} (2018), \penalty0 4087--4116.

\bibitem[Black(2019)]{TB2018_NA}{https://doi.org/10.1016/j.na.2018.10.003}
T.~Black.
\newblock \emph{Global Solvability of Chemotaxis-Fluid Systems with Nonlinear
  Diffusion and Matrix-Valued Sensitivities in Three Dimensions}.
\newblock Nonlinear Anal., {\bfseries 180} (2019), \penalty0 129--153.

\bibitem[B\"{o}gelein et~al.(2022)B\"{o}gelein, Duzaar, and
  Scheven]{MR4500278}{https://doi.org/10.1007/s42985-022-00204-0}
V.~B\"{o}gelein, F.~Duzaar, and C.~Scheven.
\newblock \emph{Higher integrability for doubly nonlinear parabolic systems}.
\newblock Partial Differ. Equ. Appl., {\bfseries 3} (2022), \penalty0 Paper No.
  74, 41.

\bibitem[Chagas et~al.()Chagas, Diehl, and
  Guidolin]{chagasPropertiesSteklovAverages2017}{}
J.~Q. Chagas, N.~M.~L. Diehl, and P.~L. Guidolin.
\newblock \emph{Some Properties for the {{Steklov}} Averages}, Preprint (2017),
  arXiv:1707.06368.

\bibitem[Chung et~al.(2017)Chung, Hwang, Kang, and
  Kim]{ChungHwangKangKim17}{https://doi.org/10.1016/j.jde.2017.03.042}
Y.-S. Chung, S.~Hwang, K.~Kang, and J.~Kim.
\newblock \emph{H\"older Continuity of {{Keller}}\textendash{{Segel}} Equations
  of Porous Medium Type Coupled to Fluid Equations}.
\newblock J. Differential Equations, {\bfseries 263} (2017), \penalty0
  2157--2212.

\bibitem[De~Giorgi(1957)]{DeGiorgi57}{}
E.~De~Giorgi.
\newblock \emph{Sulla differenziabilit\`a e l'analiticit\`a delle estremali
  degli integrali multipli regolari}.
\newblock Mem. Accad. Sci. Torino. Cl. Sci. Fis. Mat. Nat., {\bfseries 3}
  (1957), \penalty0 25--43.

\bibitem[DiBenedetto(1982)]{DB82}{https://doi.org/10.1007/BF01761493}
E.~DiBenedetto.
\newblock \emph{Continuity of Weak Solutions to Certain Singular Parabolic
  Equations}.
\newblock Ann. di Mat. Pura ed Appl., {\bfseries 130} (1982), \penalty0
  131--176.

\bibitem[DiBenedetto(1983)]{DiBenedetto83}{https://doi.org/10.1512/iumj.1983.32.32008}
E.~DiBenedetto.
\newblock \emph{Continuity of {{Weak Solutions}} to a {{General Porous Medium
  Equation}}}.
\newblock Indiana Univ. Math. J., {\bfseries 32} (1983), \penalty0 83--118.

\bibitem[DiBenedetto(1993)]{dibenedettoDegenerateParabolicEquations1993}{https://doi.org/10.1007/978-1-4612-0895-2}
E.~DiBenedetto.
\newblock \emph{Degenerate {{Parabolic Equations}}}.
\newblock Universitext. {Springer New York}, {New York, NY}, 1993.

\bibitem[DiBenedetto and
  Friedman(1985)]{dibenedettoHolderEstimatesNonlinear1985}{https://doi.org/10.1515/crll.1985.357.1}
E.~DiBenedetto and A.~Friedman.
\newblock \emph{H\"older Estimates for Nonlinear Degenerate Parabolic Sytems.}
\newblock J. f\"ur die Reine und Angew. Math, {\bfseries 357} (1985), \penalty0
  1--22.

\bibitem[DiBenedetto et~al.(2012)DiBenedetto, Gianazza, and
  Vespri]{BenedettoGianazzaVespri-Harnack}{https://doi.org/10.1007/978-1-4614-1584-8}
E.~DiBenedetto, U.~Gianazza, and V.~Vespri.
\newblock \emph{Harnack's {{Inequality}} for {{Degenerate}} and {{Singular
  Parabolic Equations}}}.
\newblock Springer {{Monographs}} in {{Mathematics}}. {Springer New York}, {New
  York, NY}, 2012.

\bibitem[Dombrowski et~al.(2004)Dombrowski, Cisneros, Chatkaew, Goldstein, and
  Kessler]{Dombrowski04}{https://doi.org/10.1103/PhysRevLett.93.098103}
C.~Dombrowski, L.~Cisneros, S.~Chatkaew, R.~E. Goldstein, and J.~O. Kessler.
\newblock \emph{Self-{{Concentration}} and {{Large-Scale Coherence}} in
  {{Bacterial Dynamics}}}.
\newblock Phys. Rev. Lett., {\bfseries 93} (2004).
%
\bibitem[Fiorenza(2016)]{fiorenza2016}{https://doi.org/10.1007/978-3-319-47940-8}
R.~Fiorenza.
\newblock \emph{{H\"older} and locally {H\"older} Continuous Functions, and Open Sets of Class $C^k, C^{k,\lambda}$}.
\newblock {{Frontiers}} in {{Mathematics}}, {Birkh\"auser/Springer}, 2016.

\bibitem[Fischer(2013)]{fischerAdvectionDrivenSupportShrinking2013}{https://doi.org/10.1137/120874291}
J.~Fischer.
\newblock \emph{Advection-{{Driven Support Shrinking}} in a {{Chemotaxis
  Model}} with {{Degenerate Mobility}}}.
\newblock SIAM J. Math. Anal., {\bfseries 45} (2013), \penalty0 1585--1615.

\bibitem[Gianazza et~al.(2010)Gianazza, Surnachev, and
  Vespri]{gianazzaNewProofHolder2010}{https://doi.org/10.1515/acv.2010.009}
U.~Gianazza, M.~Surnachev, and V.~Vespri.
\newblock \emph{A New Proof of the {{H\"older}} Continuity of Solutions to
  P-{{Laplace}} Type Parabolic Equations}.
\newblock Adv. Calc. Var., {\bfseries 3} (2010).

\bibitem[Hissink~Muller(2022)]{hissinkmullerInteriorHolderContinuity2022}{https://doi.org/10.1007/s00028-022-00849-9}
V.~Hissink~Muller.
\newblock \emph{Interior {{H\"older}} Continuity for Singular-Degenerate Porous
  Medium Type Equations with an Application to a Biofilm Model}.
\newblock J. Evol. Equ., {\bfseries 22} (2022), \penalty0 92.

\bibitem[Horstmann(2003)]{Ho03}{}
D.~Horstmann.
\newblock \emph{From 1970 until Present: The {{Keller-Segel}} Model in
  Chemotaxis and Its Consequences. {{I}}}.
\newblock Jahresber. Deutsch. Math.-Verein., {\bfseries 105} (2003), \penalty0
  103--165.

\bibitem[Hwang and
  Zhang(2021)]{HwangZhang19}{https://doi.org/10.1016/j.na.2021.112413}
S.~Hwang and Y.~P. Zhang.
\newblock \emph{Continuity Results for Degenerate Diffusion Equations with
  {${{L}}_t^p{{L}}_x^q$} Drifts}.
\newblock Nonlinear Anal., {\bfseries 211} (2021), \penalty0 112413.

\bibitem[Keller and
  Segel(1970)]{KS70}{https://doi.org/10.1016/0022-5193(70)90092-5}
E.~F. Keller and L.~A. Segel.
\newblock \emph{Initiation of Slime Mold Aggregation Viewed as an Instability}.
\newblock J. Theor. Biol., {\bfseries 26} (1970), \penalty0 399--415.

\bibitem[Ladyženskaja et~al.(1968)Ladyženskaja, Solonnikov, and
  Ural'ceva]{LSU}{https://doi.org/10.1090/mmono/023}
O.~A. Ladyženskaja, V.~A. Solonnikov, and N.~N. Ural'ceva.
\newblock \emph{Linear and quasilinear equations of parabolic type}.
\newblock Translations of mathematical monographs. American Mathematical
  Society, 1968.

\bibitem[Lankeit and
  Winkler(2017)]{LanWin2017}{https://doi.org/10.1007/s00030-017-0472-8}
J.~Lankeit and M.~Winkler.
\newblock \emph{A Generalized Solution Concept for the {{Keller-Segel}} System
  with Logarithmic Sensitivity: Global Solvability for Large Nonradial Data}.
\newblock NoDEA Nonlinear Differential Equations Appl., {\bfseries 24} (2017),
  \penalty0 Art. 49.

\bibitem[Lankeit and
  Winkler(2019)]{LanWin_JDMV_19}{https://doi.org/10.1365/s13291-019-00210-z}
J.~Lankeit and M.~Winkler.
\newblock \emph{Facing {{Low Regularity}} in {{Chemotaxis Systems}}}.
\newblock Jahresber. Deutsch. Math.-Verein.

\bibitem[Liao(2020)]{liaoUnifiedApproachHolder2020}{https://doi.org/10.1016/j.jde.2019.11.023}
N.~Liao.
\newblock \emph{A Unified Approach to the {{H\"older}} Regularity of Solutions
  to Degenerate and Singular Parabolic Equations}.
\newblock J. Differential Equations, {\bfseries 268} (2020), \penalty0
  5704--5750.

\bibitem[Lieberman(1996)]{lieberman}{https://doi.org/10.1142/3302}
G.~M. Lieberman.
\newblock \emph{Second Order Parabolic Differential Equations}.
\newblock {World Scientific Publishing Co., Inc., River Edge, NJ}, 1996.

\bibitem[Lorz(2010)]{lorz-M3AS10}{https://doi.org/10.1142/S0218202510004507}
A.~Lorz.
\newblock \emph{Coupled {{Chemotaxis Fluid Model}}}.
\newblock Math. Models Methods Appl. Sci., {\bfseries 20} (2010), \penalty0
  987--1004.

\bibitem[Marras et~al.(2023)Marras, Ragnedda, {Vernier-Piro}, and
  Vespri]{MRVV22}{https://doi.org/10.1016/j.jde.2023.04.013}
M.~Marras, F.~Ragnedda, S.~{Vernier-Piro}, and V.~Vespri.
\newblock \emph{H\"older Estimates of Weak Solutions to Degenerate Chemotaxis
  Systems with a Source Term}.
\newblock J. Differential Equations, {\bfseries 366} (2023), \penalty0 42--70.

\bibitem[Naumann(1984)]{MR0771666}{}
J.~Naumann.
\newblock \emph{Einf\"{u}hrung in die {T}heorie parabolischer
  {V}ariationsungleichungen}, Volume~64 of \emph{Teubner-Texte zur Mathematik
  [Teubner Texts in Mathematics]}.
\newblock BSB B. G. Teubner Verlagsgesellschaft, Leipzig, 1984.
\newblock With English, French and Russian summaries.

\bibitem[Osaki et~al.(2002)Osaki, Tsujikawa, Yagi, and
  Mimura]{Os02-chemologatract}{https://doi.org/10.1016/S0362-546X(01)00815-X}
K.~Osaki, T.~Tsujikawa, A.~Yagi, and M.~Mimura.
\newblock \emph{Exponential Attractor for a Chemotaxis-Growth System of
  Equations}.
\newblock Nonlinear Anal. Theory Methods Appl., {\bfseries 51} (2002),
  \penalty0 119--144.

\bibitem[Painter and Hillen(2002)]{HP-volumefilling-CAMQ02}{}
K.~J. Painter and T.~Hillen.
\newblock \emph{Volume-Filling and Quorum-Sensing in Models for Chemosensitive
  Movement}.
\newblock Can. Appl. Math. Q., {\bfseries 10} (2002), \penalty0 501--543.

\bibitem[Porzio and
  Vespri(1993)]{PorzVesp93}{https://doi.org/10.1006/jdeq.1993.1045}
M.~M. Porzio and V.~Vespri.
\newblock \emph{Holder estimates for local solutions of some doubly nonlinear
  degenerate parabolic equations}.
\newblock J. Differential Equations, {\bfseries 103} (1993), \penalty0
  146--178.

\bibitem[Stevens and
  Winkler(2022)]{stevensTaxisdrivenPersistentLocalization2022}{https://doi.org/10.1080/03605302.2022.2122836}
A.~Stevens and M.~Winkler.
\newblock \emph{Taxis-Driven Persistent Localization in a Degenerate
  {{Keller-Segel}} System}.
\newblock Comm. Partial Differential Equations, {\bfseries 47} (2022),
  \penalty0 2341--2362.

\bibitem[Tao and
  Winkler(2012)]{TaoWin-quasilinear_JDE12}{https://doi.org/10.1016/j.jde.2011.08.019}
Y.~Tao and M.~Winkler.
\newblock \emph{Boundedness in a Quasilinear Parabolic-Parabolic
  {{Keller-Segel}} System with Subcritical Sensitivity}.
\newblock J. Differential Equations, {\bfseries 252} (2012), \penalty0
  692--715.

\bibitem[Tello and
  Winkler(2007)]{TW07}{https://doi.org/10.1080/03605300701319003}
J.~I. Tello and M.~Winkler.
\newblock \emph{A Chemotaxis System with Logistic Source}.
\newblock Comm. Partial Differential Equations, {\bfseries 32} (2007),
  \penalty0 849--877.

\bibitem[Tuval et~al.(2005)Tuval, Cisneros, Dombrowski, Wolgemuth, Kessler, and
  Goldstein]{tuval2005bacterial}{https://doi.org/10.1073/pnas.0406724102}
I.~Tuval, L.~Cisneros, C.~Dombrowski, C.~W. Wolgemuth, J.~O. Kessler, and R.~E.
  Goldstein.
\newblock \emph{Bacterial Swimming and Oxygen Transport near Contact Lines}.
\newblock Proc. Natl. Acad. Sci. U.S.A., {\bfseries 102} (2005), \penalty0
  2277--2282.

\bibitem[V{\'a}zquez(2007)]{Vazquez-PME-07}{}
J.~L. V{\'a}zquez.
\newblock \emph{The Porous Medium Equation}.
\newblock Oxford {{Mathematical Monographs}}. {The Clarendon Press, Oxford
  University Press, Oxford}, 2007.

\bibitem[Wang and
  Li(2017)]{WangLi-ZAMP17}{https://doi.org/10.1007/s00033-017-0773-0}
Y.~Wang and X.~Li.
\newblock \emph{Boundedness for a {{3D}} Chemotaxis-{{Stokes}} System with
  Porous Medium Diffusion and Tensor-Valued Chemotactic Sensitivity}.
\newblock Z. Angew. Math. Phys., {\bfseries 68} (2017), \penalty0 Art. 29, 23.

\bibitem[Winkler(2011)]{Winkler_MMAS11}{https://doi.org/10.1002/mma.1346}
M.~Winkler.
\newblock \emph{Global Solutions in a Fully Parabolic Chemotaxis System with
  Singular Sensitivity}.
\newblock Math. Methods Appl. Sci., {\bfseries 34} (2011), \penalty0 176--190.

\bibitem[Winkler(2015{\natexlab{a}})]{Win-ct_fluid_3d-CPDE15}{https://doi.org/10.1007/s00526-015-0922-2}
M.~Winkler.
\newblock \emph{Boundedness and Large Time Behavior in a Three-Dimensional
  Chemotaxis-{{Stokes}} System with Nonlinear Diffusion and General
  Sensitivity}.
\newblock Calc. Var. Partial Differential Equations, {\bfseries 54} (2015),
  \penalty0 3789--3828.

\bibitem[Winkler(2015{\natexlab{b}})]{win15_chemorot}{https://doi.org/10.1137/140979708}
M.~Winkler.
\newblock \emph{Large-Data Global Generalized Solutions in a Chemotaxis System
  with Tensor-Valued Sensitivities}.
\newblock SIAM J. Math. Anal., {\bfseries 47} (2015), \penalty0 3092--3115.

\bibitem[Winkler(2016)]{Win16CS1}{https://doi.org/10.1142/S0218202516500238}
M.~Winkler.
\newblock \emph{The Two-Dimensional {{Keller-Segel}} System with Singular
  Sensitivity and Signal Absorption: Global Large-Data Solutions and Their
  Relaxation Properties}.
\newblock Math. Models Methods Appl. Sci., {\bfseries 26} (2016), \penalty0
  987--1024.

\bibitem[Xu et~al.(2020)Xu, Ji, Mei, and
  Yin]{xuChemotaxisModelDegenerate2020}{https://doi.org/10.1016/j.jde.2019.08.013}
T.~Xu, S.~Ji, M.~Mei, and J.~Yin.
\newblock \emph{On a Chemotaxis Model with Degenerate Diffusion: {{Initial}}
  Shrinking, Eventual Smoothness and Expanding}.
\newblock J. Differential Equations, {\bfseries 268} (2020), \penalty0
  414--446.

\end{thebibliography}
\end{document}